\documentclass[a4paper,oneside]{amsbook}

\usepackage{lscape}
\usepackage{multirow}
\usepackage[figuresright]{rotating}
\usepackage{amsmath}
\usepackage{amsthm}
\usepackage{amsfonts}
\usepackage{amssymb}
\usepackage[polish,english]{babel}
\usepackage{epsfig}
\usepackage{ifthen}
\usepackage{fancyhdr}
\usepackage{amscd}

\newcommand{\headfont}[1]{\scriptsize{\MakeUppercase{#1}}}

\ifthenelse{\boolean{@twoside}}{
\renewcommand{\chaptermark}[1]{\markboth{\headfont{\thechapter.\
#1}}{\headfont{\thechapter.\ #1}}}

\fancyhead[LE]{\scriptsize{\thepage}} \fancyhead[CE]{\leftmark}
\fancyhead[CO]{\iffloatpage{\leftmark}{\rightmark}}
\fancyhead[RO]{\scriptsize{\thepage}}}{  
\renewcommand{\chaptermark}[1]{\markboth{\headfont{\thechapter.\ #1}}{}}
\fancyhead[CO]{\leftmark}\fancyhead[RO]{\scriptsize{\thepage}}}

\allowdisplaybreaks 

\newtheorem*{eqp}{The equivalence problem}
\newtheorem*{gp}{The geometry problem}
\newtheorem*{theo}{Theorem}

\newtheorem{theorem}{Theorem}[chapter]
\newtheorem{lemma}[theorem]{Lemma}
\newtheorem{corollary}[theorem]{Corollary}
\newtheorem{proposition}[theorem]{Proposition}

\theoremstyle{definition}
\newtheorem{definition}[theorem]{Definition}
\newtheorem{example}[theorem]{Example}
\newtheorem{remark}[theorem]{Remark}

\DeclareMathOperator{\sgn}{sgn}
\DeclareMathOperator{\Ad}{Ad}
\DeclareMathOperator{\ad}{ad}
\DeclareMathOperator{\Hom}{Hom}

\newcommand{\be}{\begin{equation}}
\newcommand{\ee}{\end{equation}}
\newcommand{\ben}{\begin{equation*}}
\newcommand{\een}{\end{equation*}}
\newcommand{\bal}{\begin{aligned}}
\newcommand{\eal}{\end{aligned}}
\newcommand{\bma}{\begin{pmatrix}}
\newcommand{\ema}{\end{pmatrix}}

\newcommand{\inv}[2]{\mathbf{#1}_{\mathbf{#2}}}
\newcommand{\inw}[1]{\mathbf{#1}}

\newcommand{\wt}[1]{\widetilde{#1}}
\newcommand{\wh}[1]{\widehat{#1}}
\newcommand{\goth}[1]{\mathfrak{#1}}
\newcommand{\cc}[1]{\bar{#1}}

\newcommand{\inc}[2]{\mathbf{#1}^c_{\mathbf{#2}}}
\newcommand{\hc}[1]{\theta^{#1}}
\newcommand{\vc}[1]{\Omega_{#1}}

\newcommand{\inp}[2]{\mathbf{#1}^p_{\mathbf{#2}}}
\newcommand{\hp}[1]{\theta^{#1}}
\newcommand{\vp}[1]{\Omega_{#1}}

\renewcommand{\inf}[2]{\mathbf{#1}^f_{\mathbf{#2}}}
\newcommand{\hf}[1]{\theta^{#1}}
\newcommand{\vf}[1]{\Omega_{#1}}

\newcommand{\der}{{\rm d}}
\newcommand{\w}{{\scriptstyle\wedge}\,}
\newcommand{\hook}{\lrcorner}
\newcommand{\semi}[1]{\oplus_{#1}}

\newcommand{\C}{\mathcal{C}}
\newcommand{\co}{\goth{co}}
\newcommand{\conf}{\goth{conf}}
\newcommand{\D}{\mathcal{D}}
\newcommand{\Der}{\goth{D}}

\newcommand{\g}{\goth{g}}
\renewcommand{\gg}{\mathbf{g}}
\newcommand{\gl}{\goth{gl}}
\renewcommand{\H}{\mathcal{H}}
\newcommand{\h}{\goth{h}}
\newcommand{\J}{{\mathcal J}}
\newcommand{\M}{{\mathcal M}}
\newcommand{\m}{\goth{m}}
\renewcommand{\O}{\mathcal{O}}
\renewcommand{\o}{\goth{o}}
\renewcommand{\P}{{\mathcal P}}
\newcommand{\R}{{\rm R}}
\newcommand{\real}{\mathbb{R}}
\newcommand{\Ric}{{\rm Ric}}
\renewcommand{\S}{\mathcal{S}}
\newcommand{\so}{\goth{so}}
\renewcommand{\sp}{\goth{sp}}
\newcommand{\T}{{\bf T}}
\renewcommand{\u}{\goth{u}}
\newcommand{\V}{{\mathcal V}}

\begin{document}
\frontmatter
\begin{titlepage}

\begin{center}
 {\Large \textsc{ University of Warsaw} \\
  \textsc{Faculty of Physics}\\
 \textsc{Institute of Theoretical Physics}}

  \vspace*{\stretch{0.67}}
  {\bf \LARGE  Geometry of Third-Order
  Ordinary Differential Equations and Its Applications \smallskip \\ in General Relativity}

  \vspace*{\stretch{1}}
  {\LARGE Micha\l\ Godli\'nski}

  \vspace*{\stretch{1}}
  \includegraphics[width = .4\linewidth, height = .4\linewidth, angle = 0]{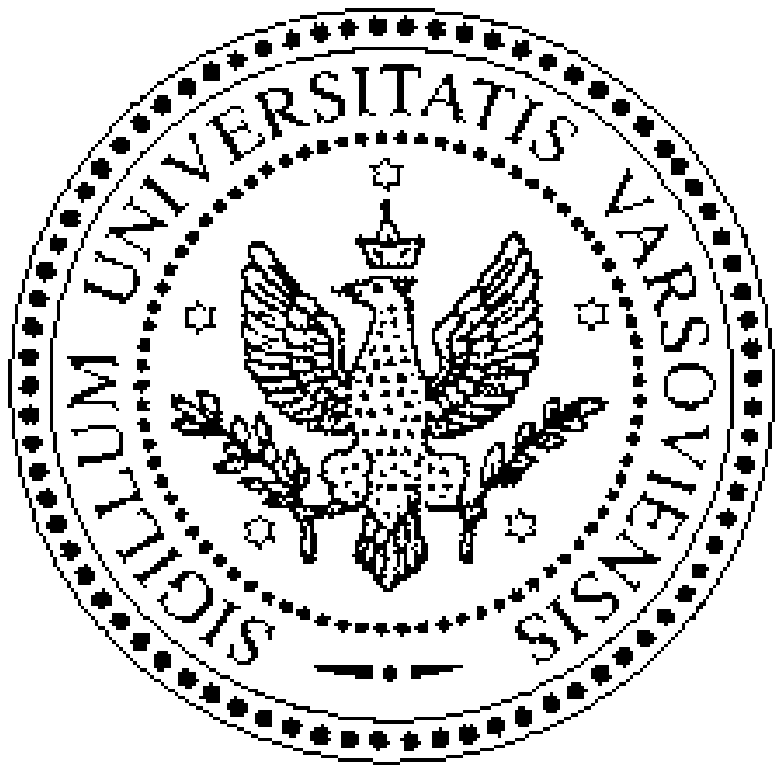} \\
   \vspace*{\stretch{1}}
  { \Large PhD thesis written \\ at
 the Chair of Theory of Relativity and Gravitation \\
  under supervision of \\
  \medskip
  \LARGE Dr. hab. Pawe\l\ Nurowski} \\
   \vspace*{\stretch{2}}
  {\Large \textsc{Warsaw  2008}}
\end{center}
\end{titlepage}

\setcounter{page}{2} \ifthenelse{\boolean{@twoside}}{%
\thispagestyle{empty}
$\phantom{empty}$\\
\clearpage}{}

\thispagestyle{empty}
\selectlanguage{polish}

\begin{center} \vspace*{\stretch{1}}
{\large \it  Dzi"ekuj"e dr. hab. Paw"lowi Nurowskiemu za opiek"e
naukow"a, za temat tej rozprawy i liczne wskaz"owki, bez kt"orych
nigdy by ona nie powsta"la.
\\ \bigskip

Dzi"ekuj"e mojej Rodzinie za mi"lo"s"c i wsparcie. Szczeg"olnie
dzi"ekuj"e Rodzicom, Siostrom i Ma"lej R"aczce. \\ \bigskip

Dzi"ekuj"e Bogu.}

\vspace*{\stretch{2}}
\end{center}

\selectlanguage{english}

 \tableofcontents

\mainmatter

\chapter{Introduction}\label{ch.problem}
 \noindent This thesis addresses the problems of equivalence
and geometry of third order ordinary differential equations (ODEs)
which are stated as follows.
\begin{eqp} Given two real differential equations
\begin{equation}
y'''=F(x,y,y',y'')\label{e10}
\end{equation}
and \be
 y'''=\bar{F}(x,y,y',y'')\label{e20}
\ee for a real function $y=y(x)$, establish whether or not there
exists a local transformation of variables of a suitable type that
transforms \eqref{e10} into \eqref{e20}. \end{eqp}
\begin{gp}
Determine geometric structures defined by a class of equations
$y'''=F(x,y,y',y'')$ equivalent under certain type of
transformations. Find relations between invariants of the ODEs and
invariants of the geometric structures.
\end{gp}

One may consider equivalence with respect to several types of
transformations, in this work we focus on three best known types:
contact, point and fibre-preserving transformations. The
fibre-preserving transformations are those which transform the
independent variable $x$ and the dependent variable $y$ in such a
way that the notion of the independent variable is retained, that
is the transformation of $x$ is a function of $x$ only:
\be\label{e.fp}
   x\mapsto\bar{x}=\chi(x),\quad\quad\quad\quad  y\mapsto\bar{y}=\phi(x,y).
\ee The transformation rules for the derivatives are already
uniquely defined by above formulae. Let us define the total
derivative to be \ben
  \Der=\partial_x+y'\partial_y+y''\partial_{y'}+y'''\partial_{y''}.
\een Then
\begin{subequations}\label{e.tprol}
\begin{align}
 y'&\mapsto\frac{\der\bar{y}}{\der\bar{x}}=\frac{\Der \phi}{\Der \chi}, \label{e.tprol1}  \\
 y''&\mapsto\frac{\der^2\bar{y}}{\der\bar{x}^2}= \frac{\Der}{\Der\chi}\left(\frac{\Der\phi}{\Der \chi}\right),\label{e.tprol2} \\
 y'''&\mapsto \frac{\der^3\bar{y}}{\der\bar{x}^3}=\frac{\Der}{\Der\chi}\left(\frac{\Der}{\Der\chi}\left(\frac{\Der\phi}{\Der \chi}\right)\right).\label{e.tprol3}
\end{align}
\end{subequations}

The point transformations of variables mix $x$ and $y$ in an
arbitrary way \be \label{e.point}
   x\mapsto\bar{x}=\chi(x,y),\quad\quad\quad\quad  y\mapsto\bar{y}=\phi(x,y),
\ee with the derivatives transforming as in \eqref{e.tprol}.

The contact transformations are more general yet. Not only they
augment the independent and the dependent variables but also the
first derivative \be\label{e.cont}
\begin{aligned}
   x&\mapsto\bar{x}=\chi(x,y,y'),\notag\\
   y&\mapsto\bar{y}=\phi(x,y,y'),\\
   y'&\mapsto\frac{\der \bar{y}}{\der\bar{x}}=\psi(x,y,y').\notag
\end{aligned}
\ee However, the functions $\chi$, $\phi$ and $\psi$ are not
arbitrary here but subjecting to \eqref{e.tprol1} which now yields
two additional constraints \ben
\psi=\frac{\Der\phi}{\Der\chi}\quad \iff \quad
\bal &\psi\chi_{y'}=\phi_{y'}, \\
  &\psi(\chi_x+y'\chi_y)=&\phi_x+y'\phi_y, \eal
\een guaranteeing that $\der\bar{y}/\der\bar{x}$ really transforms
like first derivative. With these conditions fulfilled second and
third derivative transform through \eqref{e.tprol2} --
\eqref{e.tprol3}. Of course, fibre-preserving transformations are
a subclass of point ones, as well as point transformations form a
subclass within contact ones.

We always assume in this work that ODEs are defined locally by a
smooth real function $F$ and are considered apart from
singularities. The transformations are always assumed to be local
diffeomorphisms.

\begin{example} In order to illustrate the equivalence problem let us consider
whether or not the equations
$$ y'''=0 \quad\quad\text{and}\quad\quad y'''=3\frac{y''^2}{y'}$$
are equivalent. As one can easily check they are contact
equivalent, since the transformation \ben
\begin{aligned}
   x&\mapsto-2y',\notag\\
   y&\mapsto2xy'^2-2yy',\\
   y'&\mapsto-2xy'+y \notag
\end{aligned} \een applied to $y'''=0$ brings
it to the other equation. However, they are not fibre-preserving
equivalent, since the quantity
$$I(x,y,y',y'')=\frac{\partial^2}{\partial y''^{\,2}}F(x,y,y',y'')$$ vanishes for every
equation which is fibre-preserving equivalent to $y'''=0$ but does
not vanish for $y'''=3\frac{y''^2}{y'}$. In order to see this we
apply a fibre-preserving transformation of general form
\eqref{e.point} to $y'''=0$ and check that $I=0$ for the resulting
equation.
\end{example}

The above example shows  importance of relative invariants in the
equivalence problem. A relative invariant is a function of $F$ and
its derivatives such that if it vanishes for an equation
$y'''=F(x,y,y',y'')$ then it also vanishes for every equation
equivalent to it, thereby each relative invariant provides us with
a necessary condition for equivalence. Moreover, with the help of
adequate number of relative invariants we can also formulate
sufficient conditions for equivalence of ODEs, although this
construction is complicated and can hardly ever be carried out to
the very end because of difficult calculations. The second
problem, the problem of geometry, is even more fundamental and in
fact it contains the problem of equivalence, for once a geometry
associated with ODEs is constructed and a relationship between
local invariants of the geometry and of ODEs is found then one can
study the equivalence of ODEs via objects of the associated
geometry.

A pioneering work on geometry of ODEs of arbitrary order is Karl
W\"unschmann's PhD thesis \cite{Wun} written under supervision of
F. Engel in 1905. In this paper K. W\"unschmann observed that
solutions of an $n$th-order ODE
$y^{(n)}=F(x,y,y',\ldots,y^{(n-1)})$ may be considered as both
curves $y=y(x,c_0,c_1,\ldots,c_{n-1})$ in the $xy$ space and
points $c=(c_0,\ldots,c_{n-1})$ in the solution space $\real^n$
parameterized by values of the constants of integration $c_i$. He
defined a relation of $k$th-order contact between infinitesimally
close solutions considered as curves; two solutions $y(x)$ and
$y(x)+\der y(x)$ corresponding to $c$ and $c+\der c$ have the
$k$th-order contact if their $k$th jets coincide at some point
$(x_0, y_0)$. W\"unschmann's main question was how the property of
having $(n-2)$nd contact for $n=3,4$ and $5$ might be described in
terms of the solution space. In particular he examined third-order
ODEs and showed that there is a distinguished class of ODEs
satisfying certain condition for the function $F$, which we call
the W\"unschmann condition. For a third-order ODE in this class,
the condition of having first order contact is described by a
second order Monge equation for $\der c$.  This Monge equation is
nothing but the condition that the vector defined by two
infinitesimally close points $c$ and $c+\der c$ is null with
respect to a Lorentzian conformal metric on the solution space.
The last observation, although not contained in W\"unschmann's
work, follows immediately from his reasoning and was later made by
S.-S. Chern \cite{Chern}, who cited W\"unschmann's thesis.

 The main contribution to the issue of point and contact
geometry of third-order ODEs was made by respectively E. Cartan
and S.-S. Chern in their classical papers \cite{Car1, Car2, Car3}
and \cite{Chern}. We discuss their approach and results in section
\ref{s.i.sum} of Introduction; here we only mention that E. Cartan
\cite{Car2} proved that every third-order ODE modulo point
transformations and satisfying two differential conditions on the
function $F$, one of them being the W\"unschmann condition, has a
three-dimensional Lorentzian Einstein-Weyl geometry on its
solution space. E. Cartan also showed how to construct invariants
of this Weyl geometry from relative point invariants of the
underlying ODE. In the same vein S.-S. Chern constructed a
three-dimensional Lorentzian conformal geometry for third-order
ODEs considered modulo contact transformations and satisfying the
W\"unschmann condition. In both these cases the conformal metric
is precisely the metric appearing implicitly in K. W\"unschmann's
thesis. Later H. Sato and A. Yoshikawa \cite{Sat} applying N.
Tanaka's theory \cite{Tan} constructed a Cartan normal connection
for arbitrary third-order ODEs (not only of the W\"unschmann type)
and showed how its curvature is expressed by the contact relative
invariants.

Geometry of third-order ODEs appears in General Relativity and the
theory of integrable systems. E.T. Newman et al
\cite{nsf1,nsf2,nsf3} devised the Null Surface Formulation (NSF),
a description of General Relativity in terms of families of null
hypersurface generated by a Lorentzian metric. The
$2+1$-dimensional version of this formalism \cite{nsf3d} is
equivalent to Chern's conformal geometry for third-order ODEs
\cite{New1, New2, New3}, which was noticed by P. Tod \cite{Tod}
for the first time. We encapsulate results of NSF in section
\ref{s.i.nsf}. Three-dimensional Einstein-Weyl geometry was
studied mainly from the perspective of the theory of twistors and
integrable systems by N. Hitchin \cite{Hitch}, R. Ward
\cite{Ward}, C. LeBrun \cite{Leb} and P. Tod, M. Dunajski et al
\cite{Jt,Dun1,Dun2}; for discussion of the link between the
Einstein-Weyl spaces and third-order ODEs see \cite{Nur1} and
\cite{Tod}.

 P. Nurowski, following the ideas of E.
Cartan, proposed a programme of systematic study of geometries
related to differential equations, including second- and
third-order ODEs. In this programme,
\cite{New3,God3,Godode5,Nur1,Nur2,Nur3}, both new and already
known geometries associated with differential equations are
supposed to be constructed by the Cartan equivalence method and
are to be characterized in the language of Cartan connections
associated with them. In particular in \cite{Nur1} P. Nurowski
provided new examples of geometries associated with ordinary
differential equations including a conformal geometry with special
holonomy $G_2$ from ODEs of the Monge type. Partial results on
geometries of third-order ODEs were given in \cite{Nur1, Nur2,
New3, God3} but the full analysis of these geometries has not been
published so far and this thesis, which is a part of the
programme, aims to fill this gap.

Geometry of third-order ODEs is a part of broader issue of
geometry of differential equations in general. Regarding ODEs of
order two, we owe classical results including construction of
point invariants to S. Lie \cite{Lie} and M. Tresse \cite{Tre}. In
particular E. Cartan \cite{Car5} constructed a two-dimensional
projective differential geometry on the solution spaces of some
second-order ODEs. This geometry was further studied in \cite{NN}
and \cite{Nur3}, the latter paper pursues the analogy between
geometry of three-dimensional CR structures and second-order ODEs
and provides a construction of counterparts of the Fefferman
metrics for the ODEs. Classification of second-order ODEs
possessing Lie groups of fibre-preserving symmetries was done by
L. Hsu and N. Kamran \cite{Hsu}. Geometry on the solution space of
certain four-order ODEs (satisfying two differential conditions),
which is given by the four-dimensional irreducible representation
of $GL(2,\real)$ and has exotic $GL(2,\real)$ holonomy was
discovered and studied by R. Bryant \cite{Bry}, see also
\cite{Nur4}. The $GL(2,\real)$ geometry of fifth-order ODEs has
been recently studied by M. Godli\'nski and P. Nurowski
\cite{Godode5}.

The more general problems yet are existence and properties of
geometry on solution spaces of arbitrary ODEs. The problem of
existence was solved by B. Doubrov \cite{Dubgl}, who  proved that
an $n$th-order ODE $n\geq 3$, modulo contact transformations, has
 a geometry based on the irreducible
$n$-dimensional representation of $GL(2,\real)$ provided that it
satisfies  $n-2$ scalar differential conditions. An implicit
method of constructing these conditions was given in
\cite{Dubwil}. Properties of the $GL(2,\real)$ geometries of ODEs
are still an open problem; they were studied in \cite{Godode5},
where the Doubrov conditions were interpreted as higher order
counterparts of the W\"unschmann condition, and by M. Dunajski and
P. Tod \cite{Dun3}. 

Almost all the above papers deal with geometries on solution
spaces but one can also consider other geometries, including those
defined on various jet spaces. The most general result on such
geometries \cite{Dub2} comes from T. Morimoto's nilpotent geometry
\cite{Mor1,Mor2}. It concludes that with a system of ODEs there is
associated a filtration on a suitable jet space together with a
canonical Cartan connection.

\section{Geometry of third-order ODEs --- the present status} \label{s.i.sum}

\noindent We recapitulate the classical results on geometry of
third-order ODEs. We do this mostly in the original spirit of E.
Cartan, emphasizing the role of systems of one-forms and Cartan
connections. The version of Cartan's equivalence method
\cite{Car4} we employ below is explained in books by R. Gardner
\cite{Gar} and P. Olver \cite{Olv2}. Some its aspects are also
discussed by S. Sternberg \cite{Ste} and S. Kobayashi \cite{Kob1}.

\subsection*{E. Cartan's and S.-S. Chern's approach to ODEs}
They began with the space $\J^2$ of second jets of curves in
$\real^2$ (see \cite{Olv1} and \cite{Olv2} for extensive
description of jet spaces) with coordinate system $(x,y,p,q)$
where $p$ and $q$ denote the first and second derivative $y'$ and
$y''$ for a curve $x\mapsto(x,y(x))$ in $\real^2$, so that this
curve lifts to a curve $x\mapsto(x,y(x),y'(x),y''(x))$ in $\J^2$.
Any solution $y=f(x)$ of $y'''=F(x,y,y',y'')$ is uniquely defined
by a choice of $f(x_0)$, $f'(x_0)$ and $f''(x_0)$ at some $x_0$.
Since that choice is equivalent to a choice of a point in $\J^2$,
there passes exactly one solution $(x,f(x),f'(x),f''(x))$ through
any point of $\J^2$. Therefore the solutions form a (local)
congruence on $\J^2$, which can be described by its annihilating
simple ideal. Let us choose a coframe $(\omega^i)$ on $\J^2$:
\be\label{e.omega}\begin{aligned}
 \omega^1&=\der y -p\der x,  \\
 \omega^2&=\der p- q\der x,  \\
 \omega^3&=\der q- F(x,y,p,q)\der x,  \\
 \omega^4&=\der x.
\end{aligned}
\ee Each solution $y=f(x)$  is fully described by the two
conditions: forms $\omega^1$, $\omega^2$, $\omega^3$ vanish on the
curve $t\mapsto(t,f(t),f'(t),f''(t))$ and, since this defines a
solution modulo transformations of $x$, $\omega^4=\der t$ on this
curve.

Suppose now that  equation \eqref{e10} undergoes a contact, point
or fibre-preserving transformation.\footnote{Although E.Cartan and
S.-S. Chern did not examine fibre-preserving transformations we
treat them on equal footing with the others for future
references.} Then \eqref{e.omega} transform by \be\bal
 \omega^1&\mapsto \bar{\omega}^1=u_1\omega^1, \\
 \omega^2&\mapsto \bar{\omega}^2=u_2\omega^1+u_3\omega^2,  \\
 \omega^3&\mapsto \bar{\omega}^3=u_4\omega^1+u_5\omega^2+u_6\omega^3,\\
 \omega^4&\mapsto \bar{\omega}^4=u_8\omega^1+u_9\omega^2+u_7\omega^4, \\
\eal\label{e.trans_kor} \ee with some functions $u_1,\ldots,u_9$
defined on $\J^2$ and determined by a particular choice of
transformation, for instance \begin{eqnarray*}
&u_1=\phi_y-\psi\chi_y,\\
&u_7=\Der \chi,\quad u_8=\chi_y,  \quad u_9= \chi_p.
\end{eqnarray*}
In particular, $u_9=0$ in the point case and $u_8=u_9=0$ in the
fibre-preserving case. Since the transformations are
non-degenerate, the condition $u_1u_3u_5u_7\neq 0$ is always
satisfied and transformations \eqref{e.trans_kor} form groups.
Thus a class of contact equivalent third-order ODEs is a local
G-structure\footnote{We work in the local trivializations of the
G-structure} $G_c\times\J^2$, that is a local subbundle of the
bundle of linear frames on $\J^2$, defined by the property that
the coframe $(\omega^1,\omega^2,\omega^3,\omega^4)$ belongs to it
and the structural group is given by \be\label{e.G_c}
G_c= \bma u_1 & 0 & 0 & 0 \\ u_2 & u_3 & 0 & 0 \\
u_4 & u_5 & u_6& 0\\
u_8 & u_9 & 0 & u_7 \ema. \ee In a like manner, for point and
fibre-preserving transformation we have $G$-structures with groups
$G_p\subset G_c$ given by $u_9=0$ and $G_f\subset G_p$ given by
$u_8=u_9=0$.

Thereby the problem of equivalence has been formulated; two ODEs
are contact/point/fibre-preserving equivalent if and only if the
$G_c$/$G_p$/$G_f$-structures which correspond to them are
equivalent. In other words, it means that there exists a
diffeomorphism which transforms one $G$-structure into the other
$G$-structure.

On the bundles $G_c\times\J^2$, $G_p\times\J^2$ or $G_f\times\J^2$
there are four fixed, well defined one-forms
$(\theta^1,\theta^2,\theta^3,\theta^4)$, the components of the
canonical $\real^4$-valued form $\theta$ existing on the frame
bundle of $\J^2$. Let us consider the contact case as an example.
Let $(x)$ denote $x,y,p,q$ and $(g)$ be coordinates in $G_c$ given
by \eqref{e.G_c}. Let us choose a coordinate system $(x,g)$ on
$G_c\times\J^2$ compatible with the local trivialization. Then
$\theta^i$ at the point $(x,g^{-1})$ read \be
\label{e.i.canon}\bal
 \theta^1=&u_1\omega^1, \\
 \theta^2=&u_2\omega^1+u_3\omega^2,  \\
 \theta^3=&u_4\omega^1+u_5\omega^2+u_6\omega^3,\\
 \theta^4=&u_8\omega^1+u_9\omega^2+u_7\omega^4. \\
\eal \ee

The idea of Cartan's method is the following. Starting from
$G_c\times\J^2$ (or $G_p\times\J^2$ or $G_f\times\J^2$) one
constructs a new principal bundle $\P=H\times\J^2$ (with a new
group $H$) equipped with one fixed coframe
$(\theta^i,\Omega_{\mu})$ built in some geometric way and such
that it encodes all the local invariant information about
$G_c\times\J^2$ (or $G_p\times\J^2$ or $G_f\times\J^2$
respectively) through its structural equations.

\begin{remark} As an illustration of this idea consider its very special
application to a metric structure of signature $(k,l)$ on a
manifold $\M$. The structure defines locally the bundle
$O(k,l)\times\M$, an $O(k,l)$-reduction of the frame bundle. If
$(\omega^i)$ is an orthonormal coframe on $\M$, then the forms
$(\theta^i)$ at a point $(x,g)$ are given by
$\theta^i=(g^{-1})^i_{~j}\omega^j$, $g\in O(k,l)$. Moreover,
associated to $(\omega^i)$ there are the Levi-Civita connection
one-forms $\Gamma^i_{~j}$, which when lifted to $O(k,l)\times\M$
define $\tfrac12 n(n-1)$ connection forms $\Omega^i_{~j}$ by \ben
\Omega^i_{~j}(x,g)=(g^{-1})^i_{~k}\Gamma^k_{~l}(x)\,g^l_{~j}+(g^{-1})^i_{~k}\der
g^k_{~j}. \een Clearly $(\theta^i,\Omega^j_{~k})$ is a coframe on
$O(k,l)\times\M$. It is the Cartan coframe for the metric
structure. The Cartan invariants are the Riemannian curvature $R$
given by the equations
\ben \bal &\der\theta^i+\Omega^i_{~j}\w\theta^j=0, \\
&\der\Omega^i_{~j}+\Omega^j_{~k}\w\Omega^k_{~j}=\tfrac12R^i_{jkl}
\theta^k\w\theta^l \eal \een and its consecutive covariant
derivatives. Here it is easy to construct the desired coframe
because there is a distinguished $\o(k,l)$-valued connection ---
the Levi-Civita connection. \end{remark}

In our cases situation is more complicated. There are no natural
candidates for the forms $\Omega_{\mu}$ since there is no
distinguished connection on the bundles $G_{\cdot}\times J^2$. In
order to resolve such situations E. Cartan developed the method of
constructing the desired bundle $\P$ through a series of
reductions and prolongations of an initial structural group. This
procedure is too complicated to be included in Introduction; we
postpone the full discussion to section \ref{s.c.proof} of chapter
\ref{ch.contact} where we follow E. Cartan's and S.-S. Chern's
reasoning. Here we only formulate main conclusions.

\subsection*{Point equivalence} Examining the
problem of point equivalence E. Cartan constructed for every ODE a
local seven-dimensional bundle $\P$ over $J^2$, a reduction of
$G_p\times\J^2$, together with a coframe
$(\theta^1,\theta^2,\theta^3,\theta^4,\Omega_1,\Omega_2,\Omega_3)$
on $\P$. The coframe contains all the point invariant information
on the original ODE due to the following theorem, which is the
application of the equivalence method to third-order ODEs.

\begin{theorem}[E. Cartan]\label{th.i.clas} Equations $y'''=F(x,y,y',y'')$ and
$\cc{y}'''=\cc{F}(\cc{x},\cc{y},\cc{y}',\cc{y}'')$ with smooth $F$
and $\cc{F}$ are locally point equivalent if and only if there
exists a local diffeomorphism $\Phi\colon\P\to\cc{\P}$ which maps
the Cartan coframe $(\cc{\theta}^1,\ldots,\cc{\Omega}_3)$ of
$\cc{F}$ to the coframe $(\theta^1,\ldots,\Omega^3)$ of $F$, \ben
\Phi^*\cc{\theta}^i=\theta^i,\qquad
\Phi^*\cc{\Omega}_\mu=\Omega_\mu,\qquad i=1,\ldots,4,\quad
\mu=1,2,3. \een $\Phi$ projects to a point transformation
$(x,y)\mapsto(\cc{x}(x,y),\cc{y}(x,y))$ which maps $\cc{F}$ to
$F$.
\end{theorem}

The general idea of determining whether or not $\Phi$ exist is as
follows, see\cite{Olv2}. The structural equations for the Cartan
coframe are
\begin{align}
 \der\hp{1} =&\vp{1}\w\hp{1}+\hp{4}\w\hp{2},\nonumber \\
 \der\hp{2} =&\vp{2}\w\hp{1}+\vp{3}\w\hp{2}+\hp{4}\w\hp{3},\nonumber \\
 \der\hp{3}=&\vp{2}\w\hp{2}+(2\vp{3}-\vp{1})\w\hp{3}+\inv{A}{1}\hp{4}\w\hp{1},\nonumber \\
 \der\hp{4} =&(\vp{1}-\vp{3})\w\hp{4}+\inv{B}{1}\hp{2}\w\hp{1}
   +\inv{B}{2}\hp{3}\w\hp{1}, \nonumber\\
 \der\vp{1} =&-\vp{2}\w\hp{4}+(\inv{D}{1}+3\inv{B}{3})\hp{1}\w\hp{2}
   +(3\inv{B}{4}-2\inv{B}{1})\hp{1}\w\hp{3} \label{e.i.dtheta_7d}  \\
   &+(2\inv{C}{1}-\inv{A}{2})\hp{1}\w\hp{4}-\inv{B}{2}\hp{2}\w\hp{3}, \nonumber \\
 \der\vp{2}=&(\vp{3}-\vp{1})\w\vp{2}+\inv{D}{2}\hp{1}\w\hp{2}+(\inv{D}{1}+\inv{B}{3})\hp{1}\w\hp{3}
 +\inv{A}{3}\hp{1}\w\hp{4} \nonumber \\
  &+(2\inv{B}{4}-\inv{B}{1})\hp{2}\w\hp{3}+\inv{C}{1}\hp{2}\w\hp{4}, \nonumber \\
 \der\vp{3}=&(\inv{D}{1}+2\inv{B}{3})\hp{1}\w\hp{2}+2(\inv{B}{4}-\inv{B}{1})\hp{1}\w\hp{3}
  +\inv{C}{1}\hp{1}\w\hp{4}+\inv{B}{2}\hp{2}\w\hp{3}, \nonumber
\end{align}
with some explicitly given functions
$\inv{A}{1},\inv{A}{2},\inv{A}{3},\inv{B}{1},\inv{B}{2},\inv{B}{3},
\inv{B}{4},\inv{C}{1},\inv{D}{1},\inv{D}{2}$ on $\P$. Let
$(X_1,\ldots,X_7)$ denotes the frame dual to
$(\theta^1,\ldots,\Omega_3)$. Since pull-back commutes with
exterior differentiation each of these functions is a point
relative invariant of the underlying ODE. Furthermore, the coframe
derivatives $X_1(\inv{A}{1}),\ldots,X_7(\inv{D}{2})$,
$X_1(X_1(\inv{A}{1})),\ldots,X_7(X_7(\inv{D}{2})),\ldots$ of
arbitrary order are also point relative invariants. We calculate
all relevant coframe derivatives for $F$ and $\cc{F}$ up to some
finite order $n$ and gather them into respective functions
$\T\colon\P\to\real^N$ and $\cc{\T}\colon\cc{\P}\to\real^N$ with
the same target space of dimension $N$ equal to the number of the
invariants. Finally, we examine whether or not the graphs of $\T$
and $\cc{\T}$ overlap as manifolds in $\real^N$. If they overlap,
then the sought diffeomorphism $\Phi$ exists between certain
nonempty open sets $U\subset\T^{-1}(\O)$ and
$\cc{U}\subset\cc{\T}^{-1}(\O)$, where $\O$ is the overlap. A
detailed explanation of this procedure is given at the beginning
of chapter \ref{ch.class}; at this stage we only need to know that
the coframe solves the point equivalence problem for ODEs. What
joins this problem with the domain of differential geometry is the
fact that the coframe is a geometric object
--- a Cartan connection
--- on $\P\to\J^2$.

In order to introduce the notion of Cartan connection we first
show how to read the structure of principal bundle on $\P\to\J^2$
from \eqref{e.i.dtheta_7d}. Fibres of the projection $\P\to J^2$
are annihilated by the forms
$\theta^1,\theta^2,\theta^3,\theta^4$, and the vector fields
$X_5,X_6,X_7$ are tangent to the fibres. Simultaneously, the
commutation relations of these fields are isomorphic to
commutators of the three-dimensional algebra\footnote{The symbol
$\g\semi{.}\h$ denotes a semidirect product of the Lie algebras
$\g$ and $\h$} $\h=\real\oplus(\real\semi{.}\real)$ and they
define a local action of the group
$H=\real\times(\real\ltimes\real)$ on $\P$ for which $X_5,X_6,X_7$
are the fundamental fields.

Now, let us focus on the most symmetric situation when all the
relative invariants $\inv{A}{1},\ldots,\inv{D}{2}$ vanish. This
case corresponds to an ODE which is point equivalent to the
trivial $y'''=0$. For such an ODE the equations
\eqref{e.i.dtheta_7d} become the Maurer-Cartan equations for the
algebra $\co(2,1)\semi{.}\real^3$, where
$\co(2,1)=\real\oplus\so(2,1)$ is the orthogonal algebra centrally
extended by dilatations generated by the identity matrix. As a
consequence, $\P$ becomes locally the Lie group
$CO(2,1)\ltimes\real^3$, with the left-invariant fields
$X_1,\ldots,X_7$. Therefore $\J^2$ is the homogeneous space $H\to
CO(2,1)\ltimes\real^3 \to CO(2,1)\ltimes\real^3/H$ acted upon by
$CO(2,1)\ltimes\real^3$. This group acts on $\J^2$ as the group of
point symmetries of $y'''=0$, that is for any solution $y=f(x)$ of
$y'''=0$ its graph $x\to(x,y(x),y'(x),y''(x))$ in $\J^2$ is
transformed into the graph of other solution of $y'''=0$. The
coframe $(\theta^i,\Omega_\mu)$ can be arranged into the following
matrix \be\label{e.i.conp}
 \wh{\omega}=\bma \vp{3} & 0 & 0 & 0 & 0 \\
             \hp{1} & \vp{3}-\vp{1} & -\hp{4} & 0 & 0\\
             \hp{2} & -\vp{2} & 0 & -\hp{4} & 0 \\
             \hp{3} & 0 &-\vp{2} & \vp{1}-\vp{3} & 0 \\
             0 & \hp{3} & -\hp{2} & \hp{1} & -\vp{3}
        \ema,
\ee which is the Maurer-Cartan one-form  of
$CO(2,1)\ltimes\real^3$. In the language of $\wh{\omega}$ the
equations \eqref{e.i.dtheta_7d} read \ben
\der\wh{\omega}+\tfrac12[\wh{\omega},\wh{\omega}]=0. \een

Turning to an arbitrary situation, we see that non-trivial cases
can not be described in this language, since the invariants
$\inv{A}{1},\ldots,\inv{D}{2}$ do not vanish in general and the
last equation does not hold. We need a new object, the Cartan
connection, defined here after \cite{Kob1}.
\begin{definition}
Let $H\to\P\to\M$ be a principal bundle and let $G$ be a Lie group
such that $H$ is its closed subgroup and $\dim G=\dim\P$. A Cartan
connection of type $(G,H)$ on $\P$ is a one-form $\wh{\omega}$
taking values in the Lie algebra $\goth{g}$ of $G$ and satisfying
the following conditions:
\begin{itemize}
\item[i)] $\wh{\omega}_u:T_u\P\to\goth{g}$ for every $u\in\P$ is
an isomorphism of vector spaces \item[ii)]
$A^*\hook\,\wh{\omega}=A$ for every $A\in\goth{h}$ and the
corresponding fundamental field $A^*$ \item[iii)]
$R^*_h\wh{\omega}=\Ad(h^{-1})\wh{\omega}$ for $h\in H$.
\end{itemize}
\end{definition}
A Cartan connection is then an object that generalizes the notion
of the Maurer-Cartan form on a Lie group. The curvature of a
Cartan connection \ben
\wh{K}=\der\wh{\omega}+\frac{1}{2}[\wh{\omega},\wh{\omega}] \een
does not have to vanish any longer and measures, as it were, how
far $H\to\P\to\M$ `differs' from the homogeneous space $H\to G\to
G/H$.

In our case the formula \eqref{e.i.conp} defines a
$\co(2,1)\semi{.}\real^3$-valued Cartan connection, whose
non-vanishing curvature contains basic point invariants, from
which the full set of invariants can be constructed through
exterior differentiation. This is the basic relation, announced at
the beginning of Introduction, between classification of ODEs and
their geometry. In chapter \ref{ch.point} we discuss properties of
$\wh{\omega}$ in full detail.

The geometry introduced above, however, is not the most
interesting structure one can associate with the ODEs modulo point
transformations. Indeed, the crucial observation E. Cartan made in
his paper \cite{Car2} was that the connection $\wh{\omega}$ may
generate a new type of geometry on the solution space of certain
classes of ODEs.

In order to present the idea of this construction let us invoke
once more the trivial equation $y'''=0$ but now consider how the
symmetry group act on the solutions regarded as  curves in the
$xy$ plane. It is well-known that the full group of point
symmetries of $y'''=0$ is generated by the following one-parameter
groups of transformations $(x,\,y)\mapsto\Phi^i_t(x,\,y)$
\begin{align*}
 &\Phi^1_t(x,\,y)=(x,\,y+t), &  & \Phi^2_t(x,\,y)=(x,\,y+2xt),\\
 &\Phi^3_t(x,\,y)=(x,\,y+x^2t), &  & \Phi^4_t(x,\,y)=(x,\,ye^{-t}),\\
 &\Phi^5_t(x,\,y)=(xe^{t},\,ye^{t}), &  & \Phi^6_t(x,\,y)=(x+t,\,y),\\
 & \Phi^7_t(x,y)=\left(\frac{x}{1+xt},\,\frac{y}{(1+xt)^2}\right). &  &
\end{align*}
A solution to $y'''=0$ is a parabola of the form \be\label{solo}
y(x)=c_2x^2+2c_1x+c_0\ee with three arbitrary integration
constants $c_2,c_1,c_0$. Therefore a solution of $y'''=0$ may be
identified with a point $c=(c_2,c_1,c_0)^T$ in the solution space
$\S\cong\real^3$. Transforming \eqref{solo} according to the above
formulae we find that $\Phi^1_t,\Phi^2_t,\Phi^3_t$ are
translations in the solution space $\S$: \ben \Phi^1_t(c)=\bma c_2
\\ c_1 \\ c_0-t \ema, \quad \Phi^2_t(c)=\bma
c_2 \\ c_1-t \\ c_0\ema, \quad \Phi^3_t(c)= \bma c_2-t \\ c_1 \\
c_0 \ema, \een while transformations $\Phi^4_t,\ldots,\Phi^7_t$
generate the three-dimensional irreducible representation of
$CO(2,1)$:
\begin{align*}
&\Phi^4_t(c)=\exp t\bma 1 & 0 & 0 \\ 0 & 1 & 0 \\ 0 & 0 & 1 \ema
c, &
 &\Phi^5_t(c)=\exp t\bma 1 & 0 & 0 \\ 0 & 0 & 0 \\ 0 & 0 & -1 \ema c,
 \\\\
 &\Phi^6_t(c)=\exp t\bma 0 & 0 & 0 \\ 1 & 0 & 0 \\ 0 & 2 & 0 \ema c, &
 &\Phi^7_t(c)=\exp t\bma 0 & 2 & 0 \\ 0 & 0 & 1 \\ 0 & 0 & 0 \ema c.
\end{align*}
Thereby $\S$ is a homogeneous space equipped with the flat
conformal metric\footnote{The symbol $[g]$ denotes the conformal
metric, whose representative is the metric $g$. We also adapt the
following notation:
$\alpha\beta=\tfrac12(\alpha\otimes\beta+\beta\otimes\alpha)$ for
one-forms $\alpha$ and $\beta$.} $[g]$ \be\label{e.i.gflat}
g=(\der c_0)(\der c_2)-(\der c_1)^2 \ee and preserved by $CO(2,1)$
\ben M^TgM=e^{\lambda} g,\quad \text{for}\quad M\in CO(2,1)
\quad\text{and}\quad \lambda\in\real. \een

Apart from $[g]$ there is another piece of structure in $\S$. The
symmetry group acting on $\S$ is not the full ten-dimensional
conformal symmetry group of the flat metric $Conf(2,1)\cong
O(3,2)$, but merely $CO(2,1)\ltimes\real^3$, the Euclidean group
extended by dilatations, hence another object is needed to reduce
$O(3,2)$ to $CO(2,1)\ltimes\real^3$. This object is the Weyl
one-form of the flat Weyl geometry. Let us remind that a Weyl
geometry $(g,\phi)$ is a metric $g$ and a one-form $\phi$ given
modulo transformations $g\to e^{2\lambda} g$,
$\phi\to\phi+\der\lambda$, see chapter \ref{ch.point} section
\ref{s.ew.def}. In our case the Weyl geometry is flat, because
there is the flat representative \eqref{e.i.gflat} for which, in
addition, $\phi=0$.

Now a question arises: how the flat Weyl structure on $\S$ can be
reconstructed from the Cartan coframe $(\theta^i,\Omega_\mu)$? In
order to answer this question we will use the method of
construction through fibre bundles and Lie transport, which
differs from E. Cartan's original reasoning and was introduced by
P. Nurowski, cf \cite{Nur0, Nur1}.

To begin with, we observe that $\P$ is always, not only in the
trivial case, a principal bundle $CO(2,1)\to\P\to\S$. Indeed, as
we said, $\J^2$ is foliated by solutions of the underlying ODE,
which are curves in $\J^2$ annihilated by $\omega^1,\omega^2$ and
$\omega^3$. Thus we have a projection $\J^2\to\S$ with the fibres
being the solutions. As a consequence $\P$ is also a bundle over
$\S$ with leaves of the projection annihilated by the ideal
$(\theta^1,\theta^2,\theta^3)$. On the leaves of $\P\to\S$ there
act the vector fields $X_4,\ldots,X_7$, members of the frame dual
to $(\theta^i,\Omega_\mu)$. In view of \eqref{e.i.dtheta_7d} their
commutation relations are isomorphic to $\co(2,1)$ and turn each
leaf into an orbit of the free action of $CO(2,1)$ regardless
whether or not the invariants $\inv{A}{1},\ldots,\inv{D}{2}$
vanish. If in addition $\inv{A}{1}=\ldots=\inv{D}{2}=0$ then
locally $\P\cong CO(2,1)\ltimes\real^3$, which generates the
structure of homogeneous space on $\S$. Next, still in the trivial
case, let us consider the bilinear form \be\label{i.gP}
\wh{g}=2\theta^1\theta^3-(\theta^2)^2\ee and the one-form
$\Omega_3$ on $\P$. We calculate the Lie derivatives of these
objects along the fields $X_4,\ldots,X_7$ tangent to the fibres of
$\P\to\S$ and observe that

\be\label{e.i.lieg}\bal
&\wh{g}(X_j,\cdot)\equiv0,\quad\text{for}\quad j=4,5,6,7, \\
&L_{X_4}\wh{g}=0,\quad L_{X_5}\wh{g}=0,\quad L_{X_6}\wh{g}=0,\quad
L_{X_7}\wh{g}=2\wh{g}\eal \ee and \be\label{e.i.lienu} \bal
&X_4\hook\Omega_3=0,\quad X_5\hook\Omega_3=0,\quad
X_6\hook\Omega_3=0,\quad
X_7\hook\Omega_3=1,\\
&L_{X_j}\Omega_3=0,\quad\text{for}\quad j=4,5,6,7. \eal\ee
 These properties allow
us to project the pair $(\wh{g},\Omega_3)$ along $\P\to\S$ to a
Weyl structure $(g,\phi)$ on $\S$. This is precisely the flat Weyl
structure obtained above by action of the symmetry group.

The construction through Lie derivatives and the projection has an
essential advantage in comparison to the symmetry approach, since
it can be immediately generalized to non-trivial equations. In
fact, the equations \eqref{e.i.lieg}, \eqref{e.i.lienu} hold not
only in the trivial case $y'''=0$ but under much weaker conditions
$\inv{A}{1}=0$ and $\inv{C}{1}=0$. Calculating the explicit forms
of these invariants Cartan proved that if only an ODE given by a
function $F(x,y,p,q)$ satisfies
\begin{align} &\left(\D-\tfrac{2}{3}F_q\right)\left(
\tfrac{1}{6}\D F_q
-\tfrac{1}{9}F_q^2-\tfrac{1}{2}F_p\right)+F_y=0, \label{e.i.wunsch} \\ \notag\\
&\D^2F_{qq}-\D F_{qp}+F_{qy}=0,\label{e.i.cart} \\
\intertext{where}
&\D=\partial_x+p\partial_y+q\partial_p+F\partial_q,
\notag\end{align} then it has a Weyl geometry on its solution
space $\S$. The quantity on the left hand side of the condition
\eqref{e.i.wunsch}, the W\"unschmann invariant, was found in
\cite{Wun} for the first time.

 For the equations satisfying
\eqref{e.i.wunsch} and \eqref{e.i.cart} the bundle
$CO(2,1)\to\P\to\S$ is the bundle of orthonormal frames for the
conformal metric $[g]$ on $\S$, whereas $\wh{\omega}$ of
\eqref{e.i.conp} becomes a $\co(2,1)\semi{.}\real^3$-valued Cartan
connection. The $\co(2,1)$-part of $\wh{\omega}$ is the Weyl
connection, and the Weyl curvature is expressed by the  invariants
$\inv{B}{1},\inv{B}{2},\inv{B}{3},\inv{B}{4}$, which do not vanish
in general. Cartan proved that all the Weyl geometries constructed
of third-order ODEs are Einstein, that is the traceless part of
the Ricci tensor for their Weyl connections always vanishes.
Finally, he observed \cite{Car3} that these Einstein-Weyl
geometries are of general form; for each three-dimensional
Lorentzian Einstein-Weyl structure there is an ODE satisfying
\eqref{e.i.wunsch}, \eqref{e.i.cart}, whose solution space carries
this structure.

\subsection*{Contact equivalence}
Soon after E.Cartan studied the point equivalence and geometry,
the contact problem was examined by S.-S. Chern using the same
method. After a slightly more complicated construction (we discuss
it in chapter \ref{ch.contact}) Chern built a ten-dimensional
bundle $\P\to\J^2$ equipped with the coframe
$(\theta^1,\theta^2,\theta^3,\theta^4,\Omega_1,\ldots,\Omega_6)$.
For $y'''=0$ the structural equations coincide with the
Maurer-Cartan equations for the algebra $\o(3,2)\cong\sp(4,\real)$
which reflects the fact that $O(3,2)$ is the maximal group of
contact symmetries of $y'''=0$.

Next difference between Cartan's and Chern's  results is that here
the coframe $(\theta^i,\Omega_\mu)$ for a general ODE is not a
Cartan connection since it does not transform regularly along the
fibres of the bundle $\P\to\J^2$, i.e. does not fulfill condition
iii) of the definition. Lack of this regularity means that there
are nonconstant $\Omega_\alpha\w\theta^k$ terms in the structural
equations. However, from the point of view of S.-S. Chern, who was
interested in geometries on the solutions space analogous to the
Einstein-Weyl geometries, this fault was insignificant and he did
not investigate it. The Cartan connection for this problem was
given later in \cite{Sat} and in \cite{Dub2}.

Leaving aside this topic let us consider S.-S. Chern's results.
Given the ten-dimensional coframe he noticed that in the
structural equations there appears the W\"unschmann invariant and
further analysis depends on whether it vanishes or not. If the
W\"unschmann invariant vanishes then the coframe
$(\theta^i,\Omega_\mu)$ becomes an $\o(3,2)$-valued Cartan normal
conformal connection over $\S$. In this case S.-S. Chern gave an
explicit formula for the underlying Lorentzian conformal metric
$[g]$ in terms of the forms $\omega^1,\omega^2,\omega^3$ on
$\J^2$, see eq. \eqref{e.c.g} in chapter \ref{ch.contact}. Since
the contact symmetry group of $y'''=0$ is $O(3,2)$, the full
conformal group of the flat metric in $2+1$ dimensions, there is
no additional geometric object on $\S$ and we have just the
conformal geometry instead of a Weyl geometry there. Later it was
shown \cite{Nur2} that the conformal geometry can be obtained from
the bilinear symmetric field
$\wh{g}=(\theta^2)^2-2\theta^1\theta^3$ on $\P$ by virtue of
conditions analogous to \eqref{e.i.lieg}. The lowest-order
conformal invariant in three dimensions --- the Cotton tensor ---
was also expressed in terms of contact relative invariants for
ODEs.

Turning to the equations with non-vanishing W\"unschmann
invariant, S.-S. Chern continued reduction of the bundle $\P$
according to Cartan's method and obtained a five-dimensional
manifold, say $\P_5$, which is a line bundle over $\J^2$ and is
furnished with a coframe $(\theta^1,\ldots,\theta^4,\Omega)$.
Next, he recognized  that in the homogeneous case of the equation
$y'''=-y$ the bundle $\P_5$ is a group $\real^2\ltimes\real^3$ and
it generates a geometry of cones on the solution space. This is
rather an exotic kind of structure and we will not discuss it
here, referring the reader to section \ref{s.c.furred} of chapter
\ref{ch.contact}.

This closes our summary of the known results on the geometry and
classification of third-order ODEs.

\section{Null Surface Formulation}\label{s.i.nsf}
\noindent An intriguing thing about third-order ODEs is that the
Lorentzian geometry on $\S$ was rediscovered fifty years after
S.-S. Chern from completely different perspective, the perspective
of General Relativity. In a series of papers
\cite{nsf1,nsf2,nsf3,nsf4} E.T. Newman et al proposed and
developed the Null Surface Formulation (NSF), an alternate
approach to General Relativity. Ideas of this approach are the
following. Let $\M^4$ be a four-manifold coordinated by $(x^\mu)$.
Usually, the basic object in General Relativity is a Lorentzian
metric $g$ on $\M^4$ (or the Levi-Civita connection), from which
other objects are derived, the Riemann tensor, the Weyl tensor and
the Einstein tensor, on which dynamics is imposed by the Einstein
equations. One of many objects in Lorentzian geometry is a null
hypersurface. It is by definition a hypersurface $ Z(x^\mu)=const$
in $\M^4$ which satisfies the eiconal equation \be\label{i.eikon}
g(\der Z, \der Z)=0,\ee where $g$ is the inverse (contravariant)
metric. For a given $g$ the family of all null hypersurfaces is
fully defined by the above equation.

The NSF is an alternate point of view, here one begins with $\M^4$
without any metric but endowed with a two-parameter family of
hypersurfaces. Starting from these data a Lorentzian {\em
conformal} metric is constructed by the property that these
hypersurfaces are its null hypersurfaces; the metric is found by
solving the eiconal equation \eqref{i.eikon} with respect to the
components $g^{\mu\nu}$. In this approach it is the family of
hypersurfaces that is a basic concept and the metric is a derived
one. The family is defined by a sufficiently differentiable real
function $Z(x^\mu,s,s^*)$ on $\M\times S^2$, where $s$ and $s^*$
are stereographic variables on the sphere $S^2$ and $^*$ denotes
the complex conjugation. Each hypersurface is a level set
$$Z(x^\mu,s,s^*)=const,$$
with some fixed $s$. The reason why $s\in S^2$ is that the
hypersurfaces are to be null and then $S^2$ becomes the sphere of
null directions. When $x^\mu$ are fixed and $s$ sweeps out the
sphere we obtain the corresponding hypersurfaces at $x^\mu$
orthogonal to all null directions. In order to find the metric the
authors assumed that the functions \be\label{i.fi} f^0=Z,\quad
f^+=Z_s,\quad f^-=Z_{s^*},\quad f^1=Z_{ss^*} \ee form a coordinate
system on $\M^4$ for all $s$, introduced the coframe $(\der
f^0,\der f^+,\der f^-,\der f^1)$ and found explicit formulae for
components $ g^{\mu\nu}$ in this coframe. The eiconal equation
implies $g^{00}=0$ in the coframe $(\der f^A)$ but one must still
take into account vanishing of derivatives of $g^{00}$ with
respect to $s$ and $s^*$. Doing so the authors found that the
conformal metric is uniquely defined by the family of
hypersurfaces provided that the function $Z$ and a real conformal
factor $\Omega^2(x^\mu,s,s^*)$ for the Lorentzian conformal metric
satisfy two quite complicated complex differential conditions,
referred to as metricity conditions. If they are satisfied then
the components of the conformal Lorentzian metric $[g]$ are
products of $\Omega^2$ and some expressions containing derivatives
of $Z$. Moreover, both metricity conditions and the components
$g^{AB}$ only depend of $Z$ via its derivatives $Z_{ss}$ and
$Z_{s^*s^*}$ and one can eliminate the space-time coordinates in
$\Omega$ and $Z_{ss}$ through \eqref{i.fi} and obtain the
functions \ben\Lambda(f^A,s,s^*)=Z_{ss}(x^\mu(f^A,s,s^*),s,s^*),
\qquad \Omega(f^A,s,s^*)=\Omega(x^\mu(f^A,s,s^*),s,s^*).\een
 Thereby information about a Lorentzian metric is encoded by two
 complex functions
$\Omega$ and $\Lambda$ of the variables $(s, s^*, Z(s,s^*), Z_s,
Z_{s^*}, Z_{ss^*})$, satisfying the two metricity conditions and
the integrability condition $\frac{\der^2}{\der
s^{*2}}\Lambda=\frac{\der^2}{\der s^{2}}\Lambda^*$.

Next step was writing down Einstein equations
$G^{\mu\nu}=\varkappa T^{\mu\nu}$ in the new variables. The
authors proved that applying consecutive derivatives $\partial_s$
and $\partial_{s^*}$ to \be\label{i.ein}
G^{\mu\nu}Z_{\mu}Z_{\nu}=\varkappa T^{\mu\nu}Z_{\mu}Z_{\nu},\ee
$Z_\mu=\partial_\mu Z$, one obtains nine out of ten equations and
the lacking tenth equation, the trace component, is recovered with
the help of the metricity conditions. After suitable
substitutions, the vacuum version of \eqref{i.ein} reduces to one
equation
$$\Omega_{f^1f^1}-Q[\Lambda]\Omega=0.$$
In this manner the vacuum Einstein equations were reduced to the
set of four complex equations for $\Omega$ and $\Lambda$: the one
above, the two metricity conditions and the integrability
condition for $\Lambda$. Having proven this result in \cite{nsf3},
the authors moved to the analysis of solutions of the Einstein
equations, their linearization, perspectives for quantization and
other topics we will not cover here. From our point of view the
most interesting is the three-dimensional version of the NSF
\cite{nsf3d,New1, New2, New3}, which leads immediately to
third-order ODEs. We cite the construction in detail.

In the case of $2+1$-dimensional Lorentzian geometry on $\M^3$ the
space of null directions is diffeomorphic to $S^1$ and a conformal
class can be reconstructed from the one-parameter family of
surfaces
$$ Z(x^i,s)=const, $$ where $(x^i)=(x^0,x^1,x^2)\in\M^3$ and $s\in
S^1$ real. Let us introduce the functions
\be\label{i.ypq}y=Z(x^i,s),\quad p=Z_s(x^i,s),\quad
q=Z_{ss}(x^i,s)\ee and the coframe
$$\sigma^0=\der y,\quad \sigma^1=\der p, \quad \sigma^2=\der q. $$
The third derivative, $Z_{sss}(x^i,s)$, together with
\eqref{i.ypq} define a real function
$$\Lambda(s,y,p,q)=Z_{sss}(x^i(s,y,p,q),s).$$
It is obvious that the variables $(s,\,y,\,p,\,q)$ constitute a
coordinate system on $\J^2$, the space of second jets of functions
$S^1\to\real$. Moreover, the function $\Lambda$ on $\J^2$ defines
the third-order ODE \be\label{i.3ord}
Z_{sss}=\Lambda(s,Z,Z_s,Z_{ss})\ee for $Z(s)$. The function
$Z(x^i,s)$, with which we have begun, can be identified with the
general solution of this equation, $x^i$ playing the role of three
integration constants. In this manner $\M^3$ becomes the solution
space of the ODE and $\M^3\times S^1$ is identified with $\J^2$,
where the projection $\J^2\to\M^3$ is given by
$(x^i,s)\mapsto(x^i)$. The fibres of this projection are solutions
of \eqref{i.3ord} considered as curves in $J^2$. We also notice
that the total derivative $\frac{\der}{\der s}$ applied to a
function on $\J^2$ coincides with the total derivative $\D$
$$ \frac{\der}{\der s}f(s,y,p,q)=\D f= (\partial_s+p\partial_y
+q\partial_p+\Lambda\partial_q)f. $$
It follows that the construction of a conformal geometry from the
family of surfaces is fully equivalent to Chern's construction
described earlier and the metricity conditions contain the
W\"unschmann condition. Let us look at how this construction was
done. The eiconal equation
$$g^{ij}Z_{i}Z_{j}=g^{ij}y_{i}y_{j}=0$$ for the sought metric
$g=g^{ij}\partial_{x^i}\otimes\partial_{x^j}$ implies, as before,
$g^{00}=0$ in the coframe $(\sigma^i)$. Taking $\partial_s g^{00}$
gives $$ g^{ij}y_{i}p_{j}=g^{01}=0. $$ Another derivation
$\partial_s$ yields
$$ g^{ij}y_{i}q_{j}+g^{ij}p_{i}p_{j}=g^{02}+g^{11}=0. $$
Third and fourth derivatives yield \ben g^{02}\Lambda_q+3g^{12}=0
\een and \ben
g^{02}(\D\Lambda_q-\tfrac13\Lambda_q^2-3\Lambda_p)+3g^{22}=0. \een
At this point all the components of the metric are found and are
proportional to $g^{02}$, which becomes naturally the conformal
factor $\Omega^2(x^i,s)$. The conformal metric in the basis
$(\partial_y,\partial_p,\partial_q)$ is as follows \ben g^{ij}=\Omega^2\bma 0 & 0 & 1 \\
0 & -1 & -\tfrac13\Lambda_q
\\ 1 & -\tfrac13\Lambda_q &
-\tfrac13\D\Lambda_q+\tfrac19\Lambda_q^2+\Lambda_p \ema. \een
$\Omega$ is not an arbitrary function of $s$, since
$\Omega^2=g^{02}=g^{ij}y_iq_i$ and applying $\D$ we get \ben \D
\Omega=\frac13\Omega\Lambda_q.\een This is first of the metricity
conditions, it is a constraint on $\Omega$ provided that $\Lambda$
is known. Second condition is obtained by taking the fifth
derivative of the eiconal equation with respect to $s$. It reads
\begin{align*} 0&=g^{ij}(5p_i
\partial^4_s
Z_j+Z_i\partial^5_sZ_j+10\partial^2Z_i\partial^3Z_j)=\\
&=5g^{11}(\D\Lambda)_p+5g^{12}(\D\Lambda)_q+g^{02}(\D^2\Lambda)_q
+10(g^{02}\Lambda_y+g^{12}\Lambda_p+g^{22}\Lambda_q).
\end{align*}
Substituting the metric components by their explicit forms we
obtain precisely the W\"unschmann condition \eqref{e.i.wunsch} for
$\Lambda$.

The relations between third-order ODEs and $2+1$-dimensional
conformal geometry in the NSF were further studied in \cite{New1,
New2, New3}. In \cite{New1} it was shown that one can generate the
conformal metric on the solution space starting from the system
\eqref{e.omega} of one-forms $\omega^i$ on $\J^2$ associated with
an ODE. The authors defined the tensor
\be\label{i.newg}2\omega^1(\omega^3+a\omega^1+b\omega^2)-(\omega^2)^2\ee
on $\J^2$ with arbitrary functions $a$ and $b$. Next they imposed
the condition that Lie transport of the above tensor along fibres
of the projection $\J^2\to solution\,space$ is conformal. The
construction is parallel to \eqref{i.gP} and the conditions for
conformal Lie transport fix uniquely $a$ and $b$, and yield the
W\"unschmann condition in a like manner to \eqref{e.i.lieg} and
\eqref{e.i.lienu}. In \cite{New2} the authors re-proved the
invariance of the so obtained conformal geometry under contact
transformations of the related ODE, which reflects a gauge freedom
in the NSF. Finally, in \cite{New3} explicit formulae were given
for the curvature of normal conformal connection in terms of
contact invariants of third-order ODE.

Simultaneously to the research on the three-dimensional version of
the NSF, much progress was made in understanding the full
four-dimensional formalism in \cite{New2, New3, New4}. Here $(s,
s^*, Z(s,s^*), Z_s, Z_{s^*}, Z_{ss^*}, Z_{ss}, Z_{s*s*})$ are
coordinates on the second jet space $\J^2(S^2,\real)$ of functions
$S^2\to\real$. Owing to this fact, the complex function $\Lambda$
defines a pair of PDEs for a real function $Z$ of the variables
$s$ and $s^*$ through
\begin{align}
Z_{ss}=&\Lambda(s, s^*, Z, Z_s, Z_{s^*}, Z_{ss^*}), \notag
\\
Z_{s^*s^*}=&\Lambda^*(s, s^*, Z, Z_s, Z_{s^*}, Z_{ss^*}), \notag
\end{align}
while the space-time $\M^4$ is the space of solutions to this
system, $(x^\mu)$ being constants of integrations. Following ideas
of the three-dimensional construction the authors considered the
six-dimensional submanifold $L$ of $\J^2(\real^2,\real)$ given by
the PDEs. They considered the Pfaffian system on $L$ associated
with the PDEs and, following the formula \eqref{i.newg}, they
built a symmetric tensor on $L$ from the Pfaff forms. This tensor
projects to a conformal geometry on $\M^4$ provided that the
metricity conditions for $\Lambda$ are satisfied, which can be now
viewed as generalizations of the W\"unschmann condition. The
Cartan normal conformal connection for this geometry was
constructed \cite{New4}. The Null Surface Formulation is still an
ongoing project and work is being continued to obtain conformal
Einstein equations in terms of invariants of the PDEs.

\section{Results of the thesis}
\noindent  The thesis has threefold aim.
\begin{itemize}
 \item[i)] Constructing new geometries associated with third-order ODEs
modulo contact, point and fibre-preserving transformations of
variables. Considering possible applications to General
Relativity.

\item[ii)] Describing new geometries together with the already
known in the language of Cartan connections with curvature given
by respective invariants of ODEs obtained by Cartan's method.

\item[iii)] Applying Cartan invariants to classification of
certain types of third-order ODEs.
\end{itemize}

In order to construct the geometries we follow Cartan's
construction outlined in section \ref{s.i.sum} of Introduction. We
start from the equivalence problems formulated in terms of the
$G$-structures \eqref{e.i.canon} and apply Cartan's method to
obtain manifolds $\P$ equipped with the coframes encoding all the
invariant information about ODEs. Next we show how to read the
principal bundle structures of these manifolds over distinct
bases. Usually it is the structure over $\S$ which is the most
interesting, but we also consider structures over $\J^1$, $\J^2$
and certain six-dimensional manifold $\M^6$, which appears
naturally. When the structure of a bundle is established then the
invariant coframe defines a Cartan connection on $\P\to base$,
usually under additional conditions playing similar role to the
W\"unschmann condition. In order to obtain the geometries and the
W\"unschmann-like conditions we often apply the P. Nurowski method
of construction by Lie transport and projection.

In the most symmetric cases, when the underlying ODEs have
transitive symmetry groups, the bundles $\P$ become locally Lie
groups while their bases become homogeneous spaces, providing
homogeneous models for the geometries. Since the dimension of $\S$
and $\J^1$ is three and we build geometries with at least
two-dimensional structural group then the homogeneous models are
given by the ODEs with at least five-dimensional symmetry group.
Below we encapsulate our results. The geometries of the ODEs are
also gathered in table \ref{t.geom} on page \pageref{t.geom}.

Unfortunately, among fourteen geometries which we consider in this
work, there is no Lorentzian geometry or any new geometry which
could currently be applied to General Relativity.

\subsection*{Contact geometries}
Chapter \ref{ch.contact} is devoted to the geometries of the ODEs
modulo contact transformations. The only equations possessing at
least five- dimensional contact symmetry group are the linear
equations with constant coefficients\footnote{We prove this
statement in chapter \ref{ch.class}.}, that is \ben y'''=0 \een
with the symmetry group $O(3,2)$ and \ben y'''=-2\mu y'+y,\qquad
\mu\in \real, \een mutually non-equivalent for distinct $\mu$,
with the symmetry group $\real^2\ltimes\real^3$. Sections
\ref{s.c.th} to \ref{s.c6d} discuss geometries whose homogeneous
model is generated by $y'''=0$. In section \ref{s.c.th} we state
the main theorem in that chapter, theorem \ref{th.c.1}, which
describes the geometry on $\J^2$. It can be recapitulated as
follows.

\begin{theo}[Theorem \ref{th.c.1}] The contact invariant information
about an equation $y'''=F(x,y,y',y'')$ is given by the following
data
\begin{itemize}
\item[i)] The principal fibre bundle $H_6\to\P\to\J^2$, where
$\dim\P=10$ and $H_6$ is a six-dimensional subgroup of $O(3,2)$
\item[ii)] The coframe
$(\theta^1,\theta^2,\theta^3,\theta^4,\Omega_1,\Omega_2,\Omega_3,\Omega_4,\Omega_5,\Omega_6)$
on $\P$ which defines  the $\o(3,2)\cong\sp(4,\real)$ Cartan
normal connection $\wh{\omega}^c$ on $\P$.
\end{itemize}
The coframe and the connection $\wh{\omega}_c$ are given
explicitly in terms of $F$ and its derivatives. There are two
basic relative invariants for this geometry: the W\"unschmann
invariant and $F_{y''y''y''y''}$.\end{theo} This theorem is almost
identical to the result proved in \cite{Sat} and the only new
element we add here is the explicit formula for the connection.
Section \ref{s.c.proof} contains the proof of the theorem, which
repeats S.-S. Chern's construction of the coframe and the
construction of the normal connection of \cite{Sat}.

Sections \ref{s.c.conf} to \ref{s.c6d} discuss next three
geometries generated by $\wh{\omega}^c$. Section \ref{s.c.conf}
contains the construction of the Lorentzian geometry on the
solution space and, after \cite{Nur2}, gives explicit formulae for
the normal conformal connection. Section \ref{s.c-p} studies a new
type of geometry in the context of third-order ODEs, the
contact-projective structure on $\J^1$. The idea of this geometry
is the following. Consider the solutions of an ODE as a family of
curves in $\J^1$ and ask whether these curves are among geodesics
of a linear connection. The answer to this question is positive
provided that $F_{y''y''y''y''}=0$ and in this case there is a
whole family of connections for which the solutions are geodesics.
Such a family of connections in $\J^1$ is an example of a contact
projective structure, see D. Fox \cite{Fox} for the general
definition of contact projective structures. Moreover, Tanaka's
theory allows us to define a notion of normal Cartan connection
for these structures in dimension three \cite{Cap, Cap2, Fox}. It
is then a matter of straightforward calculations to check that
$\wh{\omega}^c$ is the normal connection for our contact
projective structure. In section \ref{s.c6d} we construct a
six-dimensional split signature conformal geometry on some
six-manifold $\M^6$ over which $\P$ is a bundle. Next we show that
the associated normal conformal connection $\wh{\mathbf{w}}$ has a
special conformal holonomy reduced from $\o(4,4)$ to
$\o(3,2)\semi{.}\real^5\subset\o(4,4)$ and, what is more,
$\wh{\omega}^c$ is precisely the $\o(3,2)$ part of
$\wh{\mathbf{w}}$. These results are summarized as follows.
\begin{theo}
The connection $\wh{\omega}^c$ of theorem \ref{th.c.1} has
fourfold interpretation.
\begin{itemize}
 \item[1.] It is always the normal $\o(3,2)$ Cartan connection  on
 $\J^2$ (Chern-Sato-Yoshikawa construction.)
\item[2.] If the ODE has vanishing W\"unschmann condition then
$\wh{\omega}^c$ is the normal Lorentzian conformal connection  for
the Lorentzian structure on the solution space (Chern-NSF
construction.) \item[3.] If the ODE satisfies $F_{y''y''y''y''}=0$
then $\wh{\omega}^c$ becomes the normal Cartan connection for the
contact projective structure on $\J^1$. \item[4.] $\wh{\omega}^c$
is the $\o(3,2)$-part of the $\o(4,4)$ normal conformal connection
for the six-dimensional split conformal geometry on $\M^6$ with
special holonomy $\o(3,2)\semi{.}\real^5$.
\end{itemize}
\end{theo}

In section \ref{s.c.furred} we turn to geometries, whose
homogeneous models are provided by the equation $y'''=-2\mu y'+y$.
Following Chern we reduce the bundle $\P$ to its five-dimensional
subbundle. Then we find that
\begin{theo}[Theorem \ref{cor.c.geom5d}]
Every ODE satisfying some contact invariant condition
$\inc{a}{}[F]=\mu=const$ has a $\real^2$ geometry on its solution
space together with a $\real^2$ linear connection from the
invariant coframe. The action of the algebra $\real^2$ on $\S$ is
given by \ben
   \bma
      u & v & 0 \\
      -\mu v & u & v \\
      v & -\mu v & u
   \ema.
 \een
\end{theo}
This geometry seems to be a generalization of Chern's `cone
geometry' which was associated with the equation $y'''=-y$ and
briefly mentioned to exist for arbitrary ODEs. In our construction
the action of $\real^2$ depends on the characteristic polynomial
of respective linear equation, and we get a real cone geometry
provided that it has three distinct roots.

\subsection*{Point geometries} In sections \ref{s.p.th} to \ref{s.lor3}
of  chapter \ref{ch.point} we study the geometries associated with
the ODEs modulo point transformations. Sections \ref{s.p.th} to
\ref{s.w6d} deal with the geometries modelled on $y'''=0$. Our
approach is analogous to the contact case and results are similar,
with some caveats. The algebra $\co(2,1)\semi{.}\real^3$ is not
semisimple, hence the methods of the Tanka theory fail here.
Moreover, even its generalization, the Morimoto nilpotent
geometry, does not work in this case, since the point geometry of
third-order ODEs on $\J^2$ is not a filtration. As a consequence,
we do not have a general theory about existence of Cartan
connections and in the case of $\J^1$, where we find point
projective structure, certain refinement of contact projective
structure, we are not able to find any connection and suppose that
it does not exist in general. We also construct a six-dimensional
Weyl structure in the split signature, which is related to the
six-dimensional conformal geometry of the contact case. To
summarize
\begin{theo}
The following statements hold
\begin{itemize}
\item[1.] The point invariant information about
$y'''=F(x,y,y',y'')$ is given by the seven-dimensional principal
bundle $H_3\to\P\to\J^2$ together with the coframe
$\hp{1},\hp{2},\hp{3},\hp{4},\vp{1},\vp{2},\vp{3}$ on $\P$, which
defines the $\co(2,1)\semi{.}\real^3$ Cartan connection
$\wh{\omega}^p$ (Cartan construction.) \item[2.] If the ODE has
vanishing W\"unschmann \eqref{e.i.wunsch} and Cartan
\eqref{e.i.cart} invariants then it has the Einstein-Weyl geometry
on $\S$ and the Weyl connection is given by $\wh{\omega}^p$
(Cartan construction.) \item[3.] If the ODE satisfies
$F_{y''y''y''}=0$ then it has the point-projective structure on
$\J^1$. \item[4.] For any ODE there exists the split signature
six-dimensional Weyl geometry, which is never Einstein.
\end{itemize}
\end{theo}

A new construction, which does not have a contact counterpart, is
considered in section \ref{s.lor3}. This is a Lorentzian {\em
metric structure} on the solution space $\S$. Its construction
follows immediately from the Einstein-Weyl geometry. If the Ricci
scalar of the Weyl connection is non-zero, then it is a weighted
conformal function and may be fixed to a constant by an
appropriate choice of the conformal gauge. The homogeneous models
of this geometry are associated with
$$y'''=\frac{3\,y''^2}{2\,y'}$$
if the Ricci scalar is negative, and
$$y'''=\frac{3y''^2y'}{y'^2+1}$$ if the Ricci scalar is positive.
Both these equations are contact equivalent to $y'''=0$ and their
point symmetry groups are $O(2,2)$ and $O(4)$ respectively.

\subsection*{Fibre-preserving geometries} Sections \ref{s.f}
and \ref{s.fp} of chapter \ref{ch.point} are devoted to the
geometries of ODEs modulo fibre-preserving transformations. We
obtain a seven-dimensional bundle and the
$\co(2,1)\semi{.}\real^3$-valued Cartan connection $\wh{\omega}^f$
on it. Since both $\wh{\omega}^f$ and $\wh{\omega}^p$ of the point
case take value in the same algebra these cases are very similar
to each other. Indeed, we show that one can recover
$\wh{\omega}^f$ from $\wh{\omega}^p$ just by appending one
function on the bundle. As a consequence, the geometries of the
fibre-preserving case are obtained from their point counterparts
by appending the object generated by the function.

We did not study obvious or not interesting geometries. In the
point and fibre-preserving case geometries of $y'''=-2\mu y'+y$
are the same to what we have in the contact case, since the
respective symmetry groups are the same. Also the fibre-preserving
geometry on $\J^1$ does not seem to be worth studying.

\subsection*{Classification of ODEs}
Chapter \ref{ch.class} contains the classification part of this
work. We obtain two results
\begin{itemize}
\item[i)] We characterize regular ODEs admitting large contact and
point symmetry groups, that is the groups of dimension at least
four. We give the conditions, in terms of Cartan invariants, for
an ODE to posses the large symmetries. The classification is given
in tables \ref{t.cc.1} and \ref{t.pp.1} on pages \pageref{t.cc.1}
-- \pageref{t.pp.2}. We give the criteria for contact
linearization of the ODEs. (The fibre-preserving classification of
the ODEs with large symmetries is again parallel to the point
classification and has been already done \cite{God1, Grebot}, see
also Remark \ref{rem.fp}.) \item[ii)] We characterize regular ODEs
fibre-preserving equivalent to II, IV, V, VI, VII and XI reduced
Chazy classes, which are certain polynomial ODEs with the Painleve
property. We give the explicit formulae for the transformations.
\end{itemize}

The condition of regularity assumed above is of technical nature,
cf the discussion in the beginning of chapter \ref{ch.class}, in
particular definition \ref{def.cc.reg}.

To summarize, the following material contained in this work is
new: sections \ref{s.c-p} to \ref{s.c.furred} of chapter
\ref{ch.contact} excluding theorem \ref{th.c.2}, sections
\ref{s.p-p} to \ref{s.fp} of chapter \ref{ch.point} and the whole
of chapter \ref{ch.class}. Other sections contain a re-formulation
and an extension of already known results.

All our calculations were performed or checked using the symbolic
calculations program Maple.

\section{Notation}
In what follows we use the following symbols, in particular $W$
denotes the W\"unschmann invariant.
\begin{align}
 F=&F(x,y,p,q), \notag \\
\D=&\partial_x+p\partial_y+q\partial_p+F\partial_q, \\ 
 K=&\tfrac{1}{6}\D F_q -\tfrac{1}{9}F_q^2-\tfrac{1}{2}F_p, \label{e.defK} \\
 L=&\tfrac13F_{qq}K-\tfrac13F_qK_q-K_p-\tfrac13F_{qy},  \\ 
 M=&2K_{qq}K-2K_{qy}+\tfrac13F_{qq}L-\tfrac23F_qL_q-2L_p,  \\ 
 W=&\left(\D-\tfrac{2}{3}F_q\right)K+F_y, \label{e.defW} \\
 Z=&\frac{\D W}{W} -F_q.  
\end{align}
Parentheses denote sets of objects: \ben (a_1,\ldots,a_k)\een is
the set consisting of $a_1,\ldots,a_k$. In particular this symbol
denotes bases of vector spaces as well as coordinate systems,
frames and coframes on manifolds.  The linear span of vectors or
covectors $a_1,\ldots,a_k$ is denoted by \ben <a_1,\ldots,a_k>.
\een If $a_1,\ldots,a_k$ are vector fields or one-forms on a
manifold, then the above symbol denotes the distribution or the
simple ideal generated by them. The symmetric tensor product of
two one-forms or vector fields $\alpha$ and $\beta$ is denoted by
\ben
\alpha\beta=\tfrac12(\alpha\otimes\beta+\beta\otimes\alpha).\een
The symbols $A_{(\mu\nu)}$ and $A_{[\mu\nu]}$ denote
symmetrization and antisymmetrization of a tensor $A_{\mu\nu}$
respectively. For a metric $g$ of signature $(k,l)$ the group
$CO(k,l)$ is defined to be
$$CO(k,l)=\{A\in GL(k+l,\real)\,|\quad A^TgA=e^\lambda g,\quad\lambda\in\real \}. $$
Its Lie algebra
$$\co(k,l)=\{a\in \gl(k+l,\real)\,|\quad a^Tg+ga=\lambda g,\quad\lambda\in\real \}. $$
A semidirect product of two Lie groups $G$ $H$, where $G$ acts on
$H$ is denoted by $ G\ltimes H$. A semidirect product of their Lie
algebras is denoted by $\g\semi{.}\h$. If the action depends of a
parameter $\mu$ then we add an subscript: $\ltimes_\mu$ and
$\semi{\mu}$.

\begin{landscape}
\begin{table}
 \caption{Geometries of third-order ODEs}\label{t.geom}
\renewcommand*{\arraystretch}{1.5}
\begin{tabular}{|c|c|c|c|c|}
\hline Model & Manifold & Contact &
Point & Fibre-preserving \\
\hline \hline

\multirow{8}{*}{$y'''=0$} &  \multirow{2}{*}{$\J^2$} &
\multirow{2}{*}{$\o(3,2)$ connection} &
\multicolumn{2}{|c|}{\multirow{2}{*}{$\co(2,1)\semi{.}\real^3$ connection}} \\
& & & \multicolumn{2}{|c|}{} \\ \cline{2-5}

 & \multirow{2}{*}{$\S$} & \multirow{2}{*}{Lorentzian conformal} &
\multirow{2}{*}{Lorentzian Einstein-Weyl} &
\raisebox{-0.25ex}{Lorentzian
Einstein-Weyl} \\

& & & & \raisebox{0.3ex}{with a weighted function} \\\cline{2-5}

 & \multirow{2}{*}{$\J^1$} & \multirow{2}{*}{contact projective} &
\multirow{2}{*}{point projective} & \multirow{2}{*}{{\em was not
studied}} \\

& & & & \\ \cline{2-5}

& \multirow{2}{*}{$\M^6$} & \multirow{2}{*}{split conformal}
& \multirow{2}{*}{split Weyl} & split Weyl \\

& & & & with a weighted function \\ \hline

\multirow{2}{*}{$y'''=\frac{3y'(y'')^2}{1+(y')^2}$} &
\multirow{2}{*}{$\S$} &
\multirow{2}{*}{------} & \raisebox{-0.25ex}{Lorentzian} & \multirow{2}{*}{------} \\

& & & \raisebox{0.3ex}{Ricci$=const> 0$} & \\ \hline

\multirow{2}{*}{$y'''=\frac{3(y'')^2}{2y'}$} &
\multirow{2}{*}{$\S$} & \multirow{2}{*}{------} &
\raisebox{-0.25ex}{Lorentzian} & \raisebox{-0.25ex}{Lorentzian} \\

& & & \raisebox{0.3ex}{Ricci$=const< 0$} & \raisebox{0.3ex}{with a
one-form} \\ \hline

\multirow{2}{*}{$y'''=-2\mu y'+y$} & \multirow{2}{*}{$\S$} &
\multicolumn{3}{|c|}{\multirow{2}{*}{$\real^2$ geometry}}
\\

& & \multicolumn{3}{|c|}{} \\ \hline

\end{tabular}
\end{table}
\end{landscape}

\chapter[Geometries of ODEs modulo contact transformations]
{Geometries of ODEs considered modulo contact transformations of variables}\label{ch.contact}

\section{Cartan connection on ten-dimensional bundle}\label{s.c.th}
\noindent We formulate the theorem about
the main structure which is associated with third-order ODEs
modulo contact transformations of variables, the $\o(3,2)$ Cartan
connection on the bundle $\P^c\to\J^2$. This structure will serve
as a starting point for both analyzing of geometries of ODEs and
their classification.

\begin{theorem}\label{th.c.1}
Consider an equation $y'''=F(x,y,y',y'')$. The contact invariant
information about this equation is given by the following data
\begin{itemize}
\item[i)] The principal fibre bundle $H_6\to\P^c\to\J^2$, where
$\dim\P^c=10$, $\J^2$ is the space of second jets of curves in the
$xy$-plane, and $H_6$ is the following six-dimensional subgroup of
$SP(4,\real)$ \be \label{e.c.H6}  H_6=\bma \sqrt{u_1}, &
\frac12\frac{u_2}{\sqrt{u_1}}, & -\frac12\frac{u_4}{\sqrt{u_1}}, &
 \tfrac{1}{24}\tfrac{u_2^2u_5}{u_1^{3/2}u_3}-\tfrac12\sqrt{u_1}\,u_6 \\\\
 0 & \tfrac{u_3}{\sqrt{u_1}}, & -\tfrac{u_5}{\sqrt{u_1}},
 & \tfrac12\tfrac{u_2u_5-u_3u_4}{u_1^{3/2}} \\\\
 0 & 0 & \tfrac{\sqrt{u_1}}{u_3}, &
-\tfrac12\tfrac{u_2}{\sqrt{u_1}\,\u_3} \\\\
  0 & 0 & 0 & \tfrac{1}{\sqrt{u_1}}\ema.
\ee \item[ii)] The coframe
$(\theta^1,\theta^2,\theta^3,\theta^4,\Omega_1,\Omega_2,\Omega_3,\Omega_4,\Omega_5,\Omega_6)$
on $\P^c$, which defines the $\o(3,2)$-valued Cartan normal
connection $\wh{\omega}^c$ on $\P^c$ by \be\label{e.c.conn_sp}
 \wh{\omega}^c=\bma \tfrac{1}{2}\vc{1} & \tfrac{1}{2}\vc{2} & -\tfrac{1}{2}\vc{4} & -\tfrac{1}{4}\vc{6} \\\\
             \hc{4} & \vc{3}-\tfrac{1}{2}\vc{1} & -\vc{5} & -\tfrac{1}{2}\vc{4} \\\\
             \hc{2} & \hc{3} & \tfrac{1}{2}\vc{1}-\vc{3} & -\tfrac{1}{2}\vc{2} \\\\
             2\hc{1} & \hc{2} & -\hc{4} & -\tfrac{1}{2}\vc{1}
        \ema.
\ee
\end{itemize}

Let $(x,y,p,q,u_1,u_2,u_3,u_4,u_5,u_6)$, $(x^i,u_\mu)$ for short,
be a local coordinate system on $\P^c$, which is compatible with
the local trivialization $\P^c=H_6\times\J^2$, that is
$(x^i)=(x,y,p,q)$ are coordinates in $\J^2$ and $(u_\mu)$ are
coordinates in $H_6$ as in \eqref{e.c.H6}. Then the value of
$\wh{\omega}^c$ at the point $(x^i,u_\mu)$ in $\P^c$ is given by
\ben \wh{\omega}^c(x^i,u_\mu)=u^{-1}\,\omega^c\,u+u^{-1}\der u
\een where $u$ denotes the matrix \eqref{e.c.H6} and \ben
{\omega}^c= \bma
\tfrac{1}{2}\vc{1}^0&\tfrac{1}{2}\vc{2}^0 & -\tfrac{1}{2}\vc{4}^0 & -\tfrac{1}{4}\vc{6}^0 \\\\
\omega^4 & \vc{3}^0-\tfrac{1}{2}\vc{1}^0 & -\vc{5}^0 & -\tfrac{1}{2}\vc{4}^0 \\\\
\omega^2 & \wt{\omega}^3 & \tfrac{1}{2}\vc{1}^0-\vc{3}^0 & -\tfrac{1}{2}\vc{2}^0 \\\\
2\omega^1 & \omega^2 & -\omega^4 & -\tfrac{1}{2}\vc{1}^0 \ema\een
is the connection $\wh{\omega}^c$ calculated at the point
$(x^i,u_1=1,u_2=0,u_3=1,u_4=0,u_5=0,u_6=0)$.
 The forms $\omega^1,\omega^2,\wt{\omega}^3,\omega^4$ read
\be\label{e.c.om}\begin{aligned}
\omega^1=&\der y-p\der x, \\
\omega^2=&\der p-q\der x, \\
\wt{\omega}^3=&\der q-F\der x-\tfrac13F_q(\der p-q\der x)+K(\der y-p\der x), \\
\omega^4=&\der x.
\end{aligned}
\ee The forms $\vc{1}^0,\ldots,\vc{6}^0$ read
\be\label{e.c.Om0}\begin{aligned}
\vc{1}^0=&-K_q\,\omega^1,  \\
\vc{2}^0=&\left(\tfrac13 W_q+L\right)\,\omega^1-K_q\,\omega^2-K\omega^4,  \\
\vc{3}^0=&-K_q\,\omega^1+\tfrac16F_{qq}\,\omega^2+\tfrac13F_q\omega^4, \\
\vc{4}^0=&-(\tfrac13W_{qq}+L_q)\,\omega^1+\tfrac12K_{qq}\,\omega^2,
\\
\vc{5}^0=&\tfrac12K_{qq}\,\omega^1-\tfrac16F_{qqq}\,\omega^2-\tfrac16F_{qq}\,\omega^4,
 \\
\vc{6}^0=&(\tfrac13\D(W_{qq})-\tfrac43W_{qp} -\tfrac13F_qW_{qq}
+\tfrac13F_{qqq}W+M)\,\omega^1+ \\
&+\tfrac13(F_{qqy}-F_{qqq}K-W_{qq})\,\omega^2-K_{qq}\,\wt{\omega}^3+ \\
&+(\tfrac23F_{qy}-\tfrac13F_{qq}K-2L-\tfrac43W_q)\,\omega^4.
\end{aligned}
\ee
\end{theorem}

In above theorem we used a concept of normal Cartan connection in
the sense of N. Tanaka \cite{Tan}. A normal Cartan connection is a
connection which takes value in a semisimple graded Lie algebra
and whose curvature satisfies some algebraic conditions, which
are, so to speak, a generalization of conditions for torsion in
the case of linear connections. We explain it below.
\begin{definition}
A semisimple Lie algebra $\g$ is graded if it has a vector space
decomposition \ben
\g=\g_{-k}\oplus\ldots\oplus\g_{-1}\oplus\g_0\oplus\g_1\oplus\ldots\oplus\g_{k}
\een such that \ben[\g_i,\g_j]\subset \g_{i+j}\een and
$\g_{-k}\oplus\ldots\oplus\g_{-1}$ is generated by $\g_{-1}$.
\end{definition}

Let us suppose that $\g$ is a semisimple graded Lie algebra and
denote $\m=\g_{-k}\oplus\ldots\oplus\g_{-1}$,
$\h=\g_0\oplus\ldots\oplus\g_k$. Let us consider a $\g$ valued
Cartan connection $\wh{\omega}$ on a bundle $H\to\P\to\M$, where
the Lie algebra of $H$ is $\h$. Fix a point $p\in\P$. The
decomposition $\g=\m\oplus\h$ defines in $T_p\P$ the complement
$\H_p$ of the vertical space $\V_p$. Therefore we have
$T_p\P=\V_p\oplus\H_p$, $\wh{\omega}(\V_p)=\h$ and
$\wh{\omega}(\H_p)=\m$. The curvature $\wh{K}_p=$
$(\der\wh{\omega}+\wh{\omega}\w\wh{\omega})_p$ at $p$ is then
characterized by the tensor $\kappa_p\in\Hom(\wedge^2\m,\g)$ given
by
 \be\label{e.c.kappa}\kappa_p(A,B)=\wh{K}_p(\wh{\omega}^{-1}_p(A), \wh{\omega}^{-1}_p(B)),\quad
 A,B\in\m.\ee
The function $\kappa\colon\P\to\Hom(\wedge^2\m,\g)$ is called the
structure function.

 In the space $\Hom(\wedge^2\m,\g)$ of
$\g$-valued two-forms let us
 define $\Hom^1(\wedge^2\m,\g)$ to be the space of all
$\alpha\in\Hom(\wedge^2\m,\g)$ fulfilling
$$\alpha(\g_i,\g_j)\subset \g_{i+j+1}\oplus\ldots\oplus\g_k \quad\text{for}\quad i,j<0.$$

Since the Killing form $B$ of $\g$ is non-degenerate and satisfies
$B(\g_p,\g_q)=0$ for $p\neq -q$, one can identify $\m^*$ with
$\g_1\oplus\ldots\oplus\g_k$. For a basis $(e_1,\ldots,e_m)$ of
$\m$ let $(e^*_1,\ldots,e^*_m)$ denote the unique basis of
$\g_1\oplus\ldots\oplus\g_k$ such that $B(e_i,e^*_j)=\delta_{ij}$.
Tanaka considered the following complex \ben
\ldots\longrightarrow\Hom(\wedge^{q+1}\m,\g)\overset{\partial^*}{\longrightarrow}\Hom(\wedge^q\m,\g)\longrightarrow\ldots
\een with
$\partial^*\colon\Hom(\wedge^{q+1}\m,\g)\to\Hom(\wedge^q\m,\g)$
given by the following formula \ben
\begin{aligned}
(\partial^*\alpha)&(A_1\w\ldots\w A_q)=\sum_i[e^*_i,\alpha(e_i\w
A_1\w\ldots\w A_q)]
 \\&+\tfrac{1}{2}\sum_{i,j}\alpha([e^*_j,A_i]_\m\w e_j\w A_1\w\ldots\w\hat{A_i}\w\ldots\w A_q),
\end{aligned}
\een where $\alpha\in\Hom(\wedge^q\m,\g)$, $A_1,\ldots\,A_q\in\m$,
$(e_i)$ is any basis in $\m$ and $[\,\, ,\,]_\m$ denotes the
$\m$-component of the bracket with respect to the decomposition
$\g=\m\oplus\h$. Finally, N. Tanaka \cite{Tan} introduced the
notion of normal connection, the definition below is given in the
language of \cite{Cap}.
\begin{definition}\label{def.Tanakanorm}
A Cartan connection $\wh{\omega}$ as above is normal if its
structure function $\kappa$ fulfills the following conditions \ben
\begin{aligned}
\text{i)}& &\qquad &\kappa\in\Hom^1(\wedge^2\m,\g),\\
\text{ii)}& &\qquad &\partial^*\kappa=0.
\end{aligned}
\een
\end{definition}
N. Tanaka considered above objects in a wider context of geometric
structures associated with Cartan connections. He proved the
one-to-one correspondence between normal connections on $\P\to\M$
and some $G$-structures, so called $G^{\#}_{~0}$-structures, on
$\M$. Starting from a $G^{\#}_{~0}$-structure one obtains a normal
connection by a procedure of reductions and prolongations, which
is a version of Cartan's equivalence method. It is the
correspondence with $G^\#_0$-structures which makes the normal
connections distinguished. N. Tanaka proved that a normal
connection exists and is unique if only $\g$ is a subalgebra of
$\g(\m,\g_0)$, so called prolongation of $\m$ and $\g_0$. The
notion of prolongation in Tanaka's sense would lead us to his
general theory, which is beyond the scope of this paper; in this
work we need conditions of normality to explicitly construct the
Cartan connection in the only case of theorem \ref{th.c.1}. For
details of Tanaka's theory see \cite{Tan2,Tan}.

\section{Proof of theorem 2.1} \label{s.c.proof}
\noindent We prove theorem \ref{th.c.1} by repeating Chern's
construction of Cartan connection for the contact equivalence
problem, supplemented later in \cite{Sat}.

We begin with the $G_c$-structure \eqref{e.G_c} on $\J^2$, which,
according to Introduction, encodes all the contact invariant
information about the underlying ODE. Let us fix on $\J^2\times
G_c$ a coordinate system $(x,y,p,q,u_1,\ldots,u_9)$ or
$(x^k,g^i_{j})$ for short, where $(x^k)=(x,y,p,q)$ are the
coordinates on $\J^2$ and
$$
 g^i_{j}=\bma u_1 & 0 & 0 & 0 \\ u_2 & u_3 & 0 & 0 \\ u_4 & u_5 & u_6 & 0 \\ u_8 & u_9 & 0 & u_7 \ema
$$
are coordinates on $G_c$. We remind that there are four well
defined forms on $G_c\times\J^2$, the components of the canonical
form $\theta^i$:
$$
\theta^i(x,g^{-1})=g^i_{~j}\,\omega^j(x),\qquad i=1,\ldots,4
$$
with $\omega^i$ given by \eqref{e.omega}. We seek a bundle on
which $\theta^i$ are supplemented to a coframe by certain new
one-forms $\Omega_\mu$ chosen in a well-defined geometric manner.

\subsection{S.-S. Chern's construction} The construction of the bundle $\P^c$ and
the coframe by Cartan's method is the following.
\subsection*{1)} We calculate the exterior derivatives of
$\theta^i$ on $G_c\times\J^2$ \be\label{e.c.red10}\begin{aligned}
 \der\theta^1=&\alpha_1\w\theta^1+T^1_{~jk}\theta^j\w\theta^k, \\
 \der\theta^2=&\alpha_2\w\theta^1+\alpha_3\w\theta^2+T^2_{~jk}\theta^j\w\theta^k, \\
 \der\theta^3=&\alpha_4\w\theta^1+\alpha_5\w\theta^2+\alpha_6\w\theta^3+T^3_{~jk}\theta^j\w\theta^k, \\
 \der\theta^4=&\alpha_8\w\theta^1+\alpha_9\w\theta^2+ \alpha_7\w\theta^4
 +T^4_{~jk}\theta^j\w\theta^k,
 \end{aligned} \ee
where $\alpha_\mu$ are the entries of the matrix $\der
g^i_{~k}\cdot g^{-1k}_{~j}$ and $T^i_{~jk}$ are some functions on
$G_c\times\J^2$. Next we collect $T^i_{~jk}\theta^j\w\theta^k$
terms
\begin{align}
 \der\theta^1=&\left(\alpha_1-T^1_{~12}\theta^2-T^1_{~13}\theta^3-T^1_{~14}\theta^4\right)\w\theta^1\nonumber \\
  &+T^1_{~23}\theta^2\w\theta^3+T^1_{~24}\theta^2\w\theta^4+T^1_{~34}\theta^3\w\theta^4,\nonumber \\
 \der\theta^2=&\left(\alpha_2-T^2_{~12}\theta^2-T^2_{~13}\theta^3-T^2_{~14}\theta^4\right)\w\theta^1 \nonumber\\
   &+\left(\alpha_3-T^2_{~23}\theta^3-T^2_{~24}\theta^4\right)\w\theta^2+T^2_{~34}\theta^3\w\theta^4,\label{e.c.red20} \\
 \der\theta^3=&\left(\alpha_4-T^3_{~12}\theta^2-T^3_{~13}\theta^3-T^3_{~14}\theta^4\right)\w\theta^1,\nonumber \\
   &+\left(\alpha_5-T^3_{~23}\theta^3-T^3_{~24}\theta^4\right)\w\theta^2
   +\left(\alpha_6-T^3_{~34}\theta^4\right)\w\theta^3\nonumber\\
 \der\theta^4=&\left(\alpha_8-T^4_{~12}\theta^2-T^4_{~13}\theta^3-T^4_{~14}\theta^4\right)\w\theta^1,\nonumber \\
   &+\left(\alpha_9-T^4_{~23}\theta^3-T^4_{~24}\theta^4\right)\w\theta^2
   +\left(\alpha_7+T^4_{~34}\theta^3\right)\w\theta^4\nonumber
 \end{align}
and introduce new 1-forms $\pi_\mu$ substituting the collected terms.
Eq. \eqref{e.c.red10} now read \be\label{e.c.red30}\begin{aligned}
 \der\theta^1&=\pi_1\w\theta^1+\frac{u_1}{u_3 u_7}\theta^4\w\theta^2, \\
 \der\theta^2&=\pi_2\w\theta^1+\pi_3\w\theta^2+\frac{u_3}{u_6u_7}\theta^4\w\theta^3, \\
 \der\theta^3&=\pi_4\w\theta^1+\pi_5\w\theta^2+\pi_6\w\theta^3,\\
 \der\theta^4&=\pi_8\w\theta^1+\pi_9\w\theta^2+\pi_7\w\theta^4,
\end{aligned}\ee
since $T^1_{~23}=T^1_{~34}=0$, $T^1_{~24}=-u_3u_7/u_1$ and
$T^2_{~34}=-u_3/(u_6u_7)$.

 The equations \eqref{e.c.red30} resemble structural equations
for a linear connection very much. When $\theta^i$ are components
of the canonical form and $\Gamma^j_{~k}$ are components of a
$\g$-valued connection then we have
\ben\der\theta^i+\Gamma^i_{~j}\w\theta^j=\tfrac12T^i_{~jk}\theta^j\w\theta^k\een
with a torsion $T$. In our case $\pi_\mu$ are the entries of the
matrix $\pi$ in $\g_c$
\ben\pi=\bma \pi_1 & 0 & 0 & 0 \\ \pi_2 & \pi_3 & 0 & 0 \\
\pi_4 & \pi_5 & \pi_6& 0\\
\pi_8 & \pi_9 & 0 & \pi_7 \ema. \een However, $\pi$ is not a
linear connection. This is because in this version of Cartan's
method $\pi$ usually does not transform regularly along fibres of
bundles ($G_c\times\J^2$ here) and the `curvature'
$\der\pi+\pi\w\pi$ is not necessarily a tensor, hence $\pi$ is not
a linear connection in general. We may think of $\pi$ as a
connection in a broader meaning, that is a horizontal distribution
on $G_c\times\J^2$, which is not necessarily right-invariant.
Keeping this in mind we will refer to $T^i_{~jk}$ as torsion. Thus
$\pi$ is a `connection' chosen by the demand that its torsion is
`minimal', i.e. possesses as few terms as possible.

We observe that $\pi_\mu$, which are candidates for the sought
forms $\Omega_\mu$, are not uniquely defined by equations
\eqref{e.c.red30}, for example the gauge $\pi_1\to\pi_1+f\theta^1$
leaves \eqref{e.c.red30} unchanged. Therefore our connection is
not uniquely defined by its torsion.

\subsection*{2)} In the next step we reduce the bundle $G_c\times\J^2$.
We choose its subbundle, say $\P^{(1)}$, characterized by the
property that the torsion coefficients are constant on it. We
choose $\P^{(1)}$ such that $T^1_{~24}=-1,\,T^2_{~34}=-1$ on it.
Thus $\P^{(1)}$ is defined by \be\label{e.c.red_u6u7}
  u_6=\frac{u_3^2}{u_1},\quad\quad\quad u_7=\frac{u_1}{u_3}.
\ee  It is known \cite{Gar,Ste} that such a reduction preserves
the equivalence, in other words, two bundles are equivalent if and
only if their respective reductions are. Here $\P^{(1)}$ has the
seven-dimensional structural group \ben G_c^{(1)}=\bma u_1 & 0 & 0
& 0
\\ u_2 & u_3 & 0 & 0 \\ u_4 & u_5 & \tfrac{u_3^2}{u_1} & 0 \\ u_8 & u_9 & 0
& \tfrac{u_1}{u_3} \ema. \een

\subsection*{3)} Next we pull-back $\theta^i$ and $\pi_\mu$ to
$\P^{(1)}$. But the new structural group $G^{(1)}_c$ is a
seven-dimensional subgroup of $G_c$, so
$(\theta^1,\ldots,\theta^4,\pi_1,\ldots,\pi_9)$ of
\eqref{e.c.red30} is not a coframe on $\P^{(1)}$ any longer, since
\ben \pi_6=2\pi_3-\pi_1 \mod(\theta^i), \quad\quad\quad
\pi_7=\pi_1-\pi_3
 \mod(\theta^i). \een
Taking this into account we recalculate \eqref{e.c.red30} and
gather the torsion terms. We choose the new connection
\ben\pi=\bma \pi_1 & 0 & 0 & 0 \\ \pi_2 & \pi_3 & 0 & 0 \\
\pi_4 & \pi_5 & 2\pi_3-\pi_1& 0\\
\pi_8 & \pi_9 & 0 & \pi_1-\pi_3 \ema \een so that its torsion is
minimal again. \be\label{e.c.red40}\begin{aligned}
 \der\theta^1&=\pi_1\w\theta^1+\theta^4\w\theta^2, \\
 \der\theta^2&=\pi_2\w\theta^1+\pi_3\w\theta^2+\theta^4\w\theta^3,\\
 \der\theta^3&=\pi_4\w\theta^1+\pi_5\w\theta^2+(2\pi_3-\pi_1)\w\theta^3
 +\left(\frac{3u_5}{u_3}-\frac{3u_2-u_3F_q}{u_1}\right)\theta^4\w\theta^3,\\
 \der\theta^4&=\pi_8\w\theta^1+\pi_9\w\theta^2+(\pi_1-\pi_3)\w\theta^4.
 \end{aligned}\ee

\subsection*{4)}
We repeat the steps {\bf 2)} and {\bf 3)}. Firstly we reduce
$\P^{(1)}$ to the subbundle $\P^{(2)}\subset\P^{(1)}$ defined by
the property that the only non-constant torsion coefficient
$T^3_{~34}$ in \eqref{e.c.red40} vanishes on it,
 \be\label{e.c.red_u5}
 u_5=\frac{u_3}{u_1}\left(u_2-\frac{1}{3}u_3F_q\right).
\ee Next we recalculate connection, re-collect the torsion and
make another reduction through the constant torsion condition ($K$
 is defined in \eqref{e.defK}.) \be\label{e.c.red_u4}
 u_4=\frac{u^2_3}{u_1}K+\frac{u_2^2}{2u_1}.
\ee

At this stage we have reduced the frame bundle $G_c\times\J^2$ to
the nine-dimensional subbundle $\P^{(3)}\to\J^2$, such that its
structural group is the following
\ben G^{(3)}_c=\bma u_1 & 0 & 0 & 0 \\
u_2 & u_3 & 0 & 0 \\
\tfrac{u_2^2}{u_1} & \tfrac{u_2u_3}{u_1} & \tfrac{u_3^3}{u_1} & 0 \\
u_8 & u_9 & 0 & \frac{u_1}{u_3} \ema \een and the frame dual to
$(\omega^1,\omega^2,\omega^3-\tfrac13F_q\omega^2+K\omega^1,\omega^4)$
belongs to $\P^{(3)}$. The structural equations on $\P^{(3)}$ read
after collecting \be\label{e.c.red50}\begin{aligned}
 \der\theta^1&=\pi_1\w\theta^1+\theta^4\w\theta^2, \\
 \der\theta^2&=\pi_2\w\theta^1+\pi_3\w\theta^2+\theta^4\w\theta^3, \\
 \der\theta^3&=\pi_2\w\theta^2+\left(2\pi_3-\pi_1\right)\w\theta^3
 +\frac{u_3^3}{u_1^3}W\theta^4\w\theta^1,\\
 \der\theta^4&=\pi_8\w\theta^1+\pi_9\w\theta^2+(\pi_1-\pi_3)\w\theta^4
\end{aligned}\ee
with some one-forms $\pi_1,\pi_2,\pi_3,\pi_8,\pi_9$. The function
$W$, defined in \eqref{e.defW}, is the W\"unsch\-mann invariant
which is a contact relative invariant for third-order ODEs. It
means, as we already explained in Introduction, that every contact
transformation applied to an ODE preserves the condition $W=0$ or
$W\neq 0$. It follows that third-order ODEs $y'''=F(x,y,y',y'')$
and $y'''=\cc{F}(x,y,y',y'')$, satisfying $W[F]= 0$ and
$W[\cc{F}]\neq 0$ respectively, are not contact equivalent.
Thereby, as Chern observed, third-order ODEs fall into two main
contact inequivalent branches: the ODEs satisfying $W\neq 0$, and
those satisfying $W=0$.

Equations \eqref{e.c.red50} do not still define the forms
$\pi_\mu$ uniquely but only modulo the following transformations
\begin{align}
 \pi_1&\to  \pi_1+ 2t_1\theta^1,\nonumber \\
 \pi_2&\to  \pi_2+t_1\theta^2,\nonumber \\
 \pi_3&\to  \pi_3+t_1\theta^1, \label{e.c.prol10}\\
 \pi_8&\to  \pi_8+t_2\theta^1+t_3\theta^2+t_1\theta^4,\nonumber \\
 \pi_9&\to  \pi_9+t_3\theta^1+t_4\theta^2.\nonumber
\end{align}
That is to say, the torsion in \eqref{e.c.red50} defines the
$\g^{(3)}_c$ connection $\pi$ only up to \eqref{e.c.prol10}.

At this point, there is no pattern of further reduction. If $W=0$
there are only constant torsion coefficients in \eqref{e.c.red50}
and we do not have any conditions to define a subbundle of
$\P^{(3)}$. In these circumstances we prolong $\P^{(3)}$.

\subsection*{5)} The idea of prolongation is the following.
On $\P^{(3)}$ there is no fixed coframe but only the coframe
$(\theta^1,\theta^2,\theta^3,\theta^4,\pi_1,\pi_2,\pi_3,\pi_8,\pi_9)$
given modulo \eqref{e.c.prol10}. But `a coframe given modulo $G$'
is a $G$-structure on $\P^{(3)}$. As a consequence we can deal
with this new structure on $P^{(3)}$ by means of the Cartan
method. Let us consider the bundle $\P^{(3)}\times G^{prol}$ then,
where  \ben G^{prol}=\bma 1 & 0 \\ t & 1 \ema \een reflects the
freedom \eqref{e.c.prol10} so that the block $t$ reads
\ben \bma 2t_1 & 0 & 0 & 0 \\ 0 & t_1 & 0 & 0 \\ t_1 & 0 & 0 & 0 \\
t_2 & t_3 & 0 & t_1 \\ t_3 & t_4 & 0 & 0 \ema. \een On
$\P^{(3)}\times G^{prol}$ there exist nine fixed one-forms
$\theta^1,\theta^2,\theta^3,\theta^4,\Pi_1,\Pi_2,\Pi_3,\Pi_8,\Pi_9$,
given by
\ben \bma \theta^i \\ \Pi_\mu \ema = \bma 1 & 0 \\
t & 1 \ema \bma \theta^i \\ \pi_\mu \ema, \een which is the
canonical one-form on $\P^{(3)}\times G^{prol}\to \P^{(3)}$.

\subsection*{6)} Now we apply the method of reductions to the
above structure on $\P^{(3)}\times G^{prol}$. We calculate the
exterior derivatives of $(\theta^i,\Pi_\mu)$. The derivatives of
$\theta^i$ take the form of \eqref{e.c.red50} with $\pi_\mu$
replaced by $\Pi_\mu$. The derivatives of $\Pi_\mu$, after
collecting and introducing 1-forms $\Lambda_K$ containing $\der
t_K$, read
\begin{align}
 \der\Pi_1=& \Lambda_1\w\theta^1+\Pi_8\w\theta^2-\Pi_2\w\theta^4,\nonumber \\
 \der\Pi_2=&  \tfrac{1}{2}\Lambda_1\w\theta^2-\Pi_1\w\Pi_2-\Pi_2\w\Pi_3+\Pi_8\w\theta^3+\frac{u_3^3}{u_1^3}W\Pi_9\w\theta^1\nonumber \\
   &+2f_1\theta^1\w\theta^3+f_4\theta^1\w\theta^4+f_2\theta^2\w\theta^3+f_5\theta^2\w\theta^4,\label{e.c.prol30}\\
 \der\Pi_3=&\tfrac{1}{2}\Lambda_1\w\theta^1+\Pi_8\w\theta^2+\Pi_9\w\theta^3+f_1\theta^1\w\theta^2+f_2\theta^1\w\theta^3+f_5\theta^1\w\theta^4+f_3\theta^2\w\theta^3,\nonumber \\
 \der\Pi_8=&\Lambda_2\w\theta^1+\Lambda_3\w\theta^1+\tfrac{1}{2}\Lambda_1\w\theta^4+\Pi_9\w\Pi_2+\Pi_8\w\Pi_3+f_2\theta^3\w\theta^4,\nonumber\\
 \der\Pi_9=&\Lambda_3\w\theta^1+\Lambda_4\w\theta^2+\Pi_1\w\Pi_9-2\Pi_3\w\Pi_9+\Pi_8\w\theta^4-f_1\theta^1\w\theta^4+f_3\theta^3\w\theta^4. \nonumber
\end{align}
where $f_1,f_2,f_3,f_4,f_5$ are functions. This time the forms
$\Lambda_K$ are interpreted as connection forms. We compute
$f_1,\ldots,f_5$ and choose the subbundle $\P^c$ of
$\P^{(3)}\times G^{prol}$ by the condition that $f_1,f_2,f_3$ are
equal to zero on $\P^c$. This is done by appropriate specifying of
parameters $t_2,t_3,t_4$ as functions of
$(x$,$y$,$p$,$q$,$u_1,u_2$,$u_3$,$u_8$,$u_9$,$t_1)$. We skip
writing these complicated formulae. The structural equations on
$\P^c$ read \be\label{e.c.prol40}\begin{aligned}
  \der\theta^1 =&\Pi_1\w\theta^1+\theta^4\w\theta^2, \\
  \der\theta^2 =&\Pi_2\w\theta^1+\Pi_3\w\theta^2+\theta^4\w\theta^3, \\
  \der\theta^3 =&\Pi_2\w\theta^2+(2\Pi_3-\Pi_1)\w\theta^3+A\,\theta^4\w\theta^1,  \\
  \der\theta^4 =&\Pi_8\w\theta^1+\Pi_9\w\theta^2+(\Pi_1-\Pi_2)\w\theta^4,  \\
  \der\Pi_1 =&\Lambda_1\w\theta^1+\Pi_8\w\theta^2-\Pi_2\w\theta^4,  \\
  \der\Pi_2 =&(\Pi_3-\Pi_1)\w\Pi_2+A\,\Pi_9\w\theta^1+\tfrac{1}{2}\Lambda_1\w\theta^2
   +\Pi_8\w\theta^3+B\,\theta^1\w\theta^4+C\,\theta^2\w\theta^4,   \\
  \der\Pi_3 =&\tfrac{1}{2}\Lambda_1\w\theta^1+\Pi_8\w\theta^2+\Pi_9\w\theta^3
   +C\,\theta^1\w\theta^4, \\
  \der\Pi_8 =&\Pi_9\w\Pi_2+\Pi_8\w\Pi_3-2C\,\Pi_9\w\theta^1+\tfrac{1}{2}\Lambda_1\w\theta^4
   +D\,\theta^1\w\theta^2+2E\,\theta^1\w\theta^3 \\
   &+G\,\theta^1\w\theta^4+H\,\theta^2\w\theta^3+J\,\theta^2\w\theta^4, \\
  \der\Pi_9 =&(\Pi_1-2\Pi_3)\w\Pi_9+\Pi_8\w\theta^4 +E\,\theta^1\w\theta^2+H\,\theta^1\w\theta^3
  +J\,\theta^1\w\theta^4 +L\,\theta^2\w\theta^3, \\
  \der\Lambda_1 =&\Lambda_1\w\Pi_1+2\Pi_8\w\Pi_2+2C\,\Pi_8\w\theta^1-2C\,\Pi_9\w\theta^2
  -A\,\Pi_9\w\theta^4+\wt{M}\,\theta^1\w\theta^2 \\
   &+2(D+AL)\,\theta^1\w\theta^3+\wt{N}\,\theta^1\w\theta^4+2E\,\theta^2\w\theta^3
   +G\,\theta^2\w\theta^4
\end{aligned}\ee
with certain functions $A,B,C,D,E,F,G,H,J,L,\wt{M},\wt{N}$ on
$\P^c$.

Above structural equations \emph{uniquely define the only
remaining auxiliary form $\Lambda_1$}. In this manner we
constructed the bundle $\P^c\to\J^2$ and the fixed coframe
associated to the ODEs modulo contact transformations. As we have
explained in Introduction, the functions $A,\ldots,\wt{N}$, and
their coframe derivatives are relative contact invariants for
third-order ODEs.

\subsection{Cartan normal connection from Tanaka's theory}\label{s.c.normalcon}
The above coframe is not fully satisfactory from the geometric
point of view since it does not transform equivariantly  along the
fibres of $\P^c\to\J^2$, that is to say, it does not define a
Cartan connection.

In order to see this we consider the simplest case, related to the
equation $y'''=0$, when all the functions $A,\ldots,\wt{N}$
vanish. Then \eqref{e.c.prol40} become the Maurer-Cartan equations
for the Lie algebra $\o(3,2)\cong\sp(4,\real)$ and $\P^c$ is
locally the Lie group $O(3,2)$. The Maurer-Cartan form on $\P^c$
in the four-dimensional defining representation of $\sp(4,\real)$
is given by \ben
 \wt{\omega}=\bma \tfrac{1}{2}\Pi_1 & \tfrac{1}{2}\Pi_{2} & -\tfrac{1}{2}\Pi_{8}& -\tfrac{1}{4}\Lambda_1 \\\\
             \hc{4} & \Pi_3-\tfrac{1}{2}\Pi_{1} & -\Pi_{9} & -\tfrac{1}{2}\Pi_{8} \\\\
             \hc{2} & \hc{3} & \tfrac{1}{2}\Pi_{1}-\Pi_{3} & -\tfrac{1}{2}\Pi_{2} \\\\
             2\hc{1} & \hc{2} & -\hc{4} & -\tfrac{1}{2}\Pi_{1}
        \ema.
\een Let $\h$ be the six-dimensional subalgebra of $\o(3,2)$
annihilated by the ideal
$<\theta^1,\,\theta^2,\,\theta^3,\,\theta^4>$ and let $H_6$ be the
connected simply-connected Lie subgroup of $SP(4,\real)$ with the
algebra $\h$. Then $H_6\to\P\to\J^2$ is a homogeneous space and
$\wt{\omega}$ is a flat Cartan connection of type
$(SP(4,\real),H_6)$ on $\P^c$.

However, this object is not a Cartan connection in a general case,
when $A,\ldots,\wt{N}$ do not vanish, since its curvature
$\wt{K}=\der \wt{\omega}+\wt{\omega}\w\wt{\omega}$ is not
horizontal  with respect to the fibration $\P\to\J^2$, that is
value of $\wt{K}$ on a vector tangent to a fibre of $\P\to\J^2$ is
not necessarily zero; for instance $\wt{K}^1_{~2}$ contains the
term $\tfrac12A \Pi_9\w\theta^1$. On the other hand the
horizontality of $\wt{K}$ is necessary and locally sufficient for
$\wt{\omega}$ to be a Cartan connection.

In order to resolve this problem H. Sato and Y. Yoshikawa
\cite{Sat} found the structural equations for the normal
connection in this problem by means of the Tanaka theory. We
recalculate their result in our notation and give explicit form of
the normal connection, which their paper does not contain.

Let $E^i_{~j}\in\gl(4,\real)$ denotes the matrix whose
$(i,j)$-component is equal to one and other components equal zero.
We introduce the following basis in $\sp(4,\real)$
\be\label{e.c.basis_sp}
\begin{aligned}
 &e_1=2E^4_{~1}, & &e_2=E^3_{~1}+E^4_{~2}, & &e_3=E^3_{~2} \\
 &e_4=E^2_{~1}-E^4_{~3}, & &e_5=\tfrac{1}{2}(E^1_{~1}-E^2_{~2}+E^3_{~3}-E^4_{~4}), &
 &e_6=\tfrac{1}{2}(E^1_{~2}-E^3_{~4}),\\
 &e_7=E^2_{~2}-E^3_{~3}, & &e_8=-\tfrac{1}{2}(E^1_{~3}+E^4_{~2}), & &e_9=-E^2_{~3}, \\
 && &e_{10}= -\tfrac{1}{4}E^1_{~4}. &&
\end{aligned}
\ee
The form $\wt{\omega}$ is given by
$$\wt{\omega}=\theta^1e_1+\theta^2 e_2+\theta^3 e_3+\theta^4 e_4+\Pi_1e_5+\Pi_2e_6+\Pi_3e_7+\Pi_8e_8+\Pi_9e_9+\Lambda_1e_{10}.$$
The algebra $\g$ has the following grading
$$ \o(3,2)=\g_{-3}\oplus\g_{-2}\oplus\g_{-1}\oplus\g_0\oplus\g_1\oplus\g_2\oplus\g_3, $$
where
\ben
\begin{aligned}
&\g_{-3}=<e_1>, & &\g_{-2}=<e_2>, & &\g_{-1}=<e_3,e_4>, & & \\
&\g_{0}=<e_5,e_7>, & &\g_{1}=<e_6,e_9>, & &\g_{2}=<e_8>, &
&\g_3=<e_{10}>
\end{aligned}
\een and
$\h=\g_0\oplus\g_1\oplus\g_2\oplus\g_3=<e_5,\ldots,e_{10}>$. Let
lower case Latin indices range from $1$ to $4$ and upper case
Latin indices range from $1$ to $10$ throughout this section.

Suppose now that $\wh{\omega}$ is an $\sp(4,\real)$ Cartan
connection on $\P^c$. Its structure function $\kappa$, defined in
\eqref{e.c.kappa}, decomposes into
$$\kappa=\tfrac{1}{2}\kappa^I_{~ij}\,e_I\otimes e^i\w e^j,$$
where $(e^I)$ denotes the basis  dual  to $(e_I)$ and
$\kappa^I_{~ij}=\kappa^I_{~[ij]}$ are functions. Condition i) of
definition \ref{def.Tanakanorm} reads
\begin{align*} &\kappa^1_{~23}=0,& &\kappa^1_{~24}=0, &&\kappa^1_{~34}=0, &&\kappa^2_{~34}=0. \end{align*}

\noindent We read structural constants $[e_I,e_J]=c^K_{~IJ} e_K$
for $\sp(4,\real)$, compute the Killing form $B_{IJ}$ and its
inverse $B^{IJ}$. The operator
$\partial^*\colon\Hom(\wedge^2\m,\g)\to\Hom(\m,\g)$ acts as
follows \ben
\partial^*(e_I\otimes e^i\w e^j)=
\left(2\delta^{[i}_{~m}B^{j]K}c^L_{~KI}-\delta^L_{~I}c^{[i}_{~Km}B^{j]K}\right)
e_L\otimes e^m.\een We apply $\partial^*$ to the basis
$(e_I\otimes e^i \w e^j)$ of $\Hom^1(\wedge^2\m,\g)$ and find that
the condition $\partial^*\kappa=0$ for the normality is equivalent
to vanishing of the following combinations of $\kappa^I_{~ij}$:
\begin{align}
&\bal &\kappa^1_{~14},&& \kappa^2_{~14},&&\kappa^2_{~24},\eal && \bal&\kappa^3_{~24},&&\kappa^3_{~34},&&\kappa^4_{~23},\eal \notag \\
 &\kappa^1_{~12}+\kappa^2_{~13}+\kappa^2_{~24}, && 2\kappa^1_{~12}-\kappa^4_{~24}-\kappa^5_{~34},  \notag \\
 &\kappa^1_{~12}-\kappa^3_{~23}-\kappa^7_{~34}, && \kappa^1_{~13}+\kappa^2_{~23}, \notag \\
 &\kappa^1_{~13}-\kappa^4_{~34}, && \kappa^2_{~12}+\kappa^4_{~14}+\kappa^5_{~24}, \notag \\
 &\kappa^2_{~12}-\kappa^3_{~13}-\kappa^7_{~24}, && \kappa^2_{~12}+\kappa^5_{~24}-\kappa^6_{~34}-\kappa^7_{~24},\notag \\
 &\kappa^2_{~13}+\kappa^3_{~23}+\kappa^5_{~34}-\kappa^7_{~34}, && \kappa^2_{~23}+\kappa^4_{~34},\notag \\
 &\kappa^3_{~12}-\kappa^5_{~14}+\kappa^6_{~24}+\kappa^7_{~14}, && \kappa^4_{~12}-\kappa^5_{~13}+2\kappa^7_{~13}+\kappa^9_{~24}, \notag \\
 &\kappa^4_{~12}-\kappa^6_{~23}-\kappa^8_{~34}-\kappa^9_{~24}, && \kappa^4_{~13}+\kappa^7_{~23}-\kappa^9_{~34},  \notag \\
 &\kappa^4_{~14}+\kappa^6_{~34}+\kappa^7_{~24}, && \kappa^4_{~24}-\kappa^5_{~34}+2\kappa^7_{~34},\notag \\
 &2\kappa^5_{~12}-2\kappa^8_{~24}-\kappa^{10}_{~34}, && \kappa^5_{~13}+\kappa^6_{~23}-\kappa^8_{~34}, \notag \\
 &\kappa^5_{~14}+\kappa^6_{~24},  && \kappa^5_{~23}-2\kappa^7_{~23}-\kappa^9_{~34}, \notag \\
 &2\kappa^6_{~12}+2\kappa^8_{~14}+\kappa^{10}_{~24}, && \kappa^6_{~13}+\kappa^7_{~12}+\kappa^8_{~24}+\kappa^9_{~14}. \notag
\end{align}
We write the sought normal Cartan connection $\wh{\omega}^c$ as
follows
$$\wh{\omega}^c=\theta^1e_1+\theta^2e_2+\theta^3e_3+\theta^4e_4
+\Omega_1e_5+\Omega_2e_6+\Omega_3e_7+\Omega_4e_8+\Omega_5e_9+\Omega_6e_{10}.$$
The forms $\Omega_\mu$, unknown yet, are given by
\begin{align*}
 &\Omega_1 = \Pi_1+a_i\theta^i, && \Omega_2 = \Pi_2+b_i\theta^i,
 && \Omega_3 =\Pi_3+c_i\theta^i, \nonumber \\
 &\Omega_4 =\Pi_8+f_i\theta^i, && \Omega_5 =\Pi_9+g_i\theta^i,
 && \Omega_6 = \Lambda_1+h_i\theta^i.
\end{align*}

Next, we calculate the curvature
$\wh{K}^c=\der\wh{\omega}^c+\wh{\omega}^c\w\wh{\omega}^c$, find
the components of the structure function and put them into the
normality conditions. These are only satisfied if all the
functions $a_i,b_i,c_i,e_i,f_i,g_i$ vanish except for
$c_1,h_1,h_2,h_3,h_4$ which are arbitrary and
\begin{align*}
 &a_1=2c_1, && b_1=\tfrac{2}{3}C,&& b_2=c_1, && f_1=\tfrac{2}{3}J, &&
 f_4=c_1,
\end{align*}
where $C$ and $J$ are the functions in \eqref{e.c.prol40}.
Finally, we obtain from the $e_1$-component of the Bianchi
identity
$\der\wh{K}^c=\wh{K}^c\w\wh{\omega}^c-\wh{\omega}^c\w\wh{K}^c$
that
\begin{align*}
 &c_1=0, && h_1=\tfrac{4}{3}G-\tfrac{2}{3}\wt{X}_4(J),&& h_2=\tfrac{2}{3}J, && h_3=0, &&
 h_4=-\tfrac{2}{3}C,
\end{align*}
where $\wh{X}_4$ is the vector field in the frame
$(\wt{X}_1,\wt{X}_2,\wt{X}_3,\wt{X}_4,\wt{X}_5,\wt{X}_6,\wt{X}_7,\wt{X}_8,\wt{X}_9,\wt{X}_{10})$
dual to
$(\theta^1,\theta^1,\theta^1,\theta^1,\Pi_1,\Pi_2,\Pi_3,\Pi_8,\Pi_9,\Lambda_1)$.
The normal connection $\wh{\omega}^c$ has been constructed. The
last thing we must do is renaming the coordinates $u_8\to u_4$,
$u_9\to u_5$ and choosing $u_6$ appropriately, so that formulae
\eqref{e.c.H6} -- \eqref{e.c.Om0} hold. This finishes the proof of
theorem \ref{th.c.1}.

\section{Geometries on ten-dimensional bundle}\label{s.c.bundle}
\subsection{Curvature and its interpretation}
We turn to discussion of consequences of theorem \ref{th.c.1}. As
we explained in Introduction, the basic object that contains the
contact invariants for the ODEs is the curvature
\ben\wh{K}^c=\der\wh{\omega}^c+\wh{\omega}^c\w\wh{\omega}^c. \een
It is given by the structural equations for the coframe
$\theta^1,\ldots,\Omega_6$.
\begin{align}
 \der\hc{1} =&\vc{1}\w\hc{1}+\hc{4}\w\hc{2},\nonumber \\
 \der\hc{2} =&\vc{2}\w\hc{1}+\vc{3}\w\hc{2}+\hc{4}\w\hc{3},\nonumber \\
 \der\hc{3}=&\vc{2}\w\hc{2}+(2\vc{3}-\vc{1})\w\hc{3}+\inc{A}{2}\hc{2}\w\hc{1}
    +\inc{A}{1}\hc{4}\w\hc{1}, \nonumber \\
 \der\hc{4} =&\vc{4}\w\hc{1}+\vc{5}\w\hc{2}+(\vc{1}-\vc{3})\w\hc{4},\nonumber \\
 \der\vc{1} =&\vc{6}\w\hc{1}+\vc{4}\w\hc{2}-\vc{2}\w\hc{4},\nonumber \\
 \der\vc{2}=&(\vc{3}-\vc{1})\w\vc{2}+\tfrac{1}{2}\vc{6}\w\hc{2}+\vc{4}\w\hc{3}
  +\inc{A}{3}\,\hc{1}\w\hc{2}+\inc{A}{4}\hc{1}\w\hc{4}, \nonumber \\
 \der\vc{3}=&\tfrac{1}{2}\vc{6}\w\hc{1}+\vc{4}\w\hc{2}+\vc{5}\w\hc{3}
  +\inc{A}{5}\hc{1}\w\hc{2}+\inc{A}{2}\hc{1}\w\hc{4},\label{e.c.dtheta_10d} \\
 \der\vc{4}=&\vc{5}\w\vc{2}+\vc{4}\w\vc{3}+\tfrac{1}{2}\vc{6}\w\hc{4}
  +(\inc{A}{6}+\inc{B}{2})\hc{1}\w\hc{2} +2\inc{B}{3}\hc{1}\w\hc{3} \nonumber \\
  &-\inc{A}{3}\hc{1}\w\hc{4}+\inc{B}{4}\hc{2}\w\hc{3} \nonumber \\
 \der\vc{5} =&(\vc{1}-2\vc{3})\w\vc{5}+\vc{4}\w\hc{4}
  +(\inc{A}{7}+\inc{B}{3})\hc{1}\w\hc{2}+\inc{B}{4}\hc{1}\w\hc{3} \nonumber \\
  &-\inc{A}{5}\hc{1}\w\hc{4}+\inc{B}{1}\hc{2}\w\hc{3}, \nonumber \\
  \der\vc{6} =&\vc{6}\w\vc{1}+2\vc{4}\w\vc{2}+\inc{C}{1}\hc{1}\w\hc{2}
   +2\inc{B}{2}\hc{1}\w\hc{3}+\inc{A}{8}\hc{1}\w\hc{4}+2\inc{B}{3}\hc{2}\w\hc{3}, \nonumber
\end{align}
where
$\inc{A}{1},\ldots,\inc{A}{8},\inc{B}{1},\ldots,\inc{B}{4},\inc{C}{1}$
are functions on $\P^c$.

The functions $\inc{A}{1},\ldots,\inc{C}{1}$ are contact relative
invariants of the underlying ODE, that is their vanishing or not
is a contact invariant property. The full set of contact
invariants can be constructed by consecutive differentiation of
$\inc{A}{1},\ldots,\inc{C}{1}$ with respect to the frame
$(X_1,X_2,X_3,X_4,X_5,X_6,X_7,X_8,X_9, X_{10})$ dual to
$(\theta^1,$ $\theta^2,$ $\theta^3,$ $\theta^4,$ $\Omega_1,$
$\Omega_2,$ $\Omega_3,$ $\Omega_4,$ $\Omega_5,$ $\Omega_6)$.
Utilizing the identities $\der^2\Omega_\mu=0$ we compute the
exterior derivatives of $\inc{A}{i},\inc{B}{j},\inc{C}{1}$, for
instance
$$ \der\inc{B}{1}=X_1(\inc{B}{1})\hc{1}+X_2(\inc{B}{2})\hc{2}+X_3(\inc{B}{3})\hc{3}
-2\inc{B}{4}\hc{4}+2\inc{B}{1}\vc{1}-5\inc{B}{1}\hc{3}.$$ From
these formulae it follows that i) $\inc{A}{2},\ldots\inc{A}{8}$
express by the coframe derivatives of $\inc{A}{1}$, ii)
$\inc{B}{2},\ldots\inc{B}{4}$ express by coframe derivatives of
$\inc{B}{1}$ and iii) $\inc{C}{1}$ is a function of coframe
derivatives of both $\inc{A}{1}$ and $\inc{B}{1}$. Hence there are
two basic invariants for the system \eqref{e.c.dtheta_10d} reading
\ben \inc{A}{1}=\frac{u_3^3}{u_1^3}W \qquad
\inc{B}{1}=\frac{u_1^2}{6u_3^5}F_{qqqq} \een and all other
invariants can be derived from them.\footnote{This property means
in language of Tanaka's theory that curvature of a normal
connection is generated by its harmonic part.}

The simplest case, in which all the contact invariants
$\inc{A}{1},\ldots,\inc{C}{1}$ vanish corresponds to
$W=F_{qqqq}=0$ and is characterized by
\begin{corollary}\label{cor.c.10d_flat}
For a third-order ODE $y'''=F(x,y,y',y'')$ the following
conditions are equivalent.
\begin{itemize}
 \item[i)] The ODE is contact equivalent to $y'''=0$.
 \item[ii)] It satisfies the conditions $W=0$, and $F_{qqqq}=0$.
  \item[iii)] It has the $\o(3,2)$ algebra of contact symmetry generators.
\end{itemize}
\end{corollary}

\subsection{Structure of $\P^c$}

The manifold $\P^c$ is endowed with threefold geometry of the
principal bundle over the second jet space $\J^2$, the first jet
space $J^1$ and the solution space $\S$. We discuss these
structures consecutively. Let us remind that
$(X_1,X_2,X_3,X_4,X_5,X_6,X_7,X_8,X_9, X_{10})$ denotes the dual
frame to $(\theta^1,$ $\theta^2,$ $\theta^3,$ $\theta^4,$
$\Omega_1,$ $\Omega_2,$ $\Omega_3,$ $\Omega_4,$ $\Omega_5,$
$\Omega_6)$.

First bundle structure, $H_6\to\P\to\J^2$, has been already
introduced explicitly in theorem \ref{th.c.1}. Here we only show
that it is actually defined by the coframe, since it can be
recovered from its structural equations. Indeed, we see from
\eqref{e.c.dtheta_10d} that the ideal spanned by
$\theta^1,\theta^2,\theta^3,\theta^4$ is closed \ben \der
\theta^i\w\theta^1\w\theta^2\w\theta^3\w\theta^4=0 \qquad
\text{for}\qquad i=1,2,3,4, \een and it follows that its
annihilated distribution $<X_5,X_6,X_7,X_8,X_9,X_{10}>$ is
integrable. Maximal integral leaves of this distribution are
locally fibres of the projection $\P^c\to\J^2$. Furthermore, the
commutation relations of these vector fields are isomorphic to the
commutation relations of the six-dimensional group $H_6$, hence we
can {\em define} $X_5,\ldots,X_{10}$ to be the fundamental vector
fields associated to the action $H_6$ on $\P^c$.

In order to explain how $\P^c$ is the bundle
$CO(2,1)\ltimes\real^3\to\P^c\to\S$ let us first describe the
space $\S$ itself. On $\J^2$ there is the congruence of solutions
of the ODE. A family of solutions passing through sufficiently
small open set in $\J^2$ is given by the mapping
$$(x;c_1,c_2,c_3)\mapsto
(x,f(x;c_1,c_2,c_3),f_x(x;c_1,c_2,c_3),f_{xx}(x;c_1,c_2,c_3)),$$
where $y=f(x;c_1,c_2,c_3)$ is the general solution to
$y'''=F(x,y,y',y'')$ and $(c_1,c_2,c_3)$ are constants of
integration. Thus a solution can be considered as a point in the
three-dimensional real space $\S$ parameterized by the constants
of integration. This space can be endowed with a local structure
of differentiable manifold if we {\em choose} a parametrization
$(c_1,c_2,c_3)\mapsto f(x;c_1,c_2,c_3)$ of the solutions and admit
only sufficiently smooth re-parameterizations
$(c_1,c_2,c_3)\mapsto(\tilde{c}_1,\tilde{c}_2,\tilde{c}_3)$ of the
constants. We always assume that $\S$ is locally a manifold.
Since $J^2$ is a bundle over $\S$ so is $\P^c$  and the fibres of
the projection $\P^c\to\S$ are annihilated by the closed ideal
$<\theta^1,\theta^2,\theta^3>$. On the fibres there act the vector
fields $X_4,X_5,X_6,X_7,X_8,X_9,X_{10}$, which form the Lie
algebra $\co(2,1)\semi{.}\real^3$ and thereby define the action of
$CO(2,1)\ltimes\real^3$ on $\P^c$.

Apart from the projection $\J^2\to\S$ there is also the projection
$\J^2\to\J^1$ that takes the second jet $(x,y,p,q)$ of a curve
into its first jet $(x,y,p)$. It gives rise to the third bundle
structure, $H_7\to\P^c\to\J^1$. Here the tangent distribution is
$<X_3,X_5,X_6,X_7,X_8,X_9,X_{10}>$ and it defines the action of a
seven-dimensional group $H_7$ which of course contains $H_6$.

It appears that, under some conditions, $\widehat{\omega}_c$ is
not only a Cartan connection over $\J^2$ but over $\S$ or $\J^1$
also.

\section{Conformal geometry on solution space}\label{s.c.conf}
\noindent The section on the conformal geometry on $\S$ only
contains results of P. Nurowski \cite{Nur2,Nur1}, see also
\cite{New3}. Let us define on $\P^c$ the symmetric
two-contravariant tensor field
\ben\wh{g}=(\theta^2)^2-2\theta^1\theta^3 \een of signature $(++-
\,0\,0\,0\,0\,0\,0\,0)$. The degenerate directions of $\wh{g}$ are
precisely those tangent to the fibres of $\P^c\to\S$ \ben
\wh{g}(X_i,\cdot)=0,\qquad \text{for}\qquad i=4,5,6,7,8,9,10. \een
The Lie derivatives of $\wh{g}$ along the degenerate directions
are as follows \begin{align}
 &L_{X_4}\wh{g}=\frac{u_3^3}{u_1^3}W(\theta^1)^2,
 \qquad\quad
 L_{X_7}\wh{g}=2\wh{g},  \label{e.c.lie1} \\
\intertext{and}
 &L_{X_i}\wh{g}=0 \qquad\text{for} \qquad i=5,6,8,9,10.  \label{e.c.lie2} \end{align}
Thus, if only $$W=0,$$ all the degenerate directions but $X_7$ are
isometries of $\wh{g}$ whereas $X_7$ is a conformal symmetry. This
allows us to {\em project} $\wh{g}$ to the Lorentzian conformal
metric $[g]$ on the solution space $\S$. Since the action of
$CO(2,1)\ltimes\real^3$ on $\S$ is not given explicitly, we can
not write down the explicit formula for $g$. We can only do this
in terms of the coordinate system $(x,y,p,q)$ on $\J^2$:
\begin{align}
g&=(\omega^2)^2-2\omega^1\wt{\omega}^3= \label{e.c.g} \\
&=(\der p-q\der x)^2
 -2(\der y-p\der x)(\der q-F\der x-\tfrac13F_q(\der p-q\der x)+K(\der y-p\der
 x)).\notag\end{align}
In order to find the explicit expression for $g$ in a coordinate
system $(c_1,c_2,c_3)$ on $\S$ we would have to find the general
solution $y=f(x;c_1,c_2,c_3)$ of the underlying ODE, then solve
the system \ben\left\{ \bal y=&f(x;c_1,c_2,c_3), \\
p=&f_{x}(x;c_1,c_2,c_3), \\ q=&f_{xx}(x;c_1,c_2,c_3) \eal \right.
\een with respect to $c_i$ and substitute these formulae into
\eqref{e.c.g}.

The condition $W=0$ means $\inc{A}{1}=0$ which causes
$\inc{A}{2}=\ldots=\inc{A}{8}=0$. Thus, structural equations
\eqref{e.c.dtheta_10d} do not contain the non-constant terms
proportional to $\theta^4$ and the curvature $\wh{K}^c$ is
horizontal over $\S$. As a consequence, $\wh{\omega}^c$ is a
connection over $\S$. It appears that it is nothing but the Cartan
normal conformal connection for $[g]$.

\subsection{Normal conformal connection} We introduce after
\cite{Nur2} the normal conformal connection and discuss its
curvature. Consider $\real^n$ with coordinates $(x^\mu)$,
$\mu=1,\ldots,n$ equipped with the flat metric $g=g_{\mu\nu}\der
x^\mu\otimes\der x^\nu$ of the signature $(k,l)$, $n=k+l$. The
group $Conf(k,l)$ of conformal symmetries of $g_{\mu\nu}$ consists
of
\begin{itemize}
\item[i)] the subgroup $CO(k,l)=\real\times O(k,l)$ containing the
group $O(k,l)$ of isometries of $g$ and the dilatations,
\item[ii)] the subgroup $\real^n$ of translations, \item[iii)] the
subgroup $\real^n$ of special conformal transformations.
\end{itemize}
The stabilizer of the origin in $\real^n$ is the semisimple
product of $CO(k,l)\ltimes\real^n$ of the isometries, the
dilatations, and the special conformal transformations. The flat
conformal space is the homogeneous space
$Conf(k,l)/CO(k,l)\ltimes\real^n$. To this space there is
associated the flat Cartan connection on the bundle
$CO(k,l)\ltimes\real^n\to Conf(k,l)\to \real^n$ with values in the
algebra $\conf(k,l)$.

By virtue of the M\"obius construction, the group $Conf(k,l)$ is
isomorphic to the orthogonal group
$O(k+1,l+1)$ preserving the metric \ben \bma 0 & 0 & -1 \\
0 & g_{\mu\nu} & 0 \\ -1 & 0 & 0 \ema \een on $\real^{n+2}$. This
isomorphism gives rise to the following representation of
$\conf(k,l)\cong\o(k+1,l+1)$ \be\label{e.c.o32} \bma \phi & g_{\nu\rho}\xi^\rho & 0 \\
v^\mu & \lambda^\mu_{~\nu} & \xi^\mu \\ 0 & g_{\nu\rho}v^\rho &
-\phi \ema. \ee Here the vector $(v^\mu)\in\real^n$ generates the
translations, the matrix $(\lambda^\mu_{~\nu})\in\o(k,l)$
generates the isometries, $\phi$ -- the dilatations, and
$(\xi^\mu)\in\real^n$ -- the special conformal transformations.

Let us turn to an arbitrary case of a conformal metric $[g]$ of
the signature $(k,l)$ on a $n$-dimensional manifold $\M$,
$n=k+l>2$. Let us choose a representative $g$ of $[g]$ and
consider an orthogonal coframe $(\omega^\mu)$, in which
$g=g_{\mu\nu}\,\omega^{\mu}\otimes\omega^\nu$ with constant
coefficients $g_{\mu\nu}$. We calculate the Levi-Civita connection
$\Gamma^\mu_{~\nu}$ for $g$, its Ricci tensor $R_{\mu\nu}$ and the
Ricci scalar $R$. Next we define the following one-forms \ben
P_\nu=\left( \tfrac{1}{2-n}R_{\nu\rho}+\tfrac{1}{2(n-1)(n-2)}R
g_{\nu\rho} \right)\theta^\rho. \een Given these objects we build
the following $\conf(k,l)$-valued matrix $\omega^{conf}$ on $M$
\be\label{e.cnc} \omega^{conf}=\bma
0&P_\nu&0\\
\theta^\mu&\Gamma^\mu_{~\nu}&g^{\mu\rho} P_\rho\\
0&g_{\nu\rho}\theta^\rho&0 \ema. \ee This is the normal conformal
connection on $\M$ in the natural gauge.\footnote{The gauge is
natural since we have started from the Levi-Civita connection, not
from any Weyl connection for $g$, in which case \eqref{e.cnc}
contains a Weyl potential.} Consider now the conformal bundle
$CO(k,l)\ltimes\real^n\to\P\to\M$, and choose a coordinate system
$(h,x)$ on $\P$ compatible with the local triviality $\P\cong
CO(k,l)\ltimes\real^n \times\M$, where $x$ stands for $(x^\mu)$ in
$\M$ and the matrix $h\in CO(k,l)\ltimes\real^n$ reads \ben
h=\begin{pmatrix} {\rm e}^{-\phi}&{\rm
e}^{-\phi}g_{\nu\rho}\xi^\rho &
\frac{1}{2}{\rm e}^{-\phi}\xi^\rho\xi^\sigma g_{\rho\sigma}\\
0&\Lambda^{\mu}_{~\nu}&\Lambda^\mu_{~\rho}\xi^\rho\\
0&0&{\rm e}^\phi
\end{pmatrix},
\quad\quad\Lambda^\mu_{~\rho}\Lambda^\nu_{~\sigma}g_{\mu\nu}=g_{\rho\sigma}.
 \een

The normal conformal connection for $g$ is the following
$\conf(k,l)$-valued one-form on $\P$ \ben
\wh{\omega}^{conf}(h,x)=h^{-1}\cdot\pi^*(\omega^{conf}(x))\cdot h
+h^{-1}\der h. \een The curvature of the normal conformal
connection is as follows \ben
\wh{K}^{conf}(h,x)=h^{-1}\cdot\pi^*(K^{conf}(x))\cdot h, \een
where $K^{conf}$ is the curvature for $\omega^{conf}$ on $\M$ \ben
K^{conf}=\bma
0&D P_\nu&0\\
0& C^\mu_{~\nu} &g^{\mu\rho}D P_\rho\\
0&0&0 \ema, \een and  \ben DP_\mu=\der
P_\mu+P_\nu\w\Gamma^\nu_{~\mu}=\tfrac12P_{\mu\nu\rho}\omega^\nu\w\omega^\rho.
\een The curvature contains the lowest-order conformal invariants
for $g$, namely
\begin{itemize}
\item For $n\geq4$ \ben C^\mu_{~\nu}=\tfrac12
C^\mu_{~\nu\rho\sigma}\omega^\rho\w\omega^\sigma\een is the Weyl
conformal tensor, while \ben
P_{\mu\nu\rho}=\tfrac{1}{n-3}\nabla_{\sigma}C^\sigma_{~\mu\nu\rho}\een
is its divergence. \item For $n=3$ the Weyl tensor identically
vanishes, $C^\mu_{~\nu}= 0$, and the lowest-order conformal
invariant is the Cotton tensor $P_{\mu\nu\rho}$. It has five
independent components.
\end{itemize}

The normality of conformal connections, originally defined by E.
Cartan, is the following property. The algebra $\conf(k,l)$ is
graded: $\conf(k,l)=\g_{-1}\oplus\g_0\oplus\g_1$, where
translations are the $\g_{-1}$-part, $\co(k,l)$ is the $\g_0$-part
and the special conformal transformations are the $\g_1$-part. The
normal connection for $[g]$ is the only conformal connection such
that the $\co(k,l)$-part of its curvature, $C^\mu_{~\nu}=\tfrac12
C^\mu_{~\nu\rho\sigma}\omega^\rho\w\omega^\sigma$, is traceless:
$C^\rho_{~\nu\rho\sigma}=0$. Cartan normal conformal connections
are normal in the sense of Tanaka.

\subsection{Normal conformal connection from ODEs}
As we mentioned, $\wh{\omega}^c$ is the normal conformal
connection over $\S$. In order to see this it is enough to
rearrange $\wh{\omega}^c$ according to the representation
\eqref{e.c.o32} \be\label{e.c.conn_o}
 \wh{\omega}^c=\bma \vc{3} & -\tfrac{1}{2}\vc{6} & -\vc{4} & -\vc{5} & 0 \\
             \hc{1} & \vc{3}-\vc{1} & -\hc{4} & 0 & -\vc{5} \\
             \hc{2} & -\vc{2} & 0 & -\hc{4} & \vc{4} \\
             \hc{3} & 0 &-\vc{2} & \vc{1}-\vc{3} & -\tfrac{1}{2}\vc{6} \\
             0 & \hc{3} & -\hc{2} & \hc{1} & -\Omega^3
        \ema
\ee and calculate its curvature, which is equal to \ben
 \wh{K}^c=\bma 0 & DP_1 & DP_2 & DP_3 & 0 \\
             0 & 0 & 0 & 0 & DP_3 \\
             0 & 0 & 0 & 0 & -DP_2 \\
             0 & 0 & 0 & 0 & DP_1 \\
             0 & 0 & 0 & 0 & 0
        \ema \een
with the following components of the Cotton tensor on $\P^c$
\begin{align*}
 DP_1=&-\tfrac12\inc{C}{1}\,\hc{1}\w\hc{2}-\inc{B}{2}\,\hc{1}\w\hc{3}-\inc{B}{3}\,\hc{2}\w\hc{3}, \notag \\
 DP_2=&-\inc{B}{2}\,\hc{1}\w\hc{2}-2\inc{B}{3}\,\hc{1}\w\hc{3}-\inc{B}{4}\,\hc{2}\w\hc{3},  \\
 DP_3=&-\inc{B}{3}\,\hc{1}\w\hc{2}-\inc{B}{4}\,\hc{1}\w\hc{3}-\inc{B}{1}\,\hc{2}\w\hc{3}. \notag
\end{align*}
Finally, we pull-back these formula to $\J^2$ through $u_1=u_3=1$,
$u_2=u_4=u_5=u_6=0$ and get

\begin{align*}
 DP_1=&(\tfrac12M_p+\tfrac16F_qM_q+\tfrac16F_{qqq}K_y+K_qL_q-\tfrac16 K^2
 F_{qqqq}+\notag \\
 &+\tfrac16K_qF_{qqy}-\tfrac16F_{qqyy}-\tfrac13F_{qqq}K_qK+\tfrac13 F_{qqy}K)\omega^1\w\omega^2
 \notag \\
 &+\tfrac12\left(M_q-K_{qqq}K-2K_{qq}K_q+K_{qqy}\right)\omega^1\w\wt{\omega}^3+\notag \\
 &-\tfrac12 L_{qq}\omega^2\w\wt{\omega}^3, \\
 DP_2=&\tfrac12\left(M_q-K_{qqq}K-2K_{qq}K_q+K_{qqy}\right)\omega^1\w\omega^2+\notag \\
 &-L_{qq}\omega^1\w\wt{\omega}^3+\tfrac12 K_{qqq}\omega^2\w\wt{\omega}^3, \notag \\
 DP_3=&-\tfrac12 L_{qq}\omega^1\w\omega^2+\tfrac12 K_{qqq}\omega^1\w\wt{\omega}^3
 -\tfrac16F_{qqqq}\omega^2\w\wt{\omega}^3. \notag
\end{align*}
The formulae for the conformal connection and curvature (in a
slightly different notation) are given in \cite{Nur2,New3}.

\section{Contact projective geometry on first jet
space}\label{s.c-p} \noindent The connection $\wh{\omega}^c$ gives
rise to not only the above conformal structure but also a geometry
on the first jet space $\J^1$. As we know, there are two basic
contact invariants in the curvature $\wh{K}^c$: $W$ and
$F_{qqqq}$. The condition $W=0$ yields the conformal geometry. Let
us now examine the second possibility
$$F_{qqqq}=0.$$ Quick inspection of the structural equations
\eqref{e.c.dtheta_10d} shows that in this case the curvature does
not contain $\theta^i\w\theta^4$ terms and thereby $\wh{\omega}^c$
is an $\sp(4,\real)$ Cartan connection on $H_7\to\P^c\to\J^1$.

A natural question is to what geometric structure  $\wh{\omega}^c$
is now related. In order to answer it let us look at the geometry
defined on $\J^1$ by the solutions of the underlying ODE
$y'''=F(x,y,y',y'')$. As we have already said the solutions form a
congruence on $\J^2$ since there passes exactly one solution
through a point $(x_0,y_0,p_0,q_0)\in\J^2$. The $G_c$-structure on
$\J^2$ defined by \eqref{e.G_c} preserves this congruence and also
preserves the contact invariant information about the ODE.  The
geometry on $\J^2$ is then of first order, since it is
characterized by the group $G_c$ acting on the tangent space
$T\J^2$.

 Geometry of an
ODE on $\J^1$ is of different kind. The family of solutions is no
longer a congruence. Indeed, through a fixed
$(x_0,y_0,p_0)\in\J^1$ there pass many solutions, each of them
corresponding to some value of $q_0$, which is given by a choice
of a tangent direction at $(x_0,y_0,p_0)$. Such a choice, however,
can not be made in an arbitrary way; a solution $y=f(x)$ lifts to
a curve $x\mapsto(x,f(x),f'(x))$ in $\J^1$, whose tangent vector
field $\partial_x+f'(x)\partial_y+f''(x)\partial_p$ is always
annihilated by the one-form $\omega^1=\der y-p\der x$ on $\J^1$.
In this manner all the admissible tangent vectors to the solutions
form a rank-two distribution $\C$, {\em the contact distribution
on $\J^1$}, which is annihilated by $\omega^1$.  This is first
difference between $\J^2$ and $\J^1$: we have the rank-two
distribution $\C$ instead of the congruence. Second and more
important difference is that we still need to distinguish the
family of solutions among a class of all curves with their tangent
fields in $\C$.

There are at least two basic methods of doing this. In the first
method we treat the tangent directions of curves in $\J^1$ as new
dimensions and move to the space where over each point
$(x,y,p)\in\J^1$ there is the fibre of all the possible tangent
directions. In this manner the `entangled' family of solutions in
$\J^1$ would be `stretched' to a congruence in the new space.
However, the bundle of tangent directions of $\J^1$ is nothing but
the second jet space $\J^2$ and thereby we would come back to the
description given in the theorem \ref{th.c.1}.

The other way is describing the family of solutions in $\J^1$ as a
family of unparameterized geodesics of a linear connection in
$\J^1$. This approach leads us to the notion of projective
differential geometry.

\subsection{Contact projective geometry}\label{s.cp3}
This geometry has been exhaustively analyzed in \cite{Fox}, see
also \cite{Cap,Cap2}. We will not discuss the general theory here
but focus on an application of the three-dimensional case to the
ODEs. The definition of contact-projective geometry, see D. Fox
\cite{Fox}, adapted to our situation is the following.

\begin{definition}\label{def.c.cp}
A contact projective structure on the first jet space $\J^1$ is
given by the following data.
\begin{itemize}
\item[i)] The contact distribution $\C$, that is a distribution
annihilated by $$\omega^1=\der y -p\der x.$$ \item[ii)] A family
of unparameterized curves everywhere tangent to $\C$ and such
that: a) for a given point and a direction in $\C$ there is
exactly one curve passing through that point and tangent to that
direction, b) curves of the family are among unparameterized
geodesics for some linear connection on $\J^1$.
  \end{itemize}
\end{definition}

A contact projective structure on $\J^1$ is equivalently given by
a family of linear connections, whose geodesic spray contains the
family of curves. For $\nabla$ to belong to this class one needs
\be\label{e.cp.geod}\nabla_V V=\lambda (V) V \ee along every curve
in the family, where $X$ denotes a tangent field to the considered
curve and $\lambda(X)$ is a function. Given two such connections
$\nabla$ and $\wt{\nabla}$, their difference is a $(2,1)$-type
tensor field $$A(X,Y)=\wt{\nabla}_X Y-\nabla_X Y,$$ for all $X$
and $Y$. Simultaneously we have  $A(V,V)=\mu(V) V$ for $V\in\C$,
where $\mu(V)=\wt{\lambda}(V)-\lambda(V)$ and $\mu_w$ at a point
$w\in\J^1$ is a covector on the vector space $\C_w$. By
polarization we obtain \be\label{e.c.A}
A(X,Y)+A(Y,X)=\mu(X)Y+\mu(Y)X, \qquad X,Y\in\C.\ee The connections
associated to a contact projective structure, when considered at a
point, form an affine space characterized by the above $A$.

Let us describe $A$ more explicitly. We choose a frame
$(e_1=\partial_y, e_2=\partial_p, e_3=\partial_x+p\partial_y)$ and
denote the dual frame by $(\sigma^1,\sigma^2,\sigma^3)$. In
particular $\omega^1=\sigma^1$ and $\C=<e_2,e_3>$. Let
$i,j,\ldots=1,2,3$ and $I,J,\ldots=2,3$. Now $\nabla_j e_i=
\Gamma^k_{~ij} e_k$,
$A^k_{~ij}=\wt{\Gamma}^k_{~ij}-\Gamma^k_{~ij}$, $\mu=\mu^I e_I$
and \eqref{e.c.A} reads \ben A^k_{~(IJ)}=\mu_{(I}
\delta^{k}_{~J)}. \een Relevant components are equal to
\begin{align}\label{e.c.Aco}
&A^1_{~22}=0, &&A^1_{~(23)}=0, &&A^1_{~33}=0, \notag \\
&A^2_{~22}=\mu_2, &&A^2_{~(23)}=\tfrac12\mu_3,&& A^1_{~33}=0,  \\
&A^3_{~22}=0, &&A^3_{~(23)}=\tfrac12\mu_2, && A^1_{~33}=\mu_3
\notag \end{align} and the rest of $A^k_{~ij}$ is free. Elementary
calculations assure us that the class of admissible connections is
a $20$-dimensional subspace of  $27$-dimensional space of all
linear connections on $\J^1$. Another constraint for the
connections is given by \be\label{e.c.ompr} (\nabla_V
\omega^1)V=\nabla_V(\omega^1(V))=0,\qquad V\in\C. \ee In our frame
this is equivalent to $\Gamma^k_{(IJ)}=0$. Combining eq.
\eqref{e.c.Aco} and \eqref{e.c.ompr} we obtain
\begin{proposition}\label{prop.cp.coord}
The following quantities are invariant with respect to a choice of
a connection in the class distinguished by a contact projective
structure on $\J^1$
\begin{subequations}\label{e.cp.niezm}
\begin{align}
&\Gamma^1_{~22}=0, && \Gamma^1_{~(23)}=0, && \Gamma^1_{~33}=0, && \label{e.cp.niezm1} \\
&\Gamma^3_{~22}, && 2\Gamma^3_{~(23)}-\Gamma^2_{~22}, &&
\Gamma^3_{~33}-2\Gamma^2_{~(23)}, && \Gamma^2_{~33}.
\label{e.cp.niezm2}
\end{align}
\end{subequations}
The connection coefficients are calculated in a frame $(e_i)$ such
that $\C=<e_2,e_3>$.

Values of four the unspecified combinations \eqref{e.cp.niezm2}
define a contact projective structure.
\end{proposition}

Among the above connections there is a distinguished subclass of
those connections which covariantly preserve the distribution
$\C$. We shall call them compatible connections. They satisfy not
only \eqref{e.c.ompr} but a stronger condition \ben \nabla_X
\omega^1 =\phi(X) \omega^1, \qquad \text{for all } X,\een with
some one-form $\phi$. We have the following
\begin{proposition}
A compatible connection has non-vanishing torsion
\begin{proof}
\begin{align*}
\der \omega^1(X,Y)&=\tfrac12((\nabla_X \omega^1) Y
-(\nabla_Y\omega^1)X+\omega^1(T(X,Y)))= \\
&=(\phi\w\omega^1)(X,Y)+\tfrac12\omega^1(T(X,Y))), \\
\intertext{thus}
(\der\omega^1\w\omega^1)(X,Y,Z)&=\tfrac12(\omega^1(T(X,Y))\omega^1(Z)+{\rm
cycl. perm.})\neq0.
\end{align*}
\end{proof}
\end{proposition}

\subsection{Contact projective geometries from ODEs}
It is obvious that the family of solutions of an arbitrary
third-order ODE satisfies the conditions i) and ii a) of
definition \ref{def.c.cp}. (Condition ii a) is satisfied with the
possible exception of the direction $\partial_p$, which belongs to
$\C$ but it is not tangent to any solution in general. However,
this exception is irrelevant since our consideration is local on
$T\J^1$.) We ask when the solutions form a subfamily of geodesics
for a linear connection.
\begin{lemma} \label{lem.c-p}
A third-order ODE $y'''=F(x,y,y',y'')$ defines a
contact-projective structure on $\J^1$ if and only if
$F_{qqqq}=0$. Moreover, the quantities \eqref{e.cp.niezm2} are
given by \be\label{e.cp.gam}\bal &\Gamma^3_{~22}=a_3, &&
2\Gamma^3_{~(23)}-\Gamma^2_{~22}=a_2, \\
&\Gamma^3_{~33}-2\Gamma^2_{~(23)}=a_1, && \Gamma^2_{~33}=-a_0,
\eal\ee where \ben F=a_3q^3+a_2q^2+a_1q+a_0. \een
\begin{proof}
The field $V=\tfrac{\der}{\der x}$ tangent to a solution
$(x,f(x),f'(x))$ equals $V=f''e_2+e_3$ in the frame $(e_i)$. The
geodesic equations \eqref{e.cp.geod} read
\begin{align*}
 & (f'')^2\,\Gamma^1_{~22}+2f''\Gamma^1_{~(23)}+\Gamma^1_{~33}=0, \\
 & f'''+(f'')^2\,\Gamma^2_{~22}+2f''\Gamma^2_{~(23)}+\Gamma^2_{~33}=\lambda(V) f'', \\
 & (f'')^2\,\Gamma^3_{~22}+2f''\Gamma^3_{~(23)}+\Gamma^3_{~33}=\lambda(V).
 \end{align*}
First of these equations is equivalent to \eqref{e.cp.niezm1}.
From the remaining equations we have that
$$
f'''=\Gamma^3_{~22}
{f''}^3+(2\Gamma^3_{~(23)}-\Gamma^2_{~22}){f''}^2
+(\Gamma^3_{~33}-2\Gamma^2_{~(23)})f''-\Gamma^2_{~33},
$$
is satisfied along every solution.
\end{proof}
\end{lemma}

We observe that the condition $F_{qqqq}=0$ yields
$\inc{B}{1}=\inc{B}{2}=\inc{B}{3}=\inc{B}{4}=0$, which removes all
$\theta^i\w\theta^3$ terms in the curvature and turns
$\wh{\omega}^c$ into a connection over $\J^1$, since in the
curvature there are only terms proportional to
$\theta^1,\theta^2,\theta^4$, horizontal with respect to
$\P^c\to\J^1$. Furthermore, the algebra $\sp(4,\real)\cong\o(3,2)$
has the following grading (apart from those already mentioned)
\begin{align}
&\sp(4,\real)=\g_{-2}\oplus\g_{-1}\oplus\g_0\oplus\g_1\oplus\g_2,
\label{c.gradJ1} \\ \intertext{which reads in the basis
\eqref{e.c.basis_sp}:}
&\g_{-2}=<e_1>,  \qquad \g_{-1}=<e_2,e_4 >, \notag \\
&\g_{0}=<e_3,e_5,e_7,e_9>, \notag  \\
&\g_{1}=<e_6,e_8>, \qquad \g_{2}=<e_{10}> \notag.
\end{align}
After calculating Tanaka's normality conditions by the method of
chapter \ref{ch.contact} section \ref{s.c.normalcon}, we observe
that $\wh{\omega}^c$ is the normal connection with respect to the
grading \eqref{c.gradJ1}. In this manner we have re-proved the
known fact that to a three-dimensional contact projective geometry
there is associated the unique normal $\sp(4,\real)$-valued Cartan
connection.
\begin{proposition}
If the contact projective geometry on $\J^1$ exists, then
$\wh{\omega}^c$ of theorem \ref{th.c.1} is the normal Cartan
connection for this geometry.
\end{proposition}

From $\wh{\omega}^c$ one may reconstruct the compatible
connections. To do this we just observe that first, second and
fourth equation of
\eqref{e.c.dtheta_10d} can be written as \ben \bma \der \hp{1} \\
\der \hp{2}\\ \der \hp{4} \ema
+\underbrace{\bma -\Omega_1 & 0 & 0 \\
      -\Omega_2 & -\Omega_3 & \theta^3 \\
      -\Omega_4 & -\Omega_5 & \Omega_3-\Omega_1 \ema}_{\displaystyle
      \wh{\Gamma}}
      \wedge
\bma  \hp{1} \\  \hp{2}\\  \hp{4} \ema = \bma  \hp{4}\w\hp{2}
\\\hp{4}\w\hp{3} \\ 0 \ema. \een
 The tree by three matrix denoted by $\wh{\Gamma}$ is the
$\g_0\oplus\g_1$-part of $\wh{\omega}^c$. The following
proposition holds.
\begin{proposition}For any section $s\colon\J^1\to\P^c$ the pull-back
$s^*\wh{\Gamma}$ written in the coframe
$(s^*\theta^1,s^*\theta^2,s^*\theta^4)$ is a connection compatible
with the contact projective geometry.
\begin{proof}
First we choose the section $s_0\colon\J^1\to\P^c$ given by $q=0$,
$u_1=1$, $u_3=1$ and $u_2=u_4=u_5=u_6=0$. We denote
$\Gamma=s^*_0\wh{\Gamma}$. In the coframe
$\sigma^1=s^*_0\theta^1$, $\sigma^2=s^*_0\theta^2$ and
$\sigma^3=s^*_0\theta^4$ we have
$-s^*_0\Omega_3=\Gamma^2_{~2}=\Gamma^2_{~2k}\sigma^k$,
$s^*_0\theta^3=\Gamma^2_{~3}=\Gamma^2_{~3k}\sigma^k$ and so on.
Equations \eqref{e.cp.gam} follow from \eqref{e.c.om} and
\eqref{e.c.Om0}, provided that $F_{qqqq}=0$.

Next we consider an arbitrary section $s\colon\J^1\to\P^c$. In the
local trivialization $\P^c\cong H_7\times \J^1$ we have $\P^c\ni
w=(v,x)$, where $v\in H_7$, $x\in\J^1$ and $s$ is given by
$x\mapsto (v(x),x)$. Now
$s^*\wh{\omega}^c(x)=v^{-1}(x)s^*_0\wh{\omega}^c(x)v(x)
+v^{-1}(x)\der v(x)$, and $s^*\wh{\Gamma}$ is the $\g_0\oplus\g_1$
part of $s^*\wh{\omega}^c$. Since the Lie algebra of $H_7$ is
$\g_0\oplus\g_1\oplus\g_2$, every $v(x)$ in the connected
component of the identity may be written in the form
$v(x)=v_2(x)v_1(x)=\exp{(t_2(x)A_2(x))}\exp{(t_1(x)A_1(x))}$ with
$A_2(x)\in\g_2$ and $A_1(x)\in\g_0\oplus\g_1$. It follows that
\ben
s^*\wh{\omega}^c(x)=v^{-1}(x)_2\left\{v^{-1}_1(x)s^*_0\wh{\omega}^c(x)v_1(x)
+v^{-1}_1(x)\der v_1(x)\right\}v_2(x) +v^{-1}_2(x)\der v_2(x)\een
 But the
$\g_0\oplus\g_1$ part of the quantity in the curly brackets is the
connection $\Gamma=s^*_0\wh{\Gamma}$ written in the coframe
$(s^*\theta^1,s^*\theta^2,s^*\theta^4)$ and $\ad v^{-1}(x)$
transforms it into other compatible connection, according to
\eqref{e.c.Aco}.
\end{proof}
\end{proposition}

\section{Six-dimensional conformal geometry in the split
signature}\label{s.c6d} \noindent Until now we have not proposed
any geometric structure, apart from $\wh{\omega}^c$, that could be
associated with an ODE of generic type. Motivated by S.-S. Chern's
construction we would like to build some kind of conformal
geometry starting from an arbitrary ODE, which does not
necessarily satisfy the W\"unschmann condition. 

Let us define the `inverse' of the symmetric tensor field $
\wh{g}=2\theta^1\theta^3-(\theta^2)^2$ of section \ref{s.c.conf}
to be $\wh{g}_{inv}=\wh{g}^{ij}X_i\otimes X_j=2X_1X_3-(X_2)^2$. We
take the $\o(2,1)$-part of the connection $\wh{\omega}^c$ \ben
\Gamma= \bma
\vc{3}-\vc{1} & -\hc{4} & 0  \\
-\vc{2} & 0 & -\hc{4} \\
0 &-\vc{2} & \vc{1}-\vc{3} \ema, \een and the Levi-Civita symbol
$\epsilon_{ijk}$ in three dimensions. Next we define a new
bilinear form $\wh{\gg}$ on $\P^c$ \ben \wh{\gg}
=\epsilon_{ijk}\,\wh{g}^{kl}\,\theta^i\,\Gamma^j_{~l}. \een The
above method of obtaining  $\wh{\gg}$ of the split degenerate
signature from $\wh{g}$ is called the Sparling procedure
\cite{Nur1}. The new metric reads \be\label{e.c.g33}
\wh{\gg}=2(\Omega_1-\Omega_3)\theta^2-2\Omega_2\theta^1+2\theta^4\theta^3
\ee and was given in \cite{Nur1} in a slightly different context
for the first time. We easily find that its degenerated directions
$X_8,X_9,X_{10},$
 and $X_5+X_7$ form an integrable distribution, so that one can consider
 the
 six-dimensional space $\M^6$ of its integral leaves. The degenerated
 directions $X_8,X_9,$ and $X_{10}$ are isometries
\ben L_{X_6}\wh{\gg}=L_{X_8}\wh{\gg}=L_{X_9}\wh{\gg}=0,
 \een
whereas the fourth direction, $X_5+X_7$, is a conformal
transformation  \ben L_{(X_5+X_7)}\wh{\gg}=\wh{\gg}.\een This
allows us to project $\wh{\gg}$ to the split signature conformal
metric $[\gg]$ on $\M^6$ without any assumptions about the
underlying ODE.

It is interesting to study the normal conformal connection
associated to this geometry. Since $\P^c$ is a subbundle of the
conformal bundle over $M^6$, we can calculate the $\o(4,4)$-valued
normal conformal connection \eqref{e.cnc} at once on $\P^c$. It is
as follows.
\begin{equation*}
\wh{\mathbf{w}}=\bma \tfrac12\Omega_1&0&0&\tfrac12\Omega_2&-\tfrac12\Omega_4&-\tfrac12\Omega_6&0&0\\\\  
\Omega_1-\Omega_{{3}}&\Omega_3-\tfrac12\Omega_1&\tfrac12\theta_4&\tfrac12\Omega_2&0&-\mathrm{w}^3_{~5}
&\Omega_5&-\tfrac12\Omega_4\\\\ 
-\Omega_2&\Omega_2&\tfrac12\Omega_1&\mathrm{w}^3_{~4}&\mathrm{w}^3_{~5}&0&\Omega_4&-\tfrac12\Omega_6\\\\ 
\theta_4&0&0&\Omega_3-\tfrac12\Omega_1&-\Omega_5&-\Omega_4&0&0\\\\ 
\theta_2&0&0&\theta_3&\tfrac12\Omega_1-\Omega_3&-\Omega_2&0&0\\\\ 
\theta_1&0&0&\tfrac12\theta_2&-\tfrac12\theta_4&-\tfrac12\Omega_1&0&0\\\\ 
\theta_3&-\theta_3&-\tfrac12\theta_2&0&-\tfrac12\Omega_2&-\mathrm{w}^3_{~4}&\tfrac12\Omega_1-\Omega_3&
\tfrac12\Omega_{{2}}\\\\
0&\theta_2&\theta_1&\theta_3&\Omega_1-\Omega_3&-\Omega_2&\theta_4&-\tfrac12\Omega_1
\ema,
\end{equation*}
where
\begin{align*}
\mathrm{w}^3_{~4}=&\inc{A}{4}\theta^1+\inc{A}{2}\theta^2+\inc{A}{1}\theta^4,\\
\mathrm{w}^3_{~5}=&\tfrac12\Omega_6+\inc{A}{3}\theta^1+\inc{A}{5}\theta^2+\inc{A}{2}\theta^4.
\end{align*}

It appears that this connection is of very special form. We show
that the algebra of its holonomy group is reduced to
$\o(3,2)\semi{.}\real^5$. Let us write down the connection as
\begin{align*}
\wh{\mathbf{w}}=&(\Omega_1-\Omega_3)e_1-\Omega_2e_2+\theta^4e_3+\theta^2e_4+\theta^1e_5+\theta^3e_6+ \\
  &+\Omega_1e_7+\Omega_4e_8+\Omega_5e_9+\Omega_6e_{10}+\mathrm{w}^3_{~5}e_{11}+\mathrm{w}^3_{~4}e_{12},
\end{align*}
where $e_1,\ldots,e_{12}$ are appropriate matrices in $\o(4,4)$.
The space \ben V=\,<e_1,\ldots,e_{12}>\,\subset\o(4,4) \een is not
closed under the commutation relations, hence $V$ is not a Lie
subalgebra. However, if we extend $V$ so that it contains three
commutators $e_{13}=[e_3,e_{12}]$, $e_{14}=[e_5,e_{10}]$ and
$e_{15}=[e_5,e_{12}]$ then $<e_1,\ldots,e_{15}>$ is a Lie algebra,
a certain semidirect sum of $\o(3,2)$ and $\real^5$. Bases of the
factors are the following: \ben \real^5=<e_1+2e_7-2e_{14}, e_{11},
e_{12}, e_{13}, e_{15}>, \een \ben \o(3,2)=<e_2+e_{13}, e_3, e_4,
e_5, e_6-e_{15},e_7,e_8,e_9,e_{10},e_{14}>. \een The matrix of
$\wh{\mathbf{w}}$ can be transformed to the following conjugated
representation, which reveals its structure well \ben \bma
\tfrac{1}{2}\vc{1} & \tfrac{1}{2}\vc{2} & -\tfrac{1}{2}\vc{4} &
-\tfrac{1}{4}\vc{6} &
 2\vc{2} & -2\mathrm{w}^3_{~4} & 2\mathrm{w}^3_{~5} & 0\\\\  
 \hc{4} & \vc{3}-\tfrac{1}{2}\vc{1} & -\vc{5} &-\tfrac{1}{2}\vc{4} &
 4\vc{3}-4\vc{1} & -2\vc{2} & 0 & -2\mathrm{w}^3_{~5} \\\\ 
\hc{2} & \hc{3} & \tfrac{1}{2}\vc{1}-\vc{3} & -\tfrac{1}{2}\vc{2}
&
 4\hc{3} & 0 & 2\vc{2} & 2\mathrm{w}^3_{~4} \\\\ 
2\hc{1} & \hc{2} & -\hc{4} & -\tfrac{1}{2}\vc{1} &
 0 & -4\hc{3} & 4\vc{1}-4\vc{3} & -2\vc{2}\\\\ 
0 & 0 & 0 & 0 & \tfrac{1}{2}\vc{1} & \tfrac{1}{2}\vc{2} & \tfrac{1}{2}\vc{4} & \tfrac{1}{4}\vc{6} \\\\ 
0 & 0 & 0 & 0 & \hc{4} & \vc{3}-\tfrac{1}{2}\vc{1} & \vc{5} &\tfrac{1}{2}\vc{4} \\\\ 
0 & 0 & 0 & 0 & -\hc{2} & -\hc{3} & \tfrac{1}{2}\vc{1}-\vc{3} & -\tfrac{1}{2}\vc{2}\\\\
0 & 0 & 0 & 0 & -2\hc{1} & -\hc{2} & -\hc{4} & -\tfrac{1}{2}\vc{1}
\ema. \een $\wh{\mathbf{w}}$ has the following block structure in
this representation
\ben \wh{\mathbf{w}}=\bma \wh{\omega}^c & \wh{\tau} \\\\
        0 & -\sigma \wh{\omega}^c \sigma\ema,
\een where \ben
\sigma=\bma 0 & 0 & 0 & 1 \\
               0 & 0 & 1 & 0 \\
           0 & 1 & 0 & 0 \\
           1 & 0 & 0 & 0
       \ema.\een
Surprisingly enough the $\o(3,2)$-part of $\wh{\mathbf{w}}$, given
by the diagonal blocks, is totally determined by the $\o(3,2)$
connection $\wh{\omega}^c$. In particular, this relation holds
when $W=0$ and $\wh{\omega}^c$ is a conformal connection itself.
In this case we have rather an unexpected link between conformal
connections in dimensions three and six.
\section{Further reduction and geometry on five-dimensional bundle}\label{s.c.furred}
Theorem \ref{th.c.1} is a starting point for further reduction of
the structural group since one can use the non-constant invariants
in \eqref{e.c.dtheta_10d} to eliminate more variables $u_\mu$.
From this point of view third-order ODEs fall into three main
classes:
\begin{itemize}
\item[i)] $W=0$, $F_{qqqq}=0$,
\item[ii)] $W=0$, $F_{qqqq}\neq 0$,
\item[iii)] $W\neq 0$. 
\end{itemize}

Class i) contains the equations equivalent to $y'''=0$ and is
fully characterized by the corollary \ref{cor.c.10d_flat}. Further
reduction cannot be done due to lack of non-constant structural
functions.

Class ii) is not interesting from the geometric point of view
since it does not contain equations with five-dimensional or
larger symmetry groups, as it is proved in theorems
\ref{th.cc.sym} and \ref{th.cc.W4d} of chapter \ref{ch.class}.

Class iii), which leads to a Cartan connection on a
five-dimensional bundle, is studied below. Owing to $W\neq 0$ we
continue reduction by setting $\inc{A}{1}=1$, $\inc{A}{2}=0$,
which gives \ben
 u_1=\sqrt[3]{W}u_3,  \qquad u_5=\frac{1}{3}\frac{W_q}{\sqrt[3]{W^2}}.
\een At this moment the auxiliary variable $u_6$ which was
introduced by the prolongation becomes irrelevant and may be set
equal to zero $$u_6=0.$$ In second step we choose \ben
u_2=\frac{1}{3}Zu_3 \een and finally \ben
u_4=\frac{1}{9}\frac{W_q}{\sqrt[3]{W^2}}M-\frac{1}{3}\sqrt[3]{W}Z_q.
\een The coframe and the underlying bundle $\P^c$ of the theorem
\ref{th.c.1} have been reduced to dimension five according to the
following
\begin{theorem}[S.-S. Chern]\label{th.c.2}
A third-order ODE $y'''=F(x,y,y',y'')$ satisfying the contact
invariant condition $W\neq 0$ and considered modulo contact
transformations of variables, uniquely defines a 5-dimensional
bundle $\P^c_5$ over $\J^2$ and an invariant coframe
$(\theta^1,\ldots,\theta^4,\Omega)$ on it. In local coordinates
$(x,y,p,q,u)$ this coframe is given by \be
\label{e.c.theta_5d}\begin{aligned}
 \theta^1=&\sqrt[3]{W}u\omega^1, \\
 \theta^2=&\frac{1}{3}Z u\omega^1+u\omega^2, \\
\theta^3=&\frac{u}{\sqrt[3]{W}}\left(K+\frac{1}{18}Z^2\right)\omega^1
+\frac{u}{3\sqrt[3]{W}}\left(Z-F_q\right)\omega^2+\frac{u}{\sqrt[3]{W}}\omega^3, \\
\theta^4=&\left(\frac{1}{9}\frac{W_q}{\sqrt[3]{W^2}}Z-\frac{1}{3}\sqrt[3]{W}Z_q\right)\omega^1
+\frac{W_q}{3\sqrt[3]{W^2}}\omega^2+\sqrt[3]{W}\omega^4, \\
 \Omega=& \left(\left(\frac{1}{9}W_q\D Z-\frac{1}{27}W_qZ^2+\frac{1}{9}W_p Z\right)\frac{1}{W}
 -\frac{1}{3}Z_p-\frac{1}{9}F_qZ_q\right)\omega^1 \\
 &+\left(\frac{W_p}{3W}-\frac{1}{3}Z_q\right)\omega^2+
 \frac{W_q}{3W}\,\omega^3+\frac{1}{3}F_q\,\omega^4
 +\frac{\der u}{u}.
 \end{aligned}\ee
where  $\omega^i$ are defined by the ODE via \eqref{e.omega}. The
exterior derivatives of these forms read
\begin{align}
 d\theta^1=&\Omega\w\theta^1-\theta^2\w\theta^4, \nonumber \\
 d\theta^2=&\Omega\w\theta^2+\inc{a}{}\,\theta^1\w\theta^4-\theta^3\w\theta^4,\nonumber\\
 d\theta^3=&\Omega\w\theta^3+\inc{b}{}\,\theta^1\w\theta^2+\inc{c}{}\,\theta^1\w\theta^3
  -\theta^1\w\theta^4+\inc{e}\,\theta^2\w\theta^3+\inc{a}{}\,\theta^2\w\theta^4, \label{e.c.dtheta_5d} \\
 d\theta^4=&\inc{f}{}\,\theta^1\w\theta^2+\inc{g}{}\,\theta^1\w\theta^3+\inc{h}{}\,\theta^1\w\theta^4
  +\inc{k}{}\,\theta^2\w\theta^3-\inc{e}{}\,\theta^2\w\theta^4, \nonumber  \\
 d\Omega=&\inc{l}{}\,\theta^1\w\theta^2+(\inc{f}{}-\inc{a}{}\inc{k}{})\,\theta^1\w\theta^3
 +\inc{m}{}\,\theta^1\w\theta^4+\inc{g}{}\,\theta^2\w\theta^3+\inc{h}{}\,\theta^2\w\theta^4.\nonumber
 \end{align}
\end{theorem}
The basic functions for \eqref{e.c.dtheta_5d} (i.e. generating the
full set of functions by consecutive taking of coframe
derivatives) are
$\inc{a}{},\inc{b}{},\inc{e}{},\inc{h}{},\inc{k}{}$:
\begin{align}
\inc{a}{}=& \frac{1}{\sqrt[3]{W^2}}\left(K+\frac{1}{18}Z^2+\frac{1}{9}ZF_q-\frac{1}{3}\D Z\right), \nonumber \\
\inc{b}{}=&\frac{1}{3u\sqrt[3]{W^2}}\bigg(\frac{1}{27}F_{qq}Z^2+\left(K_q-\frac{1}{3}Z_p-\frac{2}{9}F_qZ_q\right)Z+\nonumber\\
   &+\left(\frac{1}{3}\D Z-2K\right)Z_q+Z_y+F_{qq}K-3K_p-K_qF_q-F_{qy}+W_q\bigg), \nonumber \\
\inc{e}{}=&\frac{1}{u}\left(\frac{1}{3}F_{qq}+\frac{1}{W}\left(\frac{2}{9}W_qZ-\frac{2}{3}W_p-\frac{2}{9}W_qF_q \right)\right), \label{e.c.basfun_5d}\\
\inc{h}{}=&\frac{1}{3u\sqrt[3]{W}}\bigg(\left(\frac{1}{9}W_qZ^2-\frac{1}{3}W_pZ+W_y-\frac{1}{3}W_q\D Z\right)\frac{1}{W}+\nonumber  \\
 &+\D Z_q+\frac{1}{3}F_qZ_q\bigg),\nonumber \\
\inc{k}{}=&\frac{1}{u^2\sqrt[3]{W}}\left(\frac{2W_q^2}{9W}-\frac{W_{qq}}{3}\right).
\nonumber
\end{align}

Our next aim is to obtain a Cartan connection. First of all we
study the most symmetric case to find the Lie algebra of a
connection. We assume that all the functions
$\inc{a}{},\ldots,\inc{m}{}$ are constant. Having applied the
exterior derivative to \eqref{e.c.dtheta_5d} we get that
$\inc{b}{},\ldots,\inc{m}{}$ are equal to zero and $\inc{a}{}$ is
an arbitrary real constant $\mu$. In this case the equations
\eqref{e.c.dtheta_5d} become the Maurer-Cartan equations for the
algebra $\real^2\oplus_\mu\real^3$. Straightforward calculations
show that this case corresponds to a general linear equation with
constant coefficients.
\begin{corollary} \label{cor.c.5d_flat}
A third-order ODE is contact equivalent to
$$y'''=-2\mu y'+y,$$ where $\mu$ is an arbitrary constant if and only if
 it satisfies
\begin{align}
1)\quad &W\neq 0, \nonumber \\
2)\quad &\frac{1}{\sqrt[3]{W^2}}
\left(K+\frac{1}{18}Z^2+\frac{1}{9}ZF_q-\frac{1}{3}\D Z\right)=\mu \label{e400} \\
3)\quad &2W_q^2-3W_{qq}W=0. \nonumber
\end{align}
Such an equation has the five-dimensional algebra
$\real^2\semi{\mu}\real^3$ of infinitesimal contact symmetries.
The equations with different constants $\mu_1$ and $\mu_2$ are
non-equivalent.
\begin{proof}
Assume that $\inc{a}{}=\mu$, $\inc{k}{}=0$. It follows from
$\der^2\theta^i=0$ and $\der^2\Omega=0$, that this assumption
makes other functions in \eqref{e.c.dtheta_5d} vanish. Put
$y'''=-2\mu y'+y$ into the formulae of theorem \ref{th.c.2} and
check that it satisfies $\inc{a}{}=\mu$, $\inc{k}{}=0$. Every
equation satisfying $\inc{a}{}=\mu$, $\inc{k}{}=0$ is contact
equivalent to it by virtue of theorem \ref{th.i.clas} adapted to
this situation.
\end{proof}
\end{corollary}
Now we immediately find a family of Cartan connections
\begin{theorem}\label{cor.c.geom5d}
An ODE which satisfies the condition \ben
\frac{1}{\sqrt[3]{W^2}}\left(K+\frac{1}{18}Z^2+\frac{1}{9}ZF_q-\frac{1}{3}\D
Z\right)=\lambda\een has the solution space equipped with the
following $\real^2$-valued linear torsion-free connection
\ben
 \widehat{\omega}_\lambda=\bma
             -\vc{} & -\hc{4} & 0 \\
             \lambda\hc{4} & -\vc{} & -\hc{4} \\
             -\hc{4} & \lambda\hc{4} & -\vc{}
        \ema.
\een Its curvature reads
 \ben
   \bma
      R^1_{~1} & R^1_{~2} & 0 \\
      -\lambda R^1_{~2} & R^1_{~1} & R^1_{~2} \\
      R^1_{~2} & -\lambda R^1_{~2} & R^1_{~1}
   \ema
 \een
 with
\begin{align}
 R^1_{~1}=&(\lambda\inc{g}{}+\inc{k}{})\hc{1}\w\hc{2}+(\lambda \inc{k}{}-\inc{g}{})\hc{1}\w\hc{3}
 -\inc{g}{}\hc{2}\w\hc{3}, \nonumber \\
 R^1_{~2}=&-\inc{f}{}\hc{1}\w\hc{2}-\inc{g}{}\hc{1}\w\hc{3}-\inc{k}{}\hc{2}\w\hc{3}. \nonumber
\end{align}
The connection is flat if and only if the related ODE is contact
equivalent to $y'''=-2\lambda y'+y$
\begin{proof}
The condition $\inc{a}{}=\lambda$ together with its differential
consequences $\inc{b}{}=\inc{c}{}=\inc{e}{}=\inc{h}{}=\inc{m}{}=0$
and $\inc{l}{}=-\inc{k}{}-\lambda \inc{g}{}$ is the necessary and
sufficient condition for the curvature of
$\widehat{\omega}_\lambda$ to be horizontal.
\end{proof}
\end{theorem}

It follows  that every ODE as above has its space of solution equipped with a geometric structure consisting of
\begin{itemize}
\item[i)] Reduction of $\gl(3,\real)$ to $\real^2$ represented by
\ben
   \bma
      a_1 & a_2 & 0 \\
      -\lambda a_2 & a_1 & a_2 \\
      a_2 & -\lambda a_2 & a_1
   \ema.
\een \item[ii)] A linear torsion-free connection $\Gamma$ taking
values in this $\real^2$.
\end{itemize}
The structure is an example of a geometry with special holonomy.
The algebra $\real^2$ is spanned by the unit matrix and \ben
 \mathrm{m}(\lambda) = \bma
      0 & 1 & 0 \\
      -\lambda & 0 & 1 \\
      1 & -\lambda  & 0
   \ema,
\een whose action on $\S$ is more complicated. Its eigenvalue
equation
$$ \det(u\mathbf{1}-\mathrm{m}(\lambda))=u^3+2\lambda u-1
$$ is the characteristic polynomial of the linear ODE  $y'''=-2\lambda y'+y$. If
$\lambda<\tfrac{3}{4}\sqrt[3]{2}$ the polynomial has three
distinct roots and $\mathrm{m}$ is a generator of non-isotropic
dilatations acting along the eigenspaces. If
$\lambda=\tfrac{3}{4}\sqrt[3]{2}$ the characteristic polynomial
has two roots, one of them double, for other $\lambda$s there is
one eigenvalue. The action is diagonalizable only in the case of
three distinct eigenvalues.

\chapter[Geometries of ODEs modulo point and fibre-preserving\ldots]
{Geometries of ODEs considered modulo point and fibre-preserving
transformations of variables}\label{ch.point}

\section{Point case: Cartan connection on seven-dimensional
bundle}\label{s.p.th}  \noindent Following the scheme of reduction
given in chapter \ref{ch.contact} we construct Cartan connection
for ODEs modulo point transformations.

\begin{theorem}[E. Cartan]\label{th.p.1}
The point invariant information about $y'''=F(x,y,y',y'')$ is
given by the following data
\begin{itemize}
\item[i)] The principal fibre bundle $H_3\to\P^p\to\J^2$, where
$\dim\P^p=7$, and $H_3$ is the three-dimensional group \be
\label{e.p.H3}  H_3=\bma \sqrt{u_1}, &
\frac12\frac{u_2}{\sqrt{u_1}}, & 0 & 0 \\\\
 0 & \tfrac{u_3}{\sqrt{u_1}}, & 0 & 0 \\\\
 0 & 0 & \tfrac{\sqrt{u_1}}{u_3}, &
-\tfrac12\tfrac{u_2}{\sqrt{u_1}\,\u_3} \\\\
  0 & 0 & 0 & \tfrac{1}{\sqrt{u_1}}\ema.
\ee \item[ii)] The coframe
$(\theta^1,\theta^2,\theta^3,\theta^4,\Omega_1,\Omega_2,\Omega_3)$,
which defines the $\co(2,1)\oplus_{.}\real^3$-valued Cartan normal
connection $\wh{\omega}^p$ on $\P^p_{7}$ by \be\label{e.p.conn_7d}
 \wh{\omega}^p=\bma \tfrac{1}{2}\vp{1} & \tfrac{1}{2}\vp{2} & 0 & 0 \\\\
             \hp{4} & \vp{3}-\tfrac{1}{2}\vp{1} & 0 & 0 \\\\
             \hp{2} & \hp{3} & \tfrac{1}{2}\vp{1}-\vp{3} & -\tfrac{1}{2}\vp{2} \\\\
             2\hp{1} & \hp{2} & -\hp{4} & -\tfrac{1}{2}\vp{1}
        \ema.
\ee
\end{itemize}

Let $(x,y,p,q,u_1,u_2,u_3)=(x^i,u_\mu)$ be a locally trivializing
coordinate system in $\P^p$. Then the value of $\wh{\omega}^p$ at
the point $(x^i,u_\mu)$ in $\P^p$ is given by \ben
\wh{\omega}^p(x^i,u_\mu)=u^{-1}\,\omega^p\,u+u^{-1}\der u \een
where $u$ denotes the matrix \eqref{e.p.H3} and \ben\omega^p= \bma
\tfrac{1}{2}\vp{1}^0&\tfrac{1}{2}\vp{2}^0 & 0 & 0 \\\\
\wt{\omega}^4 & \vp{3}^0-\tfrac{1}{2}\vp{1}^0 & 0 & 0 \\\\
\omega^2 & \wt{\omega}^3 & \tfrac{1}{2}\vp{1}^0-\vp{3}^0 & -\tfrac{1}{2}\vp{2}^0 \\\\
2\omega^1 & \omega^2 & -\wt{\omega}^4 & -\tfrac{1}{2}\vp{1}^0
\ema\een is the connection $\wh{\omega}^p$ calculated at the point
$(x^i,u_1=1,u_2=0,u_3=1)$.
 The forms $\omega^1,\omega^2,\wt{\omega}^3,\omega^4$ read
\begin{align*}
\omega^1=&\der y-p\der x, \notag \\
\omega^2=&\der p-q\der x, \\
\wt{\omega}^3=&\der q-F\der x-\tfrac13F_q(\der p-q\der x)+K(\der y-p\der x),\notag \\
\wt{\omega}^4=&\der x +\tfrac16F_{qq}(\der y-p\der x).\notag
\end{align*}
The forms $\vc{1}^0,\ldots,\vc{6}^0$ read
\begin{align*} 
\vc{1}^0=&-(3K_q+\tfrac29F_{qq}F_q+\tfrac23F_{qp})\,\omega^1+\tfrac16F_{qq}\omega^2, \notag \\
\vc{2}^0=&\left(L+\tfrac16F_{qq}K\right)\,\omega^1
 -(2K_q+\tfrac19F_{qq}F_q+\tfrac13F_{qp})\,\omega^2
 +\tfrac16F_{qq}\wt{\omega}^3-K\wt{\omega}^4,  \\
\vc{3}^0=&-(2K_q+\tfrac16F_{qq}F_q+\tfrac13F_{qp})\,\omega^1+\tfrac13F_{qq}\,\omega^2
+\tfrac13F_q\wt{\omega}^4.\notag
\end{align*}
\begin{proof}
We begin with the $G_p$-structure on $\J^2$ of Introduction, which
encodes an ODE up to point transformations. In the usual locally
trivializing coordinate system $(x,y,p,q,u_1,\ldots,u_8)$ on
$G_p\times\J^2$ the fundamental form $\hp{i}$ is given by
\begin{align*}
 \hp{1}=& u_1\omega^1,\\
 \hp{2}=& u_2\omega^2+u_3\omega^3, \\
 \hp{3}=& u_4\omega^1+u_5\omega^2+u_6\omega^3,\\
 \hp{4}=& u_8\omega^1+u_7\omega^4.
 \end{align*}
We repeat the procedure of section \ref{s.c.proof} of chapter
\ref{ch.contact}. We choose a connection by the minimal torsion
requirement and then reduce $G_p\times\J^2$ using the constant
torsion property. We differentiate $\hp{i}$ and gather the
$\hp{j}\w\hp{k}$ terms into \be\label{e.p.red10}\begin{aligned}
\der\hp{1}&=\vp{1}\w\hp{1}+\frac{u_1}{u_3 u_7}\hp{4}\w\hp{2}, \\
 \der\hp{2}&=\vp{2}\w\hp{1}+\vp{3}\w\hp{2}+\frac{u_3}{u_6u_7}\hp{4}\w\hp{3}, \\
 \der\hp{3}&=\vp{4}\w\hp{1}+\vp{5}\w\hp{2}+\vp{6}\w\hp{3},\\
 \der\hp{4}&=\vp{8}\w\hp{1}+\vp{9}\w\hp{2}+\vp{7}\w\hp{4}
\end{aligned}\ee
with the auxiliary connection forms $\vp{\mu}$ containing the
differentials of $u_\mu$ and terms proportional to $\hp{i}$. Then
we reduce $G_p\times\J^2$ by
\ben
  u_6=\frac{u_3^2}{u_1},\quad\quad\quad u_7=\frac{u_1}{u_3}.
\een  Subsequently, we get formulae identical to
\eqref{e.c.red_u5}, \eqref{e.c.red_u4}:
\begin{align*}
 u_5=&\frac{u_3}{u_1}\left(u_2-\frac{1}{3}u_3F_q\right), \\
 u_4=&\frac{u^2_3}{u_1}K+\frac{u_2^2}{2u_1}
\end{align*}
and also
\ben
 u_8=\frac{u_1}{6u_6}F_{qq}.
\een  After these substitutions the structural equations for
$\hp{i}$ are the following
\begin{align}\label{e.p.red50}
 \der\hp{1} =&\vp{1}\w\hp{1}+\hp{4}\w\hp{2},\nonumber \\
 \der\hp{2} =&\vp{2}\w\hp{1}+\vp{3}\w\hp{2}+\hp{4}\w\hp{3},\nonumber \\
 \der\hp{3} =&\vp{2}\w\hp{2}+(2\vp{3}-\vp{1})\w\hp{3}+\inp{A}{1}\hp{4}\w\hp{1},\nonumber \\
 \der\hp{4} =&(\vp{1}-\vp{3})\w\hp{4}+\inp{B}{1}\hp{2}\w\hp{1}
   +\inp{B}{2}\hp{3}\w\hp{1}, \nonumber
 \end{align}
with some functions $\inp{A}{1},\inp{B}{1},\inp{B}{2}$. But now,
in contrast to the contact case, the forms $\vp{1},\vp{2},\vp{3}$
are defined by the above equations without any ambiguity, thus
there is no need to prolong and we have the rigid coframe on the
seven-dimensional bundle $\P^p\to\J^2$.
\end{proof}
\end{theorem}

\subsection{Point versus contact objects} Comparing  theorem
\ref{th.p.1} to theorem \ref{th.c.1} of the contact case, we see
that the contact and point objects are related as follows. The
point bundle $\P^p$ is a subbundle of $\P^c$, the embedding
$\sigma\colon\P^p\to\P^c$ given by \be\label{e.p.cp}
u_4=\frac{u_1}{6u_3}F_{qq}, \qquad u_5=0, \qquad u_6=0. \ee In
this paragraph we denote the point coframe of theorem \ref{th.p.1}
by $(\overset{p}{\hp{1}},\ldots,\overset{p}{\hp{3}})$ in order to
distinguish it from the contact coframe now denoted by
$(\overset{c}{\hc{1}},\ldots,\overset{c}{\vc{7}})$. The point
coframe on $\P^p$ can be constructed from the contact one by the
following formula. \be\label{e.p.ctop2}\begin{aligned}
&\overset{p}{\hp{i}}=\sigma^*\overset{c}{\hc{i}}, \qquad i=1,2,3,4, \\
&\overset{p}{\vp{1}}=\sigma^*\overset{c}{\vc{1}}-2f_1\,\sigma^*\overset{c}{\hc{1}}, \\
&\overset{p}{\vp{2}}=\sigma^*\overset{c}{\vc{2}}-f_2\,\sigma^*\overset{c}{\hc{1}}
 -f_1\,\sigma^*\overset{c}{\hc{2}}, \\
&\overset{p}{\vp{3}}=\sigma^*\overset{c}{\vc{3}}-f_1\,\sigma^*\overset{c}{\hc{1}},
\end{aligned}\ee
where the functions $f_1$ and $f_2$ are
\begin{align*}
f_1=&\frac{1}{u_1}(K_q+\tfrac19F_{qq}F_q+\tfrac13F_{qp})+\frac{u_2}{12\,u_1u_3}F_{qq},
\notag \\
f_2=&\sigma^*\inc{A}{2}=\frac{u_3}{3u^2_1} W_q\notag.
\end{align*}
The mapping $\sigma$ together with the functions $f_1,f_2$ is a
new piece of structure that allows us to move from the more coarse
contact coframe to the point coframe. It is obvious that $F_{qq}$
and $f_1$ appearing in the above formulae are not contact
invariants. What is more, they are not point invariants either,
for the property $F_{qq}=0$ or $f_1=0$  is not preserved under
point transformations. At the level of geometric objects the
passage from the contact case to the point case is given by the
difference
$$\sigma^*\wh{\omega}^c-\wh{\omega}^p.$$

One can also proceed in the inverse direction and construct the
contact bundle and coframe starting from the point case. In this
approach the bundle $\P^c$ is the extension
$\P^c=\P^p\times_{H_3}H_6$ of $\P^p$ and in order to get the
contact coframe we must define forms
$\overset{c}{\hc{1}},\ldots,\overset{c}{\vc{6}}$ on $\P^p$ in
terms of the point coframe, then define a matrix one-form
$\wh{\omega}^c$ on $P^p_7$ and finally lift this $\wh{\omega}^c$
to the Cartan connection on $\P^c$. Since the formulae for
$\overset{c}{\hc{1}},\ldots,\overset{c}{\vc{3}}$ are similar to
\eqref{e.p.ctop2} and the formulae for
$\overset{c}{\vc{4}},\overset{c}{\vc{5}},\overset{c}{\vc{6}}$ are
complicated we omit them.

\subsection{Curvature} Further analysis of the coframe of
theorem \ref{th.p.1} is very similar to what we have done in
chapter \ref{ch.contact}. The curvature of the connection is given
by the following exterior differentials of the coframe, cf
\eqref{e.i.dtheta_7d}
\begin{align}
 \der\hp{1} =&\vp{1}\w\hp{1}+\hp{4}\w\hp{2},\nonumber \\
 \der\hp{2} =&\vp{2}\w\hp{1}+\vp{3}\w\hp{2}+\hp{4}\w\hp{3},\nonumber \\
 \der\hp{3}=&\vp{2}\w\hp{2}+(2\vp{3}-\vp{1})\w\hp{3}+\inp{A}{1}\hp{4}\w\hp{1},\nonumber \\
 \der\hp{4} =&(\vp{1}-\vp{3})\w\hp{4}+\inp{B}{1}\hp{2}\w\hp{1}
   +\inp{B}{2}\hp{3}\w\hp{1}, \nonumber\\
 \der\vp{1} =&-\vp{2}\w\hp{4}+(\inp{D}{1}+3\inp{B}{3})\hp{1}\w\hp{2}
   +(3\inp{B}{4}-2\inp{B}{1})\hp{1}\w\hp{3} \label{e.p.dtheta_7d}  \\
   &+(2\inp{C}{1}-\inp{A}{2})\hp{1}\w\hp{4}-\inp{B}{2}\hp{2}\w\hp{3}, \nonumber \\
 \der\vp{2}=&(\vp{3}-\vp{1})\w\vp{2}+\inp{D}{2}\hp{1}\w\hp{2}+(\inp{D}{1}+\inp{B}{3})\hp{1}\w\hp{3}
 +\inp{A}{3}\hp{1}\w\hp{4} \nonumber \\
  &+(2\inp{B}{4}-\inp{B}{1})\hp{2}\w\hp{3}+\inp{C}{1}\hp{2}\w\hp{4}, \nonumber \\
 \der\vp{3}=&(\inp{D}{1}+2\inp{B}{3})\hp{1}\w\hp{2}+2(\inp{B}{4}-\inp{B}{1})\hp{1}\w\hp{3}
  +\inp{C}{1}\hp{1}\w\hp{4}-2\inp{B}{2}\hp{2}\w\hp{3}, \nonumber
\end{align}
where
$\inp{A}{1},\inp{A}{2},\inp{A}{3},\inp{B}{1},\inp{B}{2},\inp{B}{3},
\inp{B}{4},\inp{C}{1},\inp{D}{1},\inp{D}{2}$ are functions on
$\P^p$. All these functions express by the coframe derivatives of
$\inp{A}{1},\inp{B}{1},\inp{C}{1}$, which therefore constitute the
set of basic relative invariants for this problem and read
\begin{align*}
 \inp{A}{1}=&\frac{u_3^3}{u^3_1}W, \nonumber \\ \inp{B}{1}=&\frac{1}{u_3^2}\left(\frac{1}{18}F_{qqq}F_q+\frac{1}{36}F_{qq}^2+\frac{1}{6}F_{qqp}\right)
 -\frac{u_2}{6u_3^3}F_{qqq}, \\
 \inp{C}{1}=&\frac{u_3}{u_1^2}\left(2F_{qq}K+\frac{2}{3}F_qF_{qp}-2F_{qy}+F_{pp} +2W_q\right). \nonumber
\end{align*}
$\inp{B}{2}$ and $\inp{B}{4}$ are also vital:
\begin{align*}
  \inp{B}{2}=&\frac{u_1}{6u_3^3}F_{qqq}, \\
  \inp{B}{4}=&\frac{1}{u_3^2}
  \left(K_{qq}+\frac19F_{qqq}F_q+\frac13 F_{qqp}+\frac{1}{12}F_{qq}^2 \right).
\end{align*}
\begin{corollary}\label{cor.p.7d_flat}
For a third-order ODE $y'''=F(x,y,y',y'')$ the following
conditions are equivalent.
\begin{itemize}
 \item[i)] The ODE is point equivalent to $y'''=0$.
 \item[ii)] It satisfies the conditions $W=0$, $F_{qqq}=0$, $F_{qq}^2+6F_{qqp}=0$
 and $$2F_{qq}K+\frac{2}{3}F_qF_{qp}-2F_{qy}+F_{pp}=0.$$
  \item[iii)] It has the $\co(2,1)\semi{.}\real^3$ algebra of infinitesimal point symmetries.
\end{itemize}
\end{corollary}


The manifold $\P^p$, like $\P^c$, is equipped with the threefold
structure of principal bundle over $\J^2$, $\J^1$ and $\S$. Let
$(X_1,X_2,X_3,X_4,X_5,X_6,X_7)$ be the frame dual to
$(\hp{1},\hp{2},\hp{3},\hp{4},\vp{1},\vp{2},\vp{3})$.
\begin{itemize}
\item $\P^p$ is the bundle $H_3\to\P^p\to\J^2$ with the
fundamental fields $X_5,X_6,X_7$. \item It is the bundle
$CO(2,1)\to\P^p\to\S$ with the fundamental fields
$X_4,X_5,X_6,X_7$ \item It is also the bundle $H_4\to\P^p\to\J^1$
with the fundamental fields $X_3,X_5,X_6,X_7$.
\end{itemize}

\section{Einstein-Weyl geometry on space of solutions}\label{s.ew}
\noindent We have already described the construction of the
Einstein-Weyl geometry on the solution space in Introduction. Here
we write it down in a more systematic manner.

\subsection{Weyl geometry}\label{s.ew.def}
A Weyl geometry on $\M^n$ is a pair $(g,\phi)$ such that $g$ is a
metric of signature $(k,l)$, $k+l=n$ and $\phi$ is a one-form and
they are given modulo the following transformations \ben
\phi\to\phi+\der\lambda, \qquad\qquad g\to e^{2\lambda} g.\een In
particular $[g]$ is a conformal geometry. For any Weyl geometry
there exists the Weyl connection; it is the unique torsion-free
connection such that \ben
    \nabla g=2\phi\otimes g.
\een The Weyl connection takes values in the algebra $\co(k,l)$ of
$[g]$. Let $(\omega^\mu)$ be an orthonormal coframe for some $g$
of $[g]$; $g=g_{\mu\nu}\omega^\mu\otimes\omega^\nu$ with all the
coefficients $g_{\mu\nu}$ being constant. The Weyl connection
one-forms $\Gamma^\mu_{~\nu}$ are uniquely defined by the
relations
\begin{eqnarray}
 &\der\omega^\mu+\Gamma^\mu_{~\nu}\w\omega^\nu=0, \notag \\
 &\Gamma_{(\mu\nu)}=-g_{\mu\nu}\,\phi, \quad \text{where} \quad
 \Gamma_{ij}=g_{jk}\Gamma^k_{~j}. \notag
\end{eqnarray}
The curvature tensor $R^\mu_{~\nu\rho\sigma}$, the Ricci tensor
$\Ric_{\mu\nu}$ and the Ricci scalar $\R$ of a Weyl connection are
defined as follows
\begin{align*} &R^\mu_{~\nu}=\der\Gamma^\mu_{~\nu}+\Gamma^\mu_{~\rho}\w\Gamma^{\rho}_{~\nu}=
R^\mu_{~\nu\rho\sigma}\omega^\rho\w\omega^\sigma \\
&\Ric_{\mu\nu}=R^\rho_{~\mu\rho\nu}, \\
&\R=\Ric_{\mu\nu}g^{\mu\nu}.
\end{align*}
The Ricci scalar has the conformal weight $-2$, that is it
transforms as $\R\to e^{-2\lambda}\R$ when $g\to e^{2\lambda}g$.
Apart from these objects there is another one, which has no
counterpart in the Riemannian geometry, the Maxwell two-form \ben
F=\der \phi. \een The Maxwell two-form $F_{\mu\nu}$ is
proportional to the antisymmetric part $\Ric_{[\mu\nu]}$ of the
Ricci tensor.

Einstein-Weyl  structures are, by definition, those Weyl
structures for which the symmetric trace-free part of the Ricci
tensor vanishes \ben
 \Ric_{(\mu\nu)}-\tfrac{1}{n}R\cdot g_{\mu\nu} =0.
\een

\subsection{Einstein-Weyl structures from ODEs} In this paragraph
we follow P. Nurowski \cite{Nur1,Nur2}. One sees from the system
\eqref{e.p.dtheta_7d}  that the pair $(\wh{g},\vp{3})$, where
\ben\wh{g}=2\hp{1}\hp{3}-(\hp{2})^2\een is Lie transported along
the fibres of $\P^p\to\S$ in the following way \ben L_{X_4}g=
\inp{A}{1}(\theta^1)^2,\qquad L_{X_5}\wh{g}=0,\qquad
L_{X_6}\wh{g}=0,\qquad L_{X_7}\wh{g}=2\wh{g},\een and \ben
L_{X_3}\Omega_3=\tfrac12\inp{C}{1}\hp{1}, \een
\ben\label{e.p.lie5} L_{X_j}\Omega_3=0,\qquad\text{for}\qquad
j=5,6,7. \een

\noindent Due to these properties  $(\wh{g},\Omega_3)$ descends
along $\P^p\to\S$ to the Weyl structure $(g,\phi)$ on the solution
space $\S$ on condition that \ben W=0 \een and \be\label{e.p.cart}
\left(\tfrac{1}{3}\D F_q
-\tfrac{2}{9}F_q^2-F_p\right)F_{qq}+\tfrac{2}{3}F_qF_{qp}-2F_{qy}+F_{pp}=0.
\ee These conditions are equivalent to E. Cartan's original
conditions \eqref{e.i.wunsch}, \eqref{e.i.cart}. Our conditions
are even simpler, because the quantity \eqref{e.p.cart} is of
first order in $\D$.

The conformal metric of the Weyl structure $(g,\phi)$ coincides
with the conformal metric of the contact case and is represented
by
\begin{align*}
g&=2\omega^1\wt{\omega}^3-(\omega^2)^2= \\
&=2(\der y-p\der x)(\der q -\tfrac13F_q\der p+ K\der y
+(\tfrac13qF_q-pK-F)\der x)-(\der p-q\der x)^2,
\end{align*} while the Weyl potential is given by
\ben\phi=-(2K_q+\tfrac19F_{qq}F_q+\tfrac13F_{qp})(\der y-p\der
x)+\tfrac13F_{qq}(\der p-q\der x) +\tfrac13F_q\der x \een

The Weyl connection for this geometry, lifted to
$CO(2,1)\to\P^p\to\S$, now the bundle of orthonormal frames, reads
\ben
 \Gamma=\bma
             -\vp{1} & -\hp{4} & 0 \\
             -\vp{2} & -\vp{3} & -\hp{4} \\
              0 &-\vp{2} & \vp{1}-2\vp{3}
        \ema.
\een The curvature  is as follows \be\label{e.ewcur}
 (R^\mu_{~\nu})=\bma
R^1_{~1}-F&R^1_{~2}&0\\
R^2_{~1}&-F&R^1_{~2}\\
0&R^2_{~1}&-R^1_{~1}-F \ema \ee with
\begin{eqnarray}
&F=\der\vp{3}=2\inp{B}{3}\hp{1}\w\hp{2}+(2\inp{B}{4}-2\inp{B}{1})\hp{1}\w\hp{3}
-2\inp{B}{2}\hp{2}\w\hp{3},\nonumber  \\
&R^1_{~1}=-\inp{B}{3}\hp{1}\w\hp{2}-\inp{B}{4}\hp{1}\w\hp{3}-\inp{B}{2}\hp{2}\w\hp{3},\nonumber \\
&R^1_{~2}=\inp{B}{1}\hp{1}\w\hp{2}+\inp{B}{2}\hp{1}\w\hp{3},\nonumber \\
&R^2_{~1}=-\inp{B}{3}\hp{1}\w\hp{3}+(\inp{B}{1}-2\inp{B}{4})\hp{2}\w\hp{3}.
\nonumber
\end{eqnarray}
The Ricci tensor reads
\ben
\Ric=\begin{pmatrix}
0&-3\inp{B}{3}&3\inp{B}{1}-5\inp{B}{4}\\
3\inp{B}{3}&2\inp{B}{4}&3\inp{B}{2}\\
-3\inp{B}{1}+\inp{B}{1}&-3\inp{B}{2}&0
\end{pmatrix}
\een and satisfies the Einstein-Weyl equations \ben
 \Ric_{(ij)}=\tfrac{1}{3}R \cdot {g}_{ij}
\een with the Ricci scalar $R=6\inp{B}{4}$. The components of the
curvature in the orthogonal coframe given by $u_1=1, u_2=0, u_3=1$
are given by
\begin{align*}
\inp{B}{1}=&\tfrac{1}{18}F_{qqq}F_q+\tfrac16F_{qqp}+\tfrac{1}{36}F_{qq}^2,
\\
\inp{B}{2}=&\tfrac16F_{qqq}, \\
\inp{B}{3}=&\tfrac16F_{qqy}-\tfrac13F_{qq}K_q-\tfrac16F_{qqq}K-\tfrac{1}{18}F_{qq}F_{qp}
-\tfrac{1}{54}F^2_{qq}F_q-L_q,\\
\inp{B}{4}=&K_{qq}+\frac19F_{qqq}F_q+\tfrac13F_{qqp}+\frac{1}{12}F_{qq}^2.
\end{align*}


\section{Geometry on first jet space}\label{s.p-p}
\noindent In section \ref{s.c-p} of chapter \ref{ch.contact} we
described how certain ODEs modulo contact transformations generate
contact projective geometry on $J^1$. The fact that point
transformations form a subclass within contact transformations
suggests that ODEs modulo point transformations define some
refined version of contact projective geometry. Indeed, the only
object that is preserved by point transformations but is not
preserved by contact transformations is the projection $\J^1\to
\text{\em xy plane}$, whose fibres are generated by $\partial_p$.
This motivate us to propose the following
\begin{definition}\label{def.p-p}
A point projective structure on $\J^1$ is a contact projective
structure, such that integral curves of the field $\partial_p$ are
geodesics of the contact projective structure.
\end{definition}
We immediately get
\begin{lemma}
The field $\partial_p$ is geodesic for the contact projective
geometry generated by an ODE provided that the ODEs satisfies
$$F_{qqq}=0.$$
\begin{proof}
In the notation of section \ref{s.cp3} of chapter \ref{ch.contact}
we have $\partial_p=e_2$ and, from proposition
\ref{prop.cp.coord}, $\Gamma^1_{~22}=0$. Thus $\nabla_2e_2=\lambda
e_2$ iff $\Gamma^3_{~22}=0$, which is equivalent to $F_{qqq}=0$ by
means of lemma \ref{lem.c-p}.
\end{proof}
\end{lemma}

However, the condition $F_{qqq}=0$ is not sufficient for the form
\eqref{e.p.conn_7d} to be a Cartan connection for the point
projective structure and we show that there does not exist any
simple way to construct a Cartan connection on $\P^p\to J^1$. The
algebra $\co(2,1)\semi{.}\real^3\subset\o(3,2)$ inherits the
following grading from \eqref{c.gradJ1}
$$\co(2,1)\semi{.}\real^3=\g_{-2}\oplus\g_{-1}\oplus\g_{0}\oplus\g_{1}$$
but it is not semisimple, so the Tanaka method cannot be
implemented. Moreover, the broadest generalization of this method
-- the Morimoto nilpotent geometry handling non-semisimple groups
also fails in this case. This is because the Morimoto approach
requires the algebra $\g$ to be equal to the prolongation of its
non-positive part algebra $\g_{-2}\oplus\g_{-1}\oplus\g_0$. The
notion of prolongation as well as algorithmic procedure for
calculating it has been introduced by N. Tanaka \cite{Tan2}. In
our case  the prolongation of $\g_{-2}\oplus\g_{-1}\oplus\g_0$ is
larger then $\co(2,1)\semi{.}\real^3$ and equals precisely
$\o(3,2)$, which yields a contact projective structure. Thus only
the contact case is solved by the methods of the nilpotent
geometry.

Lacking general theory we must search a Cartan connection in a
more direct way. Consider than an ODE satisfying $F_{qqq}=0$. It
follows that $\inp{B}{2}=0$ and $\inp{B}{4}=\inp{B}{1}$ in
equations \eqref{e.p.dtheta_7d}. We seek four one-forms
\begin{align*}
 \Xi_1=&\vp{1}+a_1\theta^1+a_2\theta^2+a_3\theta^4, \\
 \Xi_2=&\vp{2}+b_1\theta^1+b_2\theta^2+b_3\theta^4, \\
 \Xi_3=&\vp{3}+c_1\theta^1+c_2\theta^2+c_3\theta^4, \\
 \Xi_4=&\hp{3}+f_1\theta^1+f_2\theta^2+f_3\theta^4,
\end{align*}
with yet unknown functions $a_1,\ldots,f_3$, such that the matrix
\ben
\bma \tfrac{1}{2}\Xi_{1} & \tfrac{1}{2}\Xi_{2} & 0 & 0 \\\\
             \hp{4} & \Xi_{3}-\tfrac{1}{2}\Xi_{1} & 0 & 0 \\\\
             \hp{2} & \Xi_{4} & \tfrac{1}{2}\Xi_{1}-\Xi_{3} & -\tfrac{1}{2}\Xi_{2} \\\\
             2\hp{1} & \hp{2} & -\hp{4} & -\tfrac{1}{2}\Xi_{1}
        \ema
\een is a Cartan connection on $\P^7\to\J^1$. Calculating the
curvature for this connection we obtain that the horizontality
conditions yield
$$\der a_1=X_1(a_1)\hp{1}+X_2(a_1)\hp{2}+\inp{B}{1}\hp{3}+X_4(a_1)\hp{4}-a_1\vp{1}-a_2\vp{2}.$$
Unfortunately, none combination of the structural functions
$\inp{A}{1},\ldots,\inp{D}{2}$ and their coframe derivatives of
first order satisfies this condition. Therefore we are not able to
build a Cartan connection for an arbitrary point projective
structure. Moreover, since $\inp{B}{1}$ is a basic point invariant
together with $\inp{A}{1}$ and $\inp{C}{1}$, it seems to us
unlikely that among the coframe derivatives of
$\inp{A}{1},\ldots,\inp{D}{2}$ of any order there exists a
function satisfying the above condition. If such a function
existed it would mean that among the derivatives there is a more
fundamental function from which $\inp{B}{1}$ can be obtained by
differentiation.

Of course, we do have a Cartan connection for the point projective
geometry provided that in addition to $F_{qqq}=0$ the conditions
$\inp{B}{1}=\inp{D}{1}=0$ are imposed. However, the geometric
interpretation of these conditions is unclear.

\section{Six-dimensional Weyl geometry in the split
signature}\label{s.w6d} \noindent The construction of the
six-dimensional split signature conformal geometry given in
chapter \ref{ch.contact} has also its Weyl counterpart in the
point case. A similar construction was done by P. Nurowski,
\cite{Nur1} but he considered the conformal metric, not the Weyl
geometry. Here, apart from the tensor \ben
\wh{\gg}=2(\Omega_1-\Omega_3)\theta^2-2\Omega_3\theta^1+2\theta^4\theta^3
\een of \eqref{e.c.g33}, we also have the one-form \ben
\frac12\vp{3}.\een The Lie derivatives along the degenerate
direction $X_5+X_7$ of $\wh{\gg}$ are \ben
L_{(X_5+X_7)}\wh{\gg}=\wh{\gg}\qquad \text{and}\qquad
L_{(X_5+X_7)}\vp{3}=0.\een In this manner the pair
$(\wh{\gg},\tfrac12\vp{3})$ generates the six-dimensional
split-signature Weyl geometry $(\gg,\phi)$ on the six-manifold
$\M^6$ being the space of integral curves of $X_5+X_7$. The
associated Weyl connection is $\co(3,3)\semi{.}\real^6$-valued and
has the
following form. 
\ben
\Gamma^\mu_{~\nu}=\bma
0&\tfrac12\hp{4}&\Gamma^1_{~3}&0&\Gamma^1_{~5}
&\Gamma^1_{6}\\\\ 
\tfrac12\vp{2}&\tfrac12\vp{1}-\tfrac12\vp{3}&\Gamma^2_{~3}&-\Gamma^1_{~5}&0&\Gamma^2_{~6}\\\\ 
\tfrac12\hp{4}&0&\tfrac12\hp{3}-\tfrac12\hp{1}&-\Gamma^1_{~6}&-\Gamma^2_{~6}&0\\\\ 
0&-\tfrac12\hp{1}&\tfrac12\hp{3}&-\vp{3}&-\tfrac12\vp{2}&-\tfrac12\hp{4}\\\\ 
\tfrac12\hp{1}&0&\tfrac12\hp{2}&-\tfrac12\hp{4}&-\tfrac12\vp{1}-\tfrac12\vp{3}&0\\\\ 
-\tfrac12\hp{3}&-\tfrac12\hp{2}&0&-\Gamma^1_{~3}&-\Gamma^2_{~3}&\tfrac12\vp{1}-\tfrac32\vp{3}
\ema, \een where
\begin{align*}
\Gamma^1_{~3}=&\tfrac12\vp{2}+\tfrac12\inp{A}{2}\hp{1},\\
\Gamma^2_{~3}=&\inp{A}{3}\hp{1}+\tfrac12\inp{A}{2}\hp{2}+\inp{A}{1}\hp{4}, \\
\Gamma^1_{~5}=&\inp{D}{1}\hp{1}+\inp{B}{3}\hp{2}+(\tfrac32\inp{B}{4}-\inp{B}{1})\hp{3}
+(\inp{C}{1}-\tfrac12\inp{A}{2})\hp{4}, \\
\Gamma^1_{~6}=&(\tfrac12\inp{B}{4}-\inp{B}{1})\hp{1}-\inp{B}{2}\hp{2},\\
\Gamma^2_{~6}=&(\inp{B}{3}+\inp{D}{1})\hp{1}+\tfrac12\inp{B}{4}\hp{2}+\inp{B}{2}\hp{4}.
\end{align*}
Contrary to section \ref{s.c6d} of chapter \ref{ch.contact}, this
connection seems to have full holonomy $CO(3,3)\ltimes\real^6$ and
it is not generated by $\wh{\omega}^p$ in any simple way. It is
also never Einstein-Weyl.

\section{Three-dimensional Lorentzian geometry on solution
space}\label{s.lor3} \noindent The geometries of sections
\ref{s.ew} to \ref{s.w6d} are counterparts of respective
geometries of the contact case. The point classification, however,
contains another geometry, which is new when compared to the
contact case. This is owing to the fact that the Einstein-Weyl
geometry of section \ref{s.ew} has in general the non-vanishing
Ricci scalar, which is a weighted function and can be fixed to a
constant by an appropriate choice of the conformal gauge. Thereby
the Weyl geometry on $\S$ is reduced to a
Lorentzian metric geometry. 

These properties of the Weyl geometry are reflected at the level
of the ODEs by the fact that the equations \be\label{so22}
y'''=\frac{3}{2}\frac{(y'')^2}{y'} \ee
 and \be\label{so4}
y'''=\frac{3y'(y'')^2}{{y'}^2+1} \ee  are contact equivalent to
the trivial $y'''=0$ by means of corollary \ref{cor.c.10d_flat}
but they are mutually {\em non-equivalent} under point
transformations and possess the $\o(2,2)$ and $\o(4)$ algebra of
point symmetries respectively. Both of them generate the same flat
conformal geometry but their Weyl geometries differ. After
calculating equations \eqref{e.ewcur} we see that the only
non-vanishing component of their curvature is the Ricci scalar,
which is negative for the equation \eqref{so22} and positive for
\eqref{so4}. In this circumstances we do another reduction step in
the Cartan algorithm setting the Ricci scalar equal to $\pm 6$
respectively\footnote{We choose $\pm 6$ here to avoid large
numerical factors.}, which means $\inp{B}{4}=\pm1$, and obtain a
six-dimensional subbundle $\P^p_6$ of $\P^p$. The invariant
coframe $(\hp{1},\hp{2},\hp{3},\hp{4},\vp{1},\vp{2})$ yields the
local structure of $SO(2,2)$ or $SO(4)$ on $\P^p_6$ while the
tensor $\wh{g}=2\hp{1}\hp{3}-(\hp{2})^2$ descends to a metric
rather than a conformal class on $\S$ by means of conditions
\be\label{p.liem} L_{X_5}\wh{g}=0, \qquad L_{X_6}\wh{g}=0.\ee The
obtained metrics are locally diffeomorphic to the metrics on the
symmetric spaces $SO(2,2)/SO(2,1)$ or $SO(4)/SO(3)$.

In order to generalize this construction to a broader class of
equations we assume that the Ricci scalar of the Einstein-Weyl
geometry is non-zero \ben
6K_{qq}+\frac23F_{qqq}F_q+2F_{qqp}+\frac12F_{qq}^2 \neq 0\een and
set
$$ u_3=\sqrt{\left|6K_{qq}+\frac23F_{qqq}F_q+2F_{qqp}+\frac12F_{qq}^2\right|}$$
in the coframe of theorem \ref{th.p.1}. The
tensor $\wh{g}$ on $\P^p_6$ projects to the metric $g$ on $\S$
provided that the conditions \eqref{p.liem} still hold, which is
equivalent to \ben W=0 \quad \text{and} \quad
(\D+\tfrac23F_q)\left(
6K_{qq}+\frac23F_{qqq}F_q+2F_{qqp}+\frac12F_{qq}^2\right)=0. \een
The Cartan coframe on $\P^p_6$ is then given by
\begin{align}
 \der\hp{1} =&\vp{1}\w\hp{1}-\hp{2}\w\hp{4},\nonumber \\
 \der\hp{2} =&\vp{2}\w\hp{1}+\inv{p}{1}\hp{2}\w\hp{3}-\hp{3}\w\hp{4},\nonumber \\
 \der\hp{3}=&\vp{2}\w\hp{2}-\vp{1}\w\hp{3}+\inv{p}{2}\hp{2}\w\hp{3},\nonumber \\
 \der\hp{4} =&\vp{1}\w\hp{4}+\inv{p}{3}\hp{1}\w\hp{2}
   +\inv{p}{4}\hp{1}\w\hp{3}+\inv{p}{5}\hp{1}\w\hp{4}
   -\tfrac12\inv{p}{2}\hp{2}\w\hp{4}+\inv{p}{1}\hp{3}\w\hp{4}, \nonumber\\
 \der\vp{1} =&-\vp{2}\w\hp{4}+\inv{p}{2}\vp{2}\w\hp{1}
   +\inv{p}{6}\hp{1}\w\hp{2}+\inv{p}{7}\hp{1}\w\hp{3}+\inv{p}{4}\hp{2}\w\hp{3}
   +\inv{p}{5}\hp{2}\w\hp{4},\notag \\
 \der\vp{2}=&-\vp{1}\w\vp{2}+\inv{p}{1}\vp{2}\w\hp{3}+\inv{p}{8}\hp{1}\w\hp{2}
 +\inv{p}{9}\hp{1}\w\hp{3}+\inv{p}{10}\hp{2}\w\hp{3}
  +\inv{p}{5}\hp{3}\w\hp{4}, \nonumber
\end{align}
with some functions $\inv{p}{1},\ldots,\inv{p}{10} $ and the
Levi-Civita connection is given by \ben
 \bma
 \Gamma^1_{~1} & \Gamma^1_{~2} & 0 \\
 \Gamma^2_{~1} & 0 & \Gamma^1_{~2} \\
 0 & \Gamma^2_{~1} & -\Gamma^1_{~1}
 \ema,
\een
where
\begin{align*}
\Gamma^1_{~1}=&-\Omega_1+\tfrac12\inv{p}{2}\theta^2,\\
\Gamma^1_{~2}=&\tfrac12\inv{p}{2}\theta^1-\inv{p}{1}\theta^2-\theta^4,
\\
\Gamma^2_{~1}=&-\Omega_2+\tfrac12\inv{p}{2}\theta^3.
\end{align*}
The curvature reads \ben \bma
 R^1_{~1} & R^1_{~2} & 0 \\
 R^2_{~1} & 0 & R^1_{~2} \\
 0 & R^2_{~1} & -R^1_{~1}
 \ema,
\een

\begin{align*}
R^1_{~1}=&\tfrac12(\inv{p}{9}-\inv{p}{6})\theta^1\w\theta^2+(\tfrac14(\inv{p}{2})^2
-\inv{p}{7})\theta^1\w\theta^3+(\inv{p}{4}+X_2(\inv{p}{1})
+\tfrac12\inv{p}{1}\inv{p}{2})\theta^2\w\theta^3,\\\\
R^1_{~2}=&(\inv{p}{10}-\tfrac12
X_2(\inv{p}{2})-\tfrac14(\inv{p}{2})^2)\theta^1\w\theta^2+
(\inv{p}{4}+X_2(\inv{p}{1})+\tfrac12\inv{p}{1}\inv{p}{2})\theta^1\w\theta^3
\\
&+((\inv{p}{1})^2-X_3(\inv{p}{1}))\theta^2\w\theta^3, \\\\
R^2_{~1}=&-\inv{p}{8}\theta^1\w\theta^2+\tfrac12(\inv{p}{6}-\inv{p}{9})\theta^1\w\theta^3
+(-\inv{p}{10}+\tfrac12X_2(\inv{p}{2})+\tfrac14(\inv{p}{2})^2)\theta^2\w\theta^3.
\end{align*}


\section[Fibre-preserving case: Cartan connection on seven-dimensional\ldots]
{Fibre-preserving case: Cartan connection on seven-dimensional
bundle}\label{s.f}
\noindent The construction of a Cartan connection for the fibre
preserving case is very similar to its point counterpart. This is
due to the fact that every point symmetry of $y'''=0$ is
necessarily fibre-preserving and, as a consequence, the bundle we
will construct is also of dimension seven. Starting from the
$G_f$-structure of Introduction, which is given by the forms
\ben\bal
 \hf{1}=&u_1\omega^1, \\
 \hf{2}=&u_2\omega^1+u_3\omega^2,  \\
 \hf{3}=&u_4\omega^1+u_5\omega^2+u_6\omega^3,\\
 \hf{4}=&u_7\omega^4, \\
\eal \een and after the substitutions
\begin{align*}
 u_6=&\frac{u_3^2}{u_1},\quad\quad\quad u_7=\frac{u_1}{u_3}, \\
 u_5=&\frac{u_3}{u_1}\left(u_2-\frac{1}{3}u_3F_q\right), \\
 u_4=&\frac{u^2_3}{u_1}K+\frac{u_2^2}{2u_1}
\end{align*}
we get the following theorem.
\begin{theorem}\label{th.f.1}
An equation $y'''=F(x,y,y',y'')$ modulo fibre-preserving
transformations is described by the coframe
$(\hf{1},\hf{2},\hf{3},\hf{4},\vf{1},\vf{2},\vf{3})$ which
generates the following Cartan connection
 \be\label{e.f.conn_7d}
 \wh{\omega}^f=\bma \tfrac{1}{2}\vf{1} & \tfrac{1}{2}\vf{2} & 0 & 0 \\\\
             \hf{4} & \vf{3}-\tfrac{1}{2}\vf{1} & 0 & 0 \\\\
             \hf{2} & \hf{3} & \tfrac{1}{2}\vf{1}-\vf{3} & -\tfrac{1}{2}\vf{2} \\\\
             2\hf{1} & \hf{2} & -\hf{4} & -\tfrac{1}{2}\vf{1}
        \ema
\ee on the seven-dimensional bundle $H_3\to\P^f\to\J^2$. The group
$H_3$ is the same as in the point case
\ben H_3=\bma \sqrt{u_1}, &\frac12\frac{u_2}{\sqrt{u_1}}, & 0 & 0 \\\\
 0 & \tfrac{u_3}{\sqrt{u_1}}, & 0 & 0 \\\\
 0 & 0 & \tfrac{\sqrt{u_1}}{u_3}, &
-\tfrac12\tfrac{u_2}{\sqrt{u_1}\,\u_3} \\\\
  0 & 0 & 0 & \tfrac{1}{\sqrt{u_1}}\ema,
\een and the connection is explicitly given by \ben
\wh{\omega}^f(x^i,u_\mu)=u^{-1}\,{\omega}^f\,u+u^{-1}\der u \een
where $u\in H_3$ and \ben{\omega}^f= \bma
\tfrac{1}{2}\vp{1}^0&\tfrac{1}{2}\vp{2}^0 & 0 & 0 \\\\
\wt{\omega}^4 & \vp{3}^0-\tfrac{1}{2}\vp{1}^0 & 0 & 0 \\\\
\omega^2 & \wt{\omega}^3 & \tfrac{1}{2}\vp{1}^0-\vp{3}^0 & -\tfrac{1}{2}\vp{2}^0 \\\\
2\omega^1 & \omega^2 & -\wt{\omega}^4 & -\tfrac{1}{2}\vp{1}^0
\ema\een is given by
\begin{align*}
\omega^1=&\der y-p\der x, \notag \\
\omega^2=&\der p-q\der x, \notag \\
\wt{\omega}^3=&\der q-F\der x-\tfrac13F_q(\der p-q\der x)+K(\der y-p\der x),\notag \\
\omega^4=&\der x  \\
\vf{1}^0=&-K_q\,\omega^1+\tfrac13F_{qq}\omega^2, \notag \\
\vf{2}^0=&L\,\omega^1-K_q\,\omega^2+\tfrac13F_{qq}\wt{\omega}^3-K\omega^4, \notag \\
\vf{3}^0=&-K_q\,\omega^1+\tfrac13F_{qq}\,\omega^2
+\tfrac13F_q\omega^4.\notag
\end{align*}
\end{theorem}

The exterior differentials of the coframe are equal to
\begin{align}
 \der\hf{1} =&\vf{1}\w\hf{1}+\hf{4}\w\hf{2}+\inf{B}{1}\hf{1}\w\hf{2},\nonumber \\
 \der\hf{2} =&\vf{2}\w\hf{1}+\vf{3}\w\hf{2}+\hf{4}\w\hf{3}+\inf{B}{1}\hf{1}\w\hf{3},\nonumber \\
 \der\hf{3}=&\vf{2}\w\hf{2}+(2\vf{3}-\vf{1})\w\hf{3}+\inf{A}{1}\hf{4}\w\hf{1}+\inf{B}{1}\hf{2}\w\hf{3},\nonumber \\
 \der\hf{4} =&(\vf{1}-\vf{3})\w\hf{4}, \nonumber\\
 \der\vf{1}
 =&-\vf{2}\w\hf{4}+(\inf{D}{1}-\inf{B}{2})\hf{1}\w\hf{2}
   +\inf{B}{3}\hf{1}\w\hf{3}+ \label{e.f.dtheta_7d}  \\
   &+(2\inf{C}{1}-\inf{A}{2})\hf{1}\w\hf{4}+\inf{B}{4}\hf{2}\w\hf{3}+\inf{B}{5}\hf{2}\w\hf{4}, \nonumber \\
 \der\vf{2}=&(\vf{3}-\vf{1})\w\vf{2}+\inf{D}{2}\hf{1}\w\hf{2}+(\inf{D}{1}-2\inf{B}{2})\hf{1}\w\hf{3}
  +\inf{A}{3}\hf{1}\w\hf{4} \nonumber \\
  &+\inf{B}{6}\hf{2}\w\hf{3}+(\inf{C}{1}-\inf{A}{2})\hf{2}\w\hf{4}+\inf{B}{5}\hf{3}\w\hf{4}, \nonumber \\
 \der\vf{3}=&(\inf{D}{1}-\inf{B}{2})\hf{1}\w\hf{2}+\inf{B}{3}\hf{1}\w\hf{3}
  +(\inf{C}{1}-\inf{A}{2})\hf{1}\w\hf{4} \notag \\
  &+\inf{B}{4}\hf{2}\w\hf{3}+\tfrac12\inf{B}{5}\hf{2}\w\hf{4}, \nonumber
\end{align}
where
$\inf{A}{1},\inf{A}{2},\inf{A}{3},\inf{B}{1},\inf{B}{2},\inf{B}{3},\inf{B}{4},\inf{B}{5},\inf{B}{6},\inf{C}{1},\inf{D}{1},\inf{D}{2}$
are  functions on $\P^f$. All these invariants express by the
coframe derivatives of $\inf{A}{1},\inf{B}{1},\inf{C}{1}$, which
read
\begin{align}
 \inf{A}{1}=&\frac{u_3^3}{u^3_1}W, \nonumber \\
 \inf{B}{1}=&\frac{1}{3u_3}F_{qq}, \nonumber \\
 \inf{C}{1}=&\frac{u_2}{u_1^2}\left(\tfrac19F_{qq}F_q+\tfrac13F_{qp}+K_q\right)+ \nonumber \\
  &+\frac{u_3}{u_1^2}\left(\tfrac23F_{qq}K-\tfrac13K_qF_q-K_p-\tfrac23F_{qy}\right).\nonumber
\end{align}
If $\inf{A}{1}=0$ then $\inf{A}{2}=\inf{A}{3}=0$ and if
$\inf{B}{1}=0$ then $\inf{B}{i}=0$ for $i=2,\ldots,6$. In
particular we have \ben
\der\inf{B}{1}=-\inf{B}{2}\hf{1}+(\inf{B}{6}-\inf{B}{3})\hf{2}-\inf{B}{4}\hf{3}-\inf{B}{5}\hf{4}
-\inf{B}{1}\vf{3}.\een The flat case is given by vanishing of
$\inf{A}{1}$, $\inf{B}{1}$ and $\inf{C}{1}$.

\subsection{Fibre-preserving versus point objects}
An immediate observation about the fibre-preserving objects is
that they are closely related to their point counterparts of
theorem \ref{th.p.1}. The bundle $\P^f$ is also, like $\P^p$, a
subbundle of the contact bundle $\P^c$ given by
$\tau\colon\P^f\to\P^c$ \be\label{e.f.cf} u_4=0, \qquad u_5=0,
\qquad u_6=0, \ee and, as before, the forms
$\overset{f}{\hf{1}},\ldots,\overset{f}{\hf{4}}$ of the
fibre-preserving coframe are given by pull-backs of their contact
counterparts \ben
\overset{f}{\hf{i}}=\tau^*\overset{c}{\hf{i}},\qquad i=1,2,3,4.
\een This fact suggests that there exists a distinguished
diffeomorphism of $\P^p$ and $\P^f$ given by a similar condition.
Indeed, there is the unique diffeomorphism $\rho\colon\P^p\to\P^f$
such that
 \be\label{e.f.pf1}
\overset{p}{\hp{1}}=\rho^*\overset{f}{\hf{1}},\qquad\qquad
\overset{p}{\hp{2}}=\rho^*\overset{f}{\hf{2}},\qquad\qquad
\overset{p}{\hp{3}}=\rho^*\overset{f}{\hf{3}}. \ee It is given by
the identity map in the coordinate systems of theorems
\ref{th.p.1} and \ref{th.f.1}. The remaining one-forms are
transported as follows.\be\label{e.f.pf2} \begin{aligned}
\overset{p}{\hp{4}}=&\rho^*(\overset{f}{\hf{4}}+\tfrac12\inf{B}{1}\,\overset{f}{\hf{1}}), \\
\overset{p}{\vp{1}}=&\rho^*(\overset{f}{\vf{1}}+\inf{B}{5}\,\overset{f}{\hf{1}}
-\tfrac12\inf{B}{1}\,\overset{f}{\hf{2}}),\\
\overset{p}{\vp{2}}=&\rho^*(\overset{f}{\vf{2}}+\tfrac12\inf{B}{5}\,\overset{f}{\hf{2}}
-\tfrac12\inf{B}{1}\,\overset{f}{\hf{3}}), \\
\overset{p}{\vp{3}}=&\rho^*(\overset{f}{\vf{3}}+\tfrac12\inf{B}{5}\,\overset{f}{\hf{1}}),
\end{aligned}\ee
where \ben
\inf{B}{5}=-\frac{1}{3u_1}(\D(F_{qq})+\tfrac13F_{qq}F_q). \een The
above formulae enable us to transform easily the fibre-preserving
coframe into the point coframe.  Given the fibre-preserving
coframe $(\overset{f}{\hf{1}},\ldots,\overset{f}{\vf{3}})$ we
compute $\der \overset{f}{\hf{1}}$ and take the coefficient of the
$\overset{f}{\hf{1}}\w\overset{f}{\hf{2}}$ term. This is the
function $\inf{B}{1}$. Next we compute $\der \inf{B}{1}$,
decompose it in the fibre-preserving coframe, take minus function
that stands next to $\overset{f}{\hf{4}}$ and this is
$\inf{B}{5}$. We substitute these functions together with
$(\overset{f}{\hf{1}},\ldots,\overset{f}{\vf{3}})$ into the right
hand side of \eqref{e.f.pf1} and \eqref{e.f.pf2}, where $\rho$ is
the identity transformation of $\P^f$ and the point coframe is
explicitly constructed on $\P^f$.

Now let us consider the inverse construction, from the
fibre-preserving case to the point case. If we have  only the
point coframe $(\overset{p}{\hp{1}},\ldots,\overset{p}{\vp{3}})$
then we can not utilize eq. \eqref{e.f.pf2} since we are not able
to construct the function $\inf{B}{1}$, which is not a point
invariant\footnote{For example the point transformation
$(x,y)\to(y,x)$ destroys the condition $F_{qq}=0$.} and, as such,
does not appear among functions $\inp{A}{1},\ldots,\inp{D}{2}$ in
\eqref{e.p.dtheta_7d} or among their derivatives. However, if we
consider the point coframe {\em and} the function $\inf{B}{1}$
then the construction is possible, since $\inf{B}{5}$ is given by
the derivative $-X_4(\inf{B}{1})$ along the field $X_4$ of the
{\em point} dual frame. Therefore the passage from the point case
to the fibre-preserving case is possible if we supplement the
connection $\wh{\omega}^p$ with the function $\inf{B}{1}$. This
fact implies that each construction of the point case has its
fibre-preserving counterpart which has an additional object
generated by $\inf{B}{1}$.

\section{Fibre-preserving geometry from point
geometry}\label{s.fp}
\subsection{ Counterpart of the Einstein-Weyl geometry on $\S$}
This geometry is constructed in the following way. Let
$(\hf{1},\ldots,\vp{3})$ denotes again the fibre-preserving
coframe. Given the objects $\wh{g}=2\hf{1}\hf{3}-(\hf{2})^2$ and
\ben\wh{\phi}=\vf{3}+\frac12\inf{B}{5}\hf{1}, \een let us also
consider the function $\inf{B}{1}$, and ask under what conditions
the triple $(\wh{g},\wh{\phi},\inf{B}{1})$ can be projected to a
geometry on $\S$. There are two possibilities here, either
$\inf{B}{1}=0$ or $\inf{B}{1}\neq0$. If $\inf{B}{1}=0$ then it is
easy to see that the pair $(\wh{g},\wh{\phi})$ generates the
Einstein-Weyl geometry if only $\inf{A}{1}=\inf{C}{1}=0$, which
means that we are in the trivial case $y'''=0$.

Suppose $\inf{B}{1}\neq0$ then. For the geometry on $\S$ to exist
we need not only the conditions for the Lie transport of $\wh{g}$
and $\wh{\phi}$ but also \be\label{e.lieb}
L_{X_i}\inf{B}{1}=0,\quad \text{for}\quad i=4,5,6,\quad
L_{X_7}\inf{B}{1}=-\inf{B}{1}.\ee If all these conditions are
satisfied then $(\wh{g},\wh{\phi},\inf{B}{1})$ defines on $\S$ the
Einstein-Weyl geometry $(g,\phi)$ of the point case, which is
equipped with an additional object: a weighted function $f$ which
transforms $f\to e^{-\lambda}f$ when $g\to e^{2\lambda}g$ and is
given by the projection of $\inf{B}{1}$. The conditions for
existence of this geometry are $\inf{A}{1}=\inf{B}{5}=0$, that is
\be\label{f.lief} W=0 \qquad\text{and}\qquad
\D(F_{qq})+\tfrac13F_{qq}F_q=0.\ee As usual, the condition $W=0$
guarantees existence of $[g]$ and the other condition yields
\eqref{e.lieb}. The proper Lie transport of $\wh{\phi}$ along
$X_4$ is already guaranteed by the above conditions as their
differential consequence.

\subsection{Counterpart of the Weyl geometry on $\M^6$}
In the similar vein we show that the triple
$(\wh{\gg},\tfrac12\wh{\phi},\inf{B}{1})$, where \begin{eqnarray*}
&\wh{\gg}=2(\Omega_1-\Omega_3)\theta^2-2\Omega_3\theta^1+2\theta^4\theta^3,
\\
&\wh{\phi}=\vf{3}+\frac12\inf{B}{5}\hf{1}.
\end{eqnarray*}
 projects to the
six-dimensional split signature Weyl geometry $(\gg,\phi)$ of
chapter \ref{ch.point} section \ref{s.w6d} equipped with a
function $f$ of conformal weight $-2$.

\subsection{Counterpart of Lorentzian geometry on $\S$}
Given $(g,\phi,f)$ on $\S$ it is natural to fix the conformal
gauge so as\footnote{With possible change of the signature to make
$f$ positive.} $f=1$. This is equivalent to another substitution
\ben u_3=\tfrac13F_{qq} \een in the Cartan reduction algorithm,
which leads us to the bundle $\P^f_6$ with the following
differential system \be\bal
 \der\hf{1}=&\,\vf{1}\w\hf{1}+\hf{4}\w\hf{2},\\
 \der\hf{2}=&\,\vf{2}\w\hf{1}+\inv{f}{1}\hf{3}\w\hf{2}+\inv{f}{2}\hf{4}\w\hf{2}+\hf{4}\w\hf{3},\\
 \der\hf{3}=&-\vf{1}\w\hf{3}+\vf{2}\w\hf{2}+(2-2\inv{f}{3})\hf{3}\w\hf{2}+\inv{f}{4}\hf{4}\w\hf{1}+2\inv{f}{2}\hf{4}\w\hf{3},\\
 \der\hf{4}=&\,\vf{1}\w\hf{4}+\inv{f}{5}\,\hf{4}\w\hf{1}+(\inv{f}{3}-2)\,\hf{4}\w\hf{2}+\inv{f}{1}\,\hf{4}\w\hf{3},\\
 \der\vf{1}=&\,(2\inv{f}{3}-2)\vf{2}\w\hf{1}-\vf{2}\w\hf{4}+\inv{f}{6}\hf{1}\w\hf{2}+\inv{f}{7}\hf{1}\w\hf{3}+\inv{f}{8}\hf{1}\w\hf{4}
    -\inv{f}{5}\hf{2}\hf{4}, \\
 \der\vf{2}=&\,\vf{2}\w\vf{1}-\inv{f}{1}\vf{2}\w\hf{3}-\inv{f}{2}\vf{2}\w\hf{4}+\inv{f}{9}\hf{1}\w\hf{2}
  +\inv{f}{10}\hf{1}\w\hf{3}+\inv{f}{11}\hf{1}\w\hf{4}+\\
  &+\inv{f}{12}\hf{2}\w\hf{3}+\inv{f}{13}\hf{2}\w\hf{4}-\inv{f}{5}\hf{3}\w\hf{4}. \\
 \eal \ee
 If the conditions \eqref{f.lief}, now equivalent to
 $\inv{f}{4}=0=\inv{f}{2}$, are satisfied then $\wh{g}=2\hf{1}\hf{3}-(\hf{2})^2$
projects to a Lorentzian metric $g$ and \ben
\wh{\phi}=-2\inv{f}{5}\hf{1}+2\inv{f}{3}\hf{2}+2\inv{f}{1}\hf{3}
\een projects to a one-form $\phi$. With the Lorentzian metric
there is associated the Levi-Civita connection
$(\Gamma^\mu_{~\nu})$:
\begin{align*}
\Gamma^1_{~1}=&-\vf{1}+(\inv{f}{3}-1)\hf{2},\\
\Gamma^1_{~2}=&(\inv{f}{3}-1)\hf{1}+\inv{f}{1}\hf{2}-\hf{4},
\\
\Gamma^2_{~1}=&-\vf{2}+(\inv{f}{3}-1)\hf{3}.
\end{align*}
The covariant derivative of $\phi$ with respect to
$\Gamma^\mu_{~\nu}$ is as follows \ben \phi_{i;j}=\begin{pmatrix}
-\inv{f}{9}-(\inv{f}{5})^2&\tfrac{1}{2}\inv{f}{6}+\inv{f}{5}(\inv{f}{3}-2)&\inv{f}{12}-\inv{f}{3}(\inv{f}{3}-3)\\
\tfrac{1}{2}\inv{f}{6}+\inv{f}{5}\inv{f}{3}&2\inv{f}{12}-2\inv{f}{3}(\inv{f}{3}-2)&X_2(\inv{f}{1})+\inv{f}{1}\\
\inv{f}{12}-\inv{f}{3}(\inv{f}{3}-1)&X_2(\inv{f}{1})-\inv{f}{1}&X_3(\inv{f}{1})
\end{pmatrix}.
\een The one-form $\phi$ and the Ricci tensor satisfy the
following identities
\begin{eqnarray*}
 &\nabla_{(i}\phi_{j)}=-\Ric_{ij}-\phi_i\phi_j+(\phi^k\phi_k+2)g_{ij},\\
 &\R=2\phi^k\phi_k+6,\\
 &\der\phi=-2\ast\phi.
\end{eqnarray*}
The homogeneous model of this geometry is associated to
$y'''=\tfrac32\tfrac{(y'')^2}{y'}$.

\chapter{Classification of third-order ODEs}\label{ch.class}

\noindent Previous chapters  provided several geometric
constructions associated to third-order ODEs, now we focus on the
application of the Cartan equivalence method to the classification
of ODEs. Following \cite{Olv2} chapters 8 -- 14 we describe the
procedure which was outlined in Introduction.

By means of theorems \ref{th.c.1}, \ref{th.c.2}, \ref{th.p.1} and
\ref{th.f.1} to an ODE modulo contact/point/fibre-preserving
transformations there is associated a Cartan coframe on a bundle
$\P\to\J^2$. From theorem \ref{th.i.clas} (with obvious changes)
we know, that the underlying ODEs are equivalent if and only if
the coframes are. Thereby the problem of equivalence of ODEs is
reduced to the equivalence problem of coframes.

Consider two smooth coframes $(\omega^i)$ on $\M^n$ and
$(\cc{\omega}^i)$ on $\cc{\M}^n$. Let $(X_i)$ and $(\cc{X}_i)$ be
the dual frames. We ask under what conditions there exists a local
diffeomorphism $\Phi\colon \M \supset U\to \cc{U}\subset \cc{M}$
such that $\Phi^*\cc{\omega}^i=\omega^i$. To answer this question
we need to define several notions. The structural equations for
$(\omega^i)$ read
$\der\omega^i=\frac12T^i_{jk}\omega^j\w\omega^k$, where
$T^i_{jk}=T^i_{[jk]}$ are smooth functions since the coframe is
smooth. The functions $T^i_{jk}$ and their coframe derivatives
$T^i_{jk|l}=X_l(T^i_{jk})$, $T^i_{jk|lm}=X_m(X_l(T^i_{jk}))$ etc.
of any order are called structural functions of the coframe.
$T^i_{~jk}$ are the structural functions of zero order,
$T^i_{jk|l}$ are the structural functions of first order and so
on. Since pull-back commutes with exterior differentiation all the
structural functions are relative invariants of the equivalence
problem of coframes. We say that smooth functions $f_1,\ldots,f_k$
are independent at a point $w$ if $\der f_1\w\ldots\w\der f_k\neq
0$ at $w$. The rank of a coframe $(\omega^i)$ at $w$ is defined to
be the maximal number of its independent structural functions at
$w$. The order $s$ at $w$ of a coframe of rank $r$ at $w$ is the
smallest natural number such that among structural functions of
order at most $s$ there are $r$ functions independent at $w$.

We call a smooth coframe regular in an open set $U$ if for each
$j\geq 0$ the number of its independent structural functions of
order $j$ is constant on $U$. In particular the rank and the order
of a regular coframe are both constant on $U$. From now on we
confine our considerations to coframes regular on sufficiently
small topologically trivial open subsets $U$ of $\M$. For a
regular coframe of order $s$ one may choose a set of $r$
independent structural functions $I_1,\ldots,I_r$, which are of
order at most $s$ and all the remaining structural functions
$T_\sigma$ of any order are expressible in terms of $I_j$, $
T_\sigma=f_\sigma(I_1,\ldots,I_r)$. Indeed, all the $(s+1)$-order
functions are of this form by definition of the order of a
coframe, whereas all the functions of order $s+2$ or greater are
their coframe derivatives. These derivatives depend on $I_j$ and
$I_{j|k}$ but $I_{j|k}$ are functions of order at most $s+1$ so
their are also functions of $I_j$.

Thereby all the structural functions are described by i) the set
$(I_1,\ldots,I_r)$ and ii) the formulae which characterize how the
$(s+1)$-order functions are expressed by $(I_1,\ldots,I_r)$. Both
these objects may be encoded in the so called classifying
function. The classifying function for a coframe of order $s$ is a
function $\T\colon\M\to\real^N$ given by all the structural
functions of order at most $s+1$, which are lexicographically
ordered with respect to their indices, and $N$ is the number of
these structural functions \ben \T:w\in U\longmapsto
(T^i_{jk},T^i_{jk|l_1},\ldots,T^i_{jk|l_1\ldots
l_{s+1}})\in\real^N. \een By smoothness and regularity of the
underlying coframe the graph $\T(U)$ is an $r$-dimensional
submanifold in $\real^N$. The last ingredient which we need is
definition of overlapping; two $r$-dimensional submanifolds of
$\real^N$ overlap if their intersection is also an $r$-dimensional
submanifold of $\real^N$. Now we are in position to cite the
following
\begin{theorem}[\cite{Olv2}, theorem 14.24]\label{th.cc.clas}
Let $(\omega^i)$ and $(\cc{\omega}^i)$ be smooth, regular coframes
defined, respectively, on $\M^n$ and $\cc{\M}^n$. There exists a
local diffeomorphism $\Phi:\M\to\cc{\M}$ such that
$\Phi^*\cc{\omega}^i=\omega^i$, $i=1,\ldots,n$, if and only if the
coframes have the same order $s=\cc{s}$ and the graphs $\T(\M)$
and $\cc{\T}(\cc{\M})$ of their classifying functions overlap.
Moreover, if $w_0\in\M$ and $\cc{w}_0 \in \cc{\M}$ are any points
mapping to the same point \ben z_0=\T(w_0)=\cc{\T}(\cc{w}_0)\in
\T(\M)\cap\cc{\T}(\cc{\M})\een on the overlap, then there is a
unique equivalence map $\Phi$ such that $\cc{w_0}=\Phi(w_0)$.
\end{theorem}
An important notion is a symmetry of a coframe. This is a
diffeomorphism $\Phi:\M\to\M$ such that $\Phi^*\omega^i=\omega^i$.
We have
\begin{theorem}[\cite{Olv2}, theorem 14.26]\label{th.cc.sym}
The symmetry group $G$ of a regular coframe $(\omega^i)$ of rank
$r$ on $\M^n$ is a local Lie group of transformations of dimension
$n-r$.
\end{theorem}
There is a class of coframes which admit a particularly simple
description in this language. These are coframes of rank zero,
whose all structural functions $T^i_{~jk}$ are constant. These
coframes are of order zero and the image of their classifying
function is a point in $\real^N$. Thus, two rank-zero coframes are
equivalent if and only if their structural constants are equal,
$T^i_{~jk}=\cc{T}^i_{~jk}$. Furthermore, a coframe of rank zero
has an $n$ dimensional Lie symmetry group $G$. The algebra of this
group is precisely the algebra with structural constants
$T^i_{~jk}$ and the coframe may be interpreted as the left
invariant coframe defining a local structure of $G$ on $\M$.

Let us turn to ODEs. In order to
 determine whether or not two given ODEs are
contact/point/fibre-preserving equivalent one constructs the
respective invariant coframes for both the ODEs, calculates the
structural functions, which are now relative invariants for the
underlying ODEs, and applies theorem \ref{th.cc.clas} provided
that the coframes are regular. In particular,
contact/point/fibre-preserving symmetry group of the underlying
ODE is isomorphic to the symmetry group of the associated coframe.

In practice there are three difficulties in this approach. First
difficulty is that the coframes and invariants we have built so
far are defined on manifolds $\P$ larger then $\J^2$. As a
consequence, the invariants contain not only $x,y,p,q$ but the
auxiliary bundle variables $u_\mu$ as well. Since it is natural to
describe ODEs in terms of $x,y,p,q$ only, first of all we must
finish the reduction of the group parameters, so that finally no
free $u_\mu$ remains and we obtain a coframe on $\J^2$ from which
we compute invariants of ODEs depending on $x,y,p,q$ only.

Second problem is that the above method requires regularity of the
Cartan coframes. Thus we have to restrict our consideration to the
class of regular ODEs defined as follows.
\begin{definition}\label{def.cc.reg}
An ODE $y'''=F(x,y,y',y'')$ is regular if it is given by a locally
smooth function $F$ such that the Cartan coframe obtained after
maximal possible reduction of the structural group is regular.
\end{definition}

Third and essential obstacle for the full classification is that
the task of finding whether or not two graphs of functions
$\T\colon\J^2\to \real^N$ overlap is highly nontrivial,
particularly in our cases, where the components of $\T$ are
compound functions depending on $x,y,p,q$ through $F$. Taking this
into account we restrict the classification to two classes of
equations for which we are able to find compact criteria for the
equivalence.
\begin{itemize}
\item[i)] The regular equations possessing large contact or point
symmetry groups, that is symmetry groups of dimension at least
four. \item[ii)] The regular equations fibre-preserving equivalent
to reduced Chazy classes II, IV -- VII and XII.
\end{itemize}

For these two families we carry out the classification to the very
end, however we also provide some partial result in the case of
totally arbitrary ODEs. This partial result, theorems
\ref{th.cc.nW4d} and \ref{th.cc.W4d},  is the explicit
construction of the invariant coframes on $\J^2$ in the contact
case without further analysis of the classifying function.

\section{Equations with large contact symmetry
group}\label{s.cc.large}

\noindent This class of ODEs is particularly convenient for
characterization since, as we shall see, to these equations we can
always associate a Cartan coframe of rank zero, with the only
exception of ODEs contact equivalent to general linear ODEs
\eqref{e.cc.glin}. These exceptional ODEs are characterized by the
fact that their symmetry group act intransitively on $\J^2$, cf
corollary \ref{cor.cc.linear} and proposition \ref{prop.cc.int}

Let us begin with the ten-dimensional coframe of theorem
\ref{th.c.1}. As we said in section \ref{s.c.furred} of chapter
\ref{ch.contact}, ODEs fall into three main classes:
\begin{itemize}
\item $W=0$ and $F_{qqqq}=0$, \item $W\neq 0$, \item $W=0$ but
$F_{qqqq}\neq 0$.

\end{itemize}
If $F_{qqqq}=W=0$ then we are in the situation of corollary
\ref{cor.c.10d_flat} and this is the only case of the ODEs with a
ten-dimensional contact symmetry group, since for any other ODE
there are non-constant relative invariants in equations
\eqref{e.c.dtheta_10d} and the dimension of the symmetry group is
less then ten. Below we consider the case $W\neq 0$.

\subsection{ODEs with $W\neq 0$}
For these equations we have the five dimensional coframe of
theorem \ref{th.c.2}:
\begin{align}
 d\theta^1=&\Omega\w\theta^1-\theta^2\w\theta^4, \nonumber \\
 d\theta^2=&\Omega\w\theta^2+\inc{a}{}\,\theta^1\w\theta^4-\theta^3\w\theta^4,\nonumber\\
 d\theta^3=&\Omega\w\theta^3+\inc{b}{}\,\theta^1\w\theta^2+\inc{c}{}\,\theta^1\w\theta^3
  -\theta^1\w\theta^4+\inc{e}{}\,\theta^2\w\theta^3+\inc{a}\,\theta^2\w\theta^4, \label{e.cc.dtheta_5d} \\
 d\theta^4=&\inc{f}{}\,\theta^1\w\theta^2+\inc{g}{}\,\theta^1\w\theta^3
 +\inc{h}{}\,\theta^1\w\theta^4
  +\inc{k}{}\,\theta^2\w\theta^3-\inc{e}{}\,\theta^2\w\theta^4, \nonumber  \\
 d\Omega=&\inc{l}{}\,\theta^1\w\theta^2+(\inc{f}{}-\inc{a}{}\inc{k}{})\,\theta^1\w\theta^3
 +\inc{m}{}\,\theta^1\w\theta^4
  +\inc{g}{}\,\theta^2\w\theta^3+\inc{h}{}\,\theta^2\w\theta^4.\nonumber
 \end{align}
 with functions $\inc{a}{},\inc{b}{},\inc{e}{},\inc{h}{},\inc{k}{}$ given by
 \begin{align}
\inc{a}{}=& \frac{1}{\sqrt[3]{W^2}}\left(K+\frac{1}{18}Z^2+\frac{1}{9}ZF_q-\frac{1}{3}\D Z\right), \nonumber \\
\inc{b}{}=&\frac{1}{3u\sqrt[3]{W^2}}\bigg(\frac{1}{27}F_{qq}Z^2+\left(K_q-\frac{1}{3}Z_p-\frac{2}{9}F_qZ_q\right)Z+\nonumber\\
   &+\left(\frac{1}{3}\D Z-2K\right)Z_q+Z_y+F_{qq}K-3K_p-K_qF_q-F_{qy}+W_q\bigg), \nonumber \\
\inc{e}{}=&\frac{1}{u}\left(\frac{1}{3}F_{qq}+\frac{1}{W}\left(\frac{2}{9}W_qZ-\frac{2}{3}W_p-\frac{2}{9}W_qF_q \right)\right), \label{e.cc.basfun_5d}\\
\inc{h}{}=&\frac{1}{3u\sqrt[3]{W}}\bigg(\left(\frac{1}{9}W_qZ^2-\frac{1}{3}W_pZ+W_y-\frac{1}{3}W_q\D Z\right)\frac{1}{W}+\nonumber  \\
 &+\D Z_q+\frac{1}{3}F_qZ_q\bigg),\nonumber \\
 \inc{k}{}=&\frac{1}{u^2\sqrt[3]{W}}\left(\frac{2W_q^2}{9W}-\frac{W_{qq}}{3}\right). \nonumber
\end{align}
This means in view of theorem \ref{th.cc.sym} that if $W\neq 0$
then we have symmetry groups of dimension at most five. We proved
in corollary \ref{cor.c.5d_flat} that the assumption
$\inc{a}{}=const$ and $\inc{k}{}=0$ implies that the underlying
ODE has a five-dimensional group of symmetry. Now we easily see
that there are no other equations with a five-dimensional symmetry
group here. In fact, this property requires that all the functions
in \eqref{e.cc.dtheta_5d} are constants, in particular
$\inc{a}{}=const$ and $\inc{k}{}=const$. But $\inc{k}{}\sim
u^{-2}$ so it is a constant function on the bundle $\P^c_5$ iff it
vanishes and we are back in corollary \ref{cor.c.5d_flat}. At this
stage we know that all the ODEs with $W\neq 0$ other than
$y'''=-2\mu y'+y$ have the symmetry group of dimension at most
four. Next we reduce the last free bundle parameter $u$. The
Cartan algorithm bifurcates at this point:
\begin{itemize}
 \item[i)] The functions $\inc{b}{},\inc{e}{},\inc{h}{}$ and $\inc{k}{}$ in \eqref{e.cc.dtheta_5d} vanish.
 \item[ii)] At least one function among $\inc{b}{},\inc{e}{},\inc{h}{},\inc{k}{}$ does not vanish.
\end{itemize}
We will discuss both  these  possibilities consecutively.
\paragraph{i)} Utilizing the Jacobi identity for
\eqref{e.cc.dtheta_5d} we get that all the functions but
$\inc{a}{}$ vanish, $\inc{a}{}=\inc{a}{}(x)$ and  $\der\inc{a}{}
=(\inc{a}{})_4\theta^4$, $\der
(\inc{a}{})_4=(\inc{a}{})_{44}\theta^4$ etc, so neither
$\inc{a}{}$ nor its derivatives of any order contain the last
auxiliary variable $u$ and the full reduction of the structural
group cannot be done. We check that the ODE \be\label{e.cc.glin}
y'''=-2\mu(x)y'+(1-\mu'(x))y, \ee with an arbitrary smooth
function $\mu(x)$, satisfies
$\inc{b}{}=\inc{e}{}=\inc{h}{}=\inc{k}{}=0$ and
$\inc{a}{}=\mu(x)$. Furthermore, a linear third-order ODE of
general form satisfies either $W=0$, in which case it is
equivalent to $y'''=0$, or $W\neq 0$, in which case
$\inc{b}{},\inc{e}{},\inc{h}{}$ and $\inc{k}{}$ vanish.  Thus case
i) describes the ODEs satisfying $W\neq 0$ and linearizable
through contact transformations, in particular the linear
equations with constant coefficients are distinguished by the
additional condition $\inc{a}{}=const$.
\begin{corollary} \label{cor.cc.linear}
The following third-order ODEs are linearizable via contact
transformations of variables
\begin{itemize}
 \item The equations satisfying $W=F_{qqqq}=0$. These are equivalent to
 \ben y'''=0 \een
 and admit the group $O(3,2)$ of contact symmetries.
 \item The equations satisfying $W\neq0$, $\inc{b}{}=\inc{e}{}=\inc{h}{}=\inc{k}{}=0$ and
 $\inc{a}{}=\mu(x)$, where $\mu(x)$ is any smooth non-constant real function. These are
 equivalent to the general linear equation
 \ben y'''=-2\mu(x)y'+(1-\mu'(x))y\een
 and admit a four-dimensional group of contact symmetries. The
 symmetry group acts on three-dimensional orbits in $\J^2$. If
 such an
 equation is given in the above form, then symmetries are
 generated by \ben \bal V_i&=f_i\partial_y+f'_i\partial_p+f''_i\partial_q, \qquad i=1,2,3\\
  V_4&=y\partial_y+p\partial_p+q\partial_q,\eal \een
where $f_1,f_2,f_3$ are any three functionally independent
solutions of the ODE.

 \item The equations satisfying $W\neq0$, $\inc{k}{}=0$ and
 $\inc{a}{}=\mu\in\real$. These are equivalent to the general linear
 equation with constant coefficients \ben y'''=-2\mu y'+y \een
 and admit the group $\real^2\ltimes_\mu\real^3$ of contact
 symmetries.
\end{itemize}
The equations with different values of $\mu(x)$ or $\mu$ are
non-equivalent.
\end{corollary}

\paragraph{ii)} In this case we use a non-vanishing function among
$\inc{b}{},\inc{e}{},\inc{h}{},\inc{k}{}$ in
\eqref{e.cc.dtheta_5d} to do the last reduction and eventually
obtain a coframe on $\J^2$. The substitution is as follows.
\begin{itemize}
\item[1.] If $\inc{k}{}\neq 0$ then we set $\inc{k}{}=\epsilon=\pm
1$ depending on the sign of the quantity
  $\tfrac{1}{\sqrt[3]{W}}(\tfrac{2W_q^2}{9W}-\tfrac{W_{qq}}{3})$, which gives the substitution
  \be\label{e.cc.red10}
   u=\frac{1}{\sqrt[6]{|W|}}\sqrt{\left|\frac{2W_q^2}{9W}-\frac{W_{qq}}{3}\right|}
  \ee
in the coframe of theorem \ref{th.c.2}. \item[2.] If $\inc{k}{}=0$
and $\inc{e}{}\neq 0$, then we set $\inc{e}{}=1$ and substitute
  \ben
   u=\frac{1}{3}F_{qq}+\frac{1}{W}\left(\frac{2}{9}W_qZ-\frac{2}{3}W_p-\frac{2}{9}W_qF_q \right).
  \een
\item[3.] If $\inc{k}{}=\inc{e}{}=0$ and $\inc{h}{}\neq 0$, then
we set $\inc{h}{}=1$ and
  \begin{align*}
   u=&\frac{1}{3\sqrt[3]{W}}\bigg(\left(\frac{1}{9}W_qZ^2-\frac{1}{3}W_pZ+W_y-\frac{1}{3}W_q\D Z\right)\frac{1}{W}+ \\
   &+\D Z_q+\frac{1}{3}F_qZ_q\bigg).\nonumber
  \end{align*}
\item[4.] Finally, if $\inc{k}{}=\inc{e}{}=\inc{h}{}=0$ and
$\inc{b}{}\neq 0$, then we set $\inc{h}{}=1$ and
  \begin{align}
u=&\frac{1}{3\sqrt[3]{W^2}}\bigg(\frac{1}{27}F_{qq}Z^2+\left(K_q-\frac{1}{3}Z_p-\frac{2}{9}F_qZ_q\right)Z\label{e.cc.red20}\\
    &+\left(\frac{1}{3}\D Z-2K\right)Z_q+Z_y+F_{qq}K-3K_p-K_qF_q-F_{qy}+W_q\bigg). \nonumber
 \end{align}
\end{itemize}
In this manner we obtain
\begin{theorem}\label{th.cc.nW4d}
Let $y'''=F(x,y,y',y'')$ be an ODE such that i) $W\neq 0$, ii) the
functions $\inc{b}{},\inc{e}{},\inc{h}{},\inc{k}{}$ in
\eqref{e.cc.basfun_5d} do not vanish simultaneously. The contact
invariant information on this ODE is given by the following Cartan
coframe $(\theta^1,\theta^2,\theta^3,\theta^4)$ on $\J^2$:
\be\label{e510}\begin{aligned}
 \theta^1=&u\sqrt[3]{W}\omega^1, \\
 \theta^2=&u\left(\frac{1}{3}Z\omega^1+\omega^2\right), \\
 \theta^3=&\frac{u}{\sqrt[3]{W}}\left(\left(K+\frac{1}{18}Z^2\right)\omega^1
 +\frac{1}{3}\left(Z-F_q\right)\omega^2+\omega^3\right), \\
 \theta^4=&\left(\frac{1}{9}\frac{W_q}{\sqrt[3]{W^2}}Z-\frac{1}{3}\sqrt[3]{W}Z_q\right)\omega^1
 +\frac{W_q}{3\sqrt[3]{W^2}}\omega^2+\sqrt[3]{W}\omega^4,
 \end{aligned}\ee
where $u$ is given by eq. \eqref{e.cc.red10} --
\eqref{e.cc.red20}, depending on which functions among
$\inc{b}{},\inc{e}{},\inc{h}{},\inc{k}{}$ are non-zero.
\end{theorem}
The exterior derivatives of the above coframe read
\be\label{e.cc.dtheta_nW4d}\begin{aligned}
  d\theta^1 =&\,a\,\theta^1\w\theta^2+\inc{I}{1}\,\theta^1\w\theta^3
    +\inc{I}{2}\,\theta^1\w\theta^4-\theta^2\w\theta^4, \\
  d\theta^2 =& \,e\,\theta^1\w\theta^2+\inc{I}{3}\,\theta^1\w\theta^4
   +\inc{I}{1}\,\theta^2\w\theta^3+\inc{I}{2}\,\theta^2\w\theta^4
   -\theta^3\w\theta^4, \\
  d\theta^3 =&\,h\,\theta^1\w\theta^2+k\,\theta^1\w\theta^3
   -\,\theta^1\w\theta^4+\inc{I}{4}\,\theta^2\w\theta^3
   +\inc{I}{3}\,\theta^2\w\theta^4+\inc{I}{2}\,\theta^3\w\theta^4,\\
  d\theta^4 =& \,m\,\theta^1\w\theta^2+\,n\,\theta^1\w\theta^3
   +\,s\,\theta^1\w\theta^4+\inv{\epsilon}{1}\,\theta^2\w\theta^3
   -(a+\inc{I}{4})\,\theta^2\w\theta^4,
 \end{aligned}\ee
where $\inv{\epsilon}{1}=\pm 1,\,0$ and $\inc{I}{1}$,
$\inc{I}{2}$, $\inc{I}{3}$, $\inc{I}{4}$, $a$, $e$, $h$, $k$, $m$,
$n$, $s$ are functions. 
The most important invariants read
\begin{align*}
 \inv{\epsilon}{1}=&\sgn\left(2W_q^2-3WW_{qq}\right),\notag \\
 \inc{I}{1}=&\frac{9W_{qqq}W^2-9W_{qq}W_qW+4W_q^3}{2(2W_q^2-3WW_{qq})\sqrt{\left|2W_q^2-3WW_{qq}\right|}},\notag \\
\inc{I}{2}=&\frac{1}{6\sqrt[3]{W}(3W_{qq}W-2W_q^2)}\Big((6WF_q-6ZW)W_{qq}+4ZW^2_q-4F_qW_q^2\notag \\
  &+(3F_{qq}W-12W_p-6WZ_q)W_q+18WW_{qp}-9W^2Z_{qq}-9W^2F_{qqq}\Big),\notag \\
 \inc{I}{3}=&\frac{1}{\sqrt[3](W^2}\left( K-\frac{1}{3}\D Z+\frac{1}{18}Z^2+\frac{1}{9}F_qZ\right), \\
 \inc{I}{4}=&\frac{1}{6\sqrt[3]{W}(2W_q^2-WW_{qq})\sqrt{\left|2W_q^2-WW_{qq}\right|}}\Big(9F_qWW_{qqq}-9ZWW_{qqq}\notag \\
   &+(15F_qWW_q-15ZWW_q-27W_pW+18WF_{qq}-18WZ_q)W_{qq} \notag \\
   &+8F_qW_q^3-8ZW_q^3+(6WZ_q+24W_p-9F_{qq}W)W_q^2 \notag \\
   &-(9W^2Z_{qq}+18WW_{qp}+9W^2F_{qqq})W_q+27W_{qqp}W^2 \Big).\notag
\end{align*}

\subsection{Equations satisfying $W\neq 0$ and admitting large contact symmetry
groups}\label{s.cc.nw} Apart from the contact linearizable
equations all the equations admitting contact symmetry group are
characterized by the coframe \eqref{e510}. By virtue of theorem
\ref{th.cc.sym}, an ODE associated with the coframe \eqref{e510}
admits a four-dimensional symmetry group if all its relative
invariants in \eqref{e.cc.dtheta_nW4d} are constant.

As we said, two ODEs associated with the coframe \eqref{e510} and
admitting a four-dimensional contact symmetry group are equivalent
if and only if their respective invariants have the same constant
value.

Below we describe our method of finding these equations. First we
assume that the invariants
$\inc{I}{1},\inc{I}{2},\inc{I}{3},\inc{I}{4}$ are constant, which
is a necessary condition for a large symmetry group. Then we close
the system \eqref{e.cc.dtheta_nW4d} and the identities
$\der^2\hc{i}=0$ give us information about the remaining
invariants, for example
\ben\begin{aligned}\der^2\hc{1}=&(\inc{I}{1}\inc{I}{3}
-\inc{I}{2}\inc{I}{4}+e+s)\hc{1}\w\hc{4}\w\hc{2}+
(a\inc{I}{1}+\inc{I}{1}\inc{I}{4}-\inv{\epsilon}{1}\inc{I}{2}+n)\hc{1}\w\hc{3}\w\hc{2}+\\
&+(a-\inc{I}{1}\inc{I}{2})\hc{1}\w\hc{4}\w\hc{3}+ \der
a\w\hc{1}\w\hc{2}=0. \end{aligned} \een This equation yields
$\der^2\hc{1}\w\hc{2}=0=(a-\inc{I}{1}\inc{I}{2})\hc{1}\w\hc{4}\w\hc{3}\w\hc{2}$,
hence $a=\inc{I}{1}\inc{I}{2}=const$ and we get that $n$ and $e+s$
are also constants expressed by $\inc{I}{j}$. Next
$\der^2\hc{2}\w\hc{2}=0$ gives $e-k+s=0$ and
$\der^2\hc{3}\w\hc{1}=0$ gives $k=const$, so $a,e,k,n$ and $s$ are
constant. Continuing this reasoning we also find that $h$ and $m$
are constant. As a consequence
$\inc{I}{1},\inc{I}{2},\inc{I}{3},\inc{I}{4}=const$ is the
necessary and sufficient condition for an ODE to admit a large
symmetry group.

In these circumstances the identities $\der^2\hc{i}=0$ become a
system of quadratic algebraic equations for $\inc{I}{1},\ldots,s$.
We solve this system by the method of consecutive substitutions.
Doing so we find  that the system has no solutions if
$\inv{\epsilon}{1}=0$. This fact implies that the ODEs which we
search exist only in the branch 1 above defined by the
normalization \eqref{e.cc.red10}. In the case
$\inv{\epsilon}{1}=\pm 1$ the system is underdetermined and we
express all the invariants by $\inv{\epsilon}{1}$ and
$\inc{I}{1}$, whose value is arbitrary except for $\inc{I}{1}=
0,\,\pm 3/\sqrt{2}$. 
 \begin{align*}
\inc{I}{3}=&\frac{3(\inc{I}{2})^2}{8(\inc{I}{1})^2}(3\inv{\epsilon}{1}-(\inc{I}{1})^2),
&
\inc{I}{4}=&\frac{\inc{I}{2}}{2\inc{I}{1}}(3\inv{\epsilon}{1}-2(\inc{I}{1})^2),\\\\
 a=&\inc{I}{1}\inc{I}{2}, &
 e=&\frac{1}{8 \inc{I}{1}}(\inc{I}{2})^2(9\inv{\epsilon}{1}-5(\inc{I}{1})^2),
 \\\\
 h=&-\frac{3\inv{\epsilon}{1}}{16(\inc{I}{1})^3}(16(\inc{I}{1})^2
 -3(\inc{I}{2})^3(\inc{I}{1})^2+9(\inc{I}{2})^3\inv{\epsilon}{1}), &
 n=&-\frac12\inv{\epsilon}{1}\inc{I}{2},\\\\
 m=&\frac{\inv{\epsilon}{1}(\inc{I}{2})^3}{8(\inc{I}{1})^2}((\inc{I}{2})^2-9\inv{\epsilon}{1}),&
 k=&\frac{5(\inc{I}{2})^2}{8\inc{I}{1}}(3\inv{\epsilon}{1}-(\inc{I}{1})^2),
 \\\\
 s=&-\frac{3\inv{\epsilon}{1}(\inc{I}{2})^2}{4\inc{I}{1}},& &
\end{align*}
and
$$\inc{I}{2}=2\sqrt[3]{\frac{(\inc{I}{1})^2}{(2(\inc{I}{1})^2-9\inv{\epsilon}{1})}}.$$

Thereby we have obtained all the structural equations which {\em
may be generated} by ODEs admitting large symmetry groups.
However, it is unknown whether every admissible values of the
invariants are realized by the ODEs. In order to find the
equations we apply two approaches. The first approach is as
follows. We choose some simple ODEs, for example $F=q^\alpha$ or
$F=e^q$, calculate the equations \eqref{e.cc.dtheta_nW4d} for them
and check if they satisfy the large symmetry conditions. By this
method we find the ODEs which realize some but not all the
admissible values of $\inv{\epsilon}{1}$ and $\inc{I}{1}$. In
order to find the remaining ODEs we use the fact that when all the
invariants are constant then the coframe $(\hc{i})$ is a
left-invariant coframe of a four-dimensional Lie group. This group
is precisely the group of symmetry of the underlying ODE which
acts on $\J^2$ and the frame dual to $(\hc{i})$ is a system of
infinitesimal symmetry generators for the differential equation.
We integrate the coframe $(\hc{i})$, that is we find the explicit
formulae for the forms $\hc{i}$ which satisfy the equations
\eqref{e.cc.dtheta_nW4d} with constant invariants
$\inc{I}{1},\ldots,s$. Having found formulae for the symmetry
generators we find their common invariant functions on $\J^2$
which are our desired ODEs. On integrating the coframe $(\hc{i})$
we used M. MacCallum's classification of real four-dimensional Lie
algebras \cite{MacC} to transform the frames into a simple
canonical form, which simplified the calculations.

In this way we have found the full list of third-order ODEs
satisfying the condition $W\neq 0$ and admitting at least four
dimensional group of contact symmetries. We gather these equations
in table \ref{t.cc.1} together with equations of the branch $W=0$,
which is analyzed below. Two equations of the same type given in
the tables are equivalent if and only if the constants they
involve are equal. For example $F=q^\mu$ and $F=q^\nu$ are
equivalent provided that $\mu=\nu$. We give the necessary and
sufficient conditions for any ODE to be contact equivalent to any
member of the list. We also attach the description of the symmetry
algebra. Our result agrees with the part of B. Doubrov and B.
Komrakov's classification \cite{Dub1} referring to third-order
ODEs, however we provide other canonical forms of the considered
equations.

\subsection{ODEs with $W=0$ and $F_{qqqq}\neq0$}

We repeat the scheme explained above. This time we start from the
coframe of theorem \ref{th.c.1}. The following substitutions \ben
\bal &u_1=u_3^2\sqrt{\left|\frac{u_3}{F_{qqqq}}\right|}, &
&u_2=-3u_3\frac{K_{qqq}}{F_{qqqq}}, \\
&u_4=u_3\sqrt{\left|\frac{u_3}{F_{qqqq}}\right|}\left(\frac{3K_{qqqq}}{F_{qqqq}}-\frac{12F_{qqqqq}K_{qqq}}{5F^2_{qqqq}}\right),
&
&u_5=-u_3\sqrt{\left|\frac{u_3}{F_{qqqq}}\right|}\frac{F_{qqqqq}}{5F_{qqqq}},
\\ &u_6=0 & &  \eal \een reduce the bundle $\P^c$ to a five-dimensional
subbundle. The last reduction is as follows.
\begin{itemize}
\item[1.] If
\be\label{e.cc.red40}
 2L_{qq}F_{qqqq}-3K^2_{qqq}\neq 0,
\ee then \be
 u_3=\sqrt[3]{9L_{qq}-\frac{27K^2_{qqq}}{2F_{qqqq}}}.
\ee
\item[2.] If  $L_{qq}F_{qqqq}-3K^2_{qqq}=0$ but
\be\label{e.cc.red50}
 5F_{qqqqqq}F_{qqqq}-6F^2_{qqqqq}\neq 0,
\ee
then
\be
 u_3=\frac{25F^3_{qqqq}}{5F_{qqqqqq}F_{qqqq}-6F^2_{qqqqq}}.
\ee

\item[3.] If $L_{qq}F_{qqqq}-3K^2_{qqq}=0$  and  $5F_{qqqqqq}F_{qqqq}-6F^2_{qqqqq}=0$, but
\be\label{e.cc.red60} \frac{1}{3}F_{qq}+\frac{1}{F_{qqqq}}\left(\frac{18}{5}K_{qqqq}+\frac{2}{5}F_{qqqqp}+\frac{2}{15}F_qF_{qqqqq}\right)-\frac{12F_{qqqqq}K_{qqq}}{5F^2_{qqqq}}\neq 0
\ee
then
\be\label{e.cc.red70} u_3=\frac{1}{3}F_{qq}+\frac{1}{F_{qqqq}}\left(\frac{18}{5}K_{qqqq}+\frac{2}{5}F_{qqqqp}+\frac{2}{15}F_qF_{qqqqq}\right)
 -\frac{12F_{qqqqq}K_{qqq}}{5F^2_{qqqq}}.
\ee
\end{itemize}
\begin{remark}
The three quantities defined in \eqref{e.cc.red40},
\eqref{e.cc.red50}, \eqref{e.cc.red60} cannot vanish
simultaneously because it would be in contradiction to the
condition $F_{qqqq}\neq 0$.
\end{remark}
Thus we have the following
\begin{theorem}\label{th.cc.W4d}
Let $y'''=F(x,y,y',y'')$ be an ODE satisfying $W=0$ and
$F_{qqqq}\neq 0$. The contact invariant information on the ODE is
given by the following Cartan coframe
$(\theta^1,\theta^2,\theta^3\theta^4)$ on $\J^2$:
\begin{align*}
 \theta^1=&u^2\sqrt{\left|\frac{u}{F_{qqqq}}\right|}\omega^1,\nonumber \\
 \theta^2=&u\left(-\frac{3K_{qqq}}{F_{qqqq}}\omega^1+\omega^2\right),\nonumber \\
 \theta^3=&\sqrt{\left|\frac{F_{qqqq}}{u}\right|}\left(\left(K+\frac{9K^2_{qqq}}{2F^2_{qqqq}}
 \right)\omega^1-\left(\frac{F_q}{3}+\frac{3K_{qqq}}{F_{qqqq}}\right)\omega^2+\omega^3\right), \\
 \theta^4=&u\sqrt{\left|\frac{u}{F_{qqqq}}\right|}\left
 (\left(\frac{3K_{qqq}}{F_{qqqq}}-\frac{12F_{qqqqq}K_{qqq}}{5F^2_{qqqq}}\right)\omega^1
 -\frac{F_{qqqqq}}{5F_{qqqq}}\omega^2+\omega^4\right),\nonumber
 \end{align*}
where $u$ is the function of $x,y,p,q$ given by the formulae
\eqref{e.cc.red40} -- \eqref{e.cc.red70}.
\end{theorem}

The exterior derivatives of the above forms are the following
 \begin{align*}
  d\theta^1 &=a\,\theta^1\w\theta^2+\inc{I}{5}\,\theta^1\w\theta^3
    +\inc{I}{6}\,\theta^1\w\theta^4-\theta^2\w\theta^4, \nonumber \\
  d\theta^2 &=e\,\theta^1\w\theta^2+\inv{\epsilon}{2}\,\theta^1\w\theta^4
    +\frac{2}{5}\inc{I}{5}\,\theta^2\w\theta^3+\frac{2}{5}\inc{I}{6}\,\theta^2\w\theta^4
    -\theta^3\w\theta^4, \nonumber \\
  d\theta^3 &=f\,\theta^1\w\theta^2+g\,\theta^1\w\theta^3+\inc{I}{7}\,\theta^2\w\theta^3
    +\inc{\epsilon}{2}\,\theta^2\w\theta^4-\frac{1}{5}\inc{I}{6}\,\theta^3\w\theta^4,\label{e610} \\
  d\theta^4 &=k\,\theta^1\w\theta^2+l\,\theta^1\w\theta^3
    +m\,\theta^1\w\theta^4+\inc{I}{8}\,\theta^2\w\theta^3
   +s\,\theta^2\w\theta^4-\frac{3}{5}\inc{I}{5}\,\theta^3\w\theta^4, \nonumber
 \end{align*}
where $\inv{\epsilon}{2}=\pm 1,0$ is defined through
$$\inv{\epsilon}{2}=\sgn (2F_{qqqq}L_{qq}-3K_{qqq}^3)$$ and
$\inc{I}{5},\inc{I}{6},\inc{I}{7},\inc{I}{8},a,e,f,g,k,l,m,s$ are
functions of $x,y,p,q$. We do not display these functions, since
they are complicated. They can be immediately calculated from
theorem \ref{th.cc.W4d}. We apply the procedure of seeking the
ODEs with four-dimensional symmetry group and insert the results
in table \ref{t.cc.1} on pages \pageref{t.cc.1} --
\pageref{t.cc.2}.

Finally, we have the following geometric description of general
linearizable ODEs.
\begin{proposition}\label{prop.cc.int}
The only smooth third-order ODEs admitting large contact symmetry
group acting intransitively on $\J^2$ are the equations contact
equivalent to  \ben y'''=-2\mu(x)y'+(1-\mu'(x))y, \een with an
arbitrary smooth function $\mu(x)\neq const$.
\begin{proof}
We showed in  corollary \ref{cor.cc.linear} that the above
equations admit a 4-dimensional contact symmetry group acting on
3-dimensional orbits in $\J^2$. All other ODEs admitting a large
contact symmetry group $G$ possess Cartan coframes of rank zero.
These coframes were explicitly constructed in theorem \ref{th.c.1}
for the case $W=0$, $F_{qqqq}=0$, in theorems \ref{th.c.2} and
\ref{th.cc.nW4d} for the case $W\neq 0$ and in theorem
\ref{th.cc.W4d} for the case $W=0$, $F_{qqqq}\neq 0$. The coframe
for the case $W=0$, $F_{qqqq}=0$ generates on $\J^2$ local
structure of the homogeneous space $SP(4,\real)/H_6$. The coframes
of theorems \ref{th.cc.nW4d} and \ref{th.cc.W4d} equip $\J^2$ with
local structure of $G$. In both these cases the action of $G$ on
$\J^2$ is transitive.
\end{proof}
\end{proposition}

\section{Equations with large point symmetry
group}\label{s.pp.large}

\noindent Given the detailed description of the ODEs with large
contact symmetry groups it is already easy to find the equations
admitting large point symmetry groups. It follows from the fact
that any equation possessing point symmetries has at least the
same number of contact symmetries (`number of symmetries' means
the dimension of the symmetry group here.) Therefore the equations
with large point  symmetry groups lie entirely within the classes
with large contact symmetry groups. As a consequence, we must only
do the full reduction for the equations contact equivalent to
those in table \ref{t.cc.1} and analyze existence of point
symmetries. The procedure of reduction for point transformations
is fully analogous to the contact case. We have the following main
branches
\begin{itemize}
\item[i)] Linear and point linearizable equations equivalent to
$y'''=0$ with the 7-dimensional algebra $\co(2,1)\semi{.}\real^3$
of point symmetries. \item[ii)] Non-linearizable equations
admitting the 6-dimensional algebras $\o(2,2)$ or $\o(4)$ of point
symmetries. They are equivalent to $y'''=\tfrac{3y''^2}{2y'}$ or
$y'''=\tfrac{3y'y''^2}{y'^2+1}$ respectively. These classes are
new when compared to the contact classification. \item[iii)] The
linear and point linearizable equations which satisfy $W\neq 0$.
They are equivalent to $y''=-2\mu(x)y'+(1-\mu'(x))y$ and have a
5-dimensional or a 4-dimensional group of symmetries.

\item[iv)] All the remaining equations satisfying $W\neq 0$.

\item[v)] All the remaining equations satisfying $W=0$.
\end{itemize}

Branches i) to iii) are relatively simple to characterize in terms
of point invariants. i) has been described in corollary
\ref{cor.p.7d_flat}. ii), which is contact equivalent to $y'''=0$,
has been discussed in section \ref{s.lor3} of chapter
\ref{ch.point}, and description of iii) is based on the reduction
to a five-dimensional coframe, as in section \ref{s.cc.nw}. After
\ben u_1=\sqrt[3]{W}u_3, \qquad\qquad u_2=\tfrac13Zu_3 \een we
obtain a five dimensional coframe $(\hp{1},\ldots,\hp{4},\vp{})$
with the following basic invariants
\begin{align}
\inc{a}{}=&
\frac{1}{\sqrt[3]{W^2}}\left(K+\frac{1}{18}Z^2+\frac{1}{9}ZF_q-\frac{1}{3}\D
Z\right), \notag \\
\inp{b}{}=&\frac{1}{3u\sqrt[3]{W^2}}\bigg(\left(\frac{1}{12}F_{qq}+\frac{1}{18}Z_q\right)Z^2
+\left(K_q-\frac{1}{3}Z_p-\frac{1}{9}F_qZ_q+\frac{1}{18}F_{qq}F_q\right)Z+\nonumber\\
   &-\frac{1}{6}F_{qq}\D Z-KZ_q+Z_y+\frac32F_{qq}K-3K_p-K_qF_q-F_{qy}\bigg), \nonumber \\
\inp{e}{}=&\frac{1}{u}\left( \frac16F_{qq} -\frac13Z_q
          +\frac{1}{W}\left(\frac29W_qZ-\frac23W_p-\frac29W_qF_q\right)\right),
          \notag \\
\inp{h}{}=&\frac{1}{3u\sqrt[3]{W}}
\bigg(\left(\frac{1}{18}W_qZ^2-\left(\frac{1}{3}W_p+\frac19W_qF_q\right)Z+W_y
-W_qK\right)\frac{1}{W}+\nonumber  \\
 &-3K_q-\frac{1}{3}F_{qq}F_q-F_{qp}\bigg),\nonumber \\
\inp{k}{}=&\frac13\frac{W_q}{\sqrt[3]{W^2}u}. \notag
\end{align}
The linearizable equations are described by conditions expressed
in terms of these invariants, as given in table \ref{t.pp.1}.

The branches iv) and v) need more thorough study. We consider iv)
first. We know from section \ref{s.cc.nw} that any equation with
large point symmetry group in this branch satisfies necessarily
$3WW_{qq}-2W_q^2\neq 0$. It implies $W_q\neq 0$, which allows us
to reduce the last free group parameter via \ben
u_3=\frac{1}{3}\frac{W_q}{\sqrt[3]{W^2}}. \een The set of basic
invariants is the following
\begin{align*} \inp{I}{1}=&-3\frac{W_{qq}W}{W_q^2}, \notag \\
\inp{I}{2}=&\frac{1}{\sqrt[3]{W}W_q}\left(3W_p+W_qF_q-W_qZ-3WZ_q-3F_{qq}W
\right), \notag \\
\inp{I}{3}=&-\frac{\sqrt[3]{W^2}}{2W_q}\left(2Z_q+F_{qq}\right),
\\
\inp{I}{4}=&\frac{1}{18\sqrt[3]{W^2}}\left(Z^2-6\D
Z+18K+2ZF_q\right),\notag
\end{align*} and the above invariants are (up to constant numbers) the
$T^1_{13}$, $T^1_{14}$, $T^2_{13}$, and $T^2_{14}$ coefficients in
the structural equations \ben \der
\hp{i}=\frac12T^i_{jk}\hp{j}\w\hp{k},\qquad T^i_{jk}=-T^i_{kj}\een
for the coframe. The ODEs with large point symmetry groups which
fall into this branch are types II.2, II.3, IV, and VI of table
\ref{t.pp.1}.

We turn to the branch v). The coframe is given either by \ben
 u_1=-\frac{3 F_{qqq}^5}{4 F_{qqqq}^3},\qquad  u_2=\frac{F_{qqq}^2N}{2F_{qqqq}},
     \qquad u_3=\frac{F_{qqq}^2}{2F_{qqqq}}, \een
provided that $F_{qqqq}\neq0$ or, if $F_{qqqq}=0$ but
$F_{qqq}\neq0$,  by
\begin{align*}
u_1=&-\frac{1}{36F_{qqq}^4}\left(6F_{qqqp}+5F_{qqq}F_{qq}\right)^3,\\
u_2=&-\frac{1}{6F_{qqq}^2}\left(6F_{qqqp}+5F_{qqq}F_{qq}\right)N,\\
u_3=&-\frac{1}{6F_{qqq}}\left(6F_{qqqp}+5F_{qqq}F_{qq}\right),
\end{align*} where \ben
N=F_{qqp}+\frac{1}{6}F_{qq}^2+\frac{1}{3}F_{qqq}F_{q}. \een For
$F_{qqqq}\neq 0$ we have the following basic invariants
\begin{align*}
\inp{I}{5}=&\frac{F_{qqq}F_{qqqqq}}{F_{qqqq}^2}, \notag \\
\inp{I}{6}=&\frac{F_{qqqq}}{F_{qqq}^4}\left(\frac83F_{qqqq}-12F_{qqq}K_{qqq}+\frac59F_{qqqq}F_{qq}^2+20F_{qqqq}K_{qq} \right), \notag \\
\inp{I}{7}=&\frac{F_{qqqq}}{F_{qqq}^4}\left(
6N_qF_{qqq}-6NF_{qqqq}+F_{qq}F_{qqq}^2\right)  \\
\inp{I}{8}=&-\frac{2}{27}\frac{F_{qqqq}^4}{F_{qqq}^8}\big(4NF_qF_{qqq}+6\D
NF_{qqq}+\notag \\
&-9N^2-F_{qq}^2N-36K_{qq}N-6F_{qqq}^2K\big), \notag
\end{align*}
which are obtained from the coefficients $T^1_{13}, T^1_{14},
T^2_{13}$ and $T^2_{14}$ in the structural equations.

For the branch $F_{qqqq}=0$, $F_{qqq}\neq0$, which contains only
one class of equations with large point symmetry groups, the
equations equivalent to $F=q^3$, we have the invariants \ben
\inp{I}{9}=T^1_{12}, \qquad \inp{I}{10}=T^1_{14}, \qquad
\inp{I}{11}=T^2_{14}.\een

Properties of the ODEs admitting a large point symmetry group are
given in table \ref{t.pp.1} on pages \pageref{t.pp.1} --
\pageref{t.pp.2}. In the point classification we also have a
counterpart of proposition \ref{prop.cc.int}.
\begin{remark}\label{rem.fp}
The fibre-preserving classification of the ODEs admitting a large
symmetry group is parallel to the point classification and has
been already done \cite{God1, Grebot}. The main difference is that
types I.3, II.3. and IX do not admit four-dimensional
fibre-preserving symmetry groups.
\end{remark}

\section{Fibre-preserving equivalence to certain reduced Chazy equations}
\noindent An ordinary differential equation in the {\em complex}
domain is said to have the Painlev\'e property if its general
solution does not have movable branch points, that is branch
points whose location depends on integration constants,
\cite{Cos}. The problem of classifying the third-order Painlev\'e
ODEs which are polynomials in $y,y'$, and $y''$ and are locally
analytic in $x$ was studied by J. Chazy \cite{Cha}, who considered
the polynomial equations modulo the following transformations \ben
x\to \chi(x),\qquad\qquad y\to \alpha(x)y+\beta(x), \een which are
a subclass of complex analytic fibre-preserving transformations.
J. Chazy found thirteen classes of the ODEs satisfying the
Painleve property. Each of these classes has a particularly simple
representative -- the reduced Chazy class -- obtained by a certain
limit procedure. Here we are interested in reduced Chazy classes
II, IV, V, VI, VII and XI $\sigma\neq11$. They are as follows
\begin{flalign*}
II: && F=&-2yq-2p^2,&&&\\
IV: && F=&-3yq-3p^2-3y^2p,&&&\\
V: &&  F=&-2yq-4p^2-2y^2p,&&&\\
VI: && F=&-yq-5p^2-y^2p,&&&\\
VII: && F=&-yq-2p^2+2y^2p,&&&\\
XI: &&
F=&-2yq-2p^2+\frac{24}{\sigma^2-1}\left(p+y^2\right)^2,&&\sigma\in\mathbf{N},\
\sigma\neq 1,6k.&
\end{flalign*}
All of them have the form
\be
 F=\kappa yq+\lambda p^2+ \mu y^2p+\nu y^4,\label{e.ch.e}
\ee with some constant numbers $\kappa$, $\lambda$, $\mu$, $\nu$.
We exclude type XI for $\sigma=11$ for technical reasons.

We aim to find necessary and sufficient conditions for a {\em
regular real} third order ODE to be fibre-preserving equivalent to
one of the above equations. We find the Cartan coframe for such
equations and the explicit formulae for the fibre-preserving
invariants. Next we find the functional relations between the
invariants, which allows us in to describe the classifying
function $\T$ explicitly and consequently describe its image.

First step is calculating the structural equations for the
fibre-preserving coframe of theorem \ref{th.f.1} for the Chazy
types. They are as follows
\begin{align}
 \der\hf{1} =&\vf{1}\w\hf{1}+\hf{4}\w\hf{2},\nonumber \\
 \der\hf{2} =&\vf{2}\w\hf{1}+\vf{3}\w\hf{2}+\hf{4}\w\hf{3},\nonumber \\
 \der\hf{3}=&\vf{2}\w\hf{2}+(2\vf{3}-\vf{1})\w\hf{3}+\inf{A}{1}\hf{4}\w\hf{1},\nonumber \\
 \der\hf{4} =&(\vf{1}-\vf{3})\w\hf{4}, \nonumber\\
 \der\vf{1} =&-\vf{2}\w\hf{4}+(2\inf{C}{1}-\inf{A}{2})\hf{1}\w\hf{4}, \nonumber \\
 \der\vf{2}=&(\vf{3}-\vf{1})\w\vf{2}+\inf{A}{3}\hf{1}\w\hf{4}
  +(\inf{C}{1}-\inf{A}{2})\hf{2}\w\hf{4}, \nonumber \\
 \der\vf{3}=&(\inf{C}{1}-\inf{A}{2})\hf{1}\w\hf{4}. \nonumber
\end{align}
In this system the functions $\inp{B}{i}$, $i=1\ldots6$,
$\inp{D}{1}$, and $\inp{D}{2}$ vanish, which is equivalent to the
following three fibre-preserving invariant conditions
\begin{equation}
F_{qq}=0,\qquad F_{qpp}=0, \qquad
F_{ppp}=2F_{qpy}-\frac23F^2_{qp}.\label{e.ch.10}
\end{equation}
The most important non-vanishing invariants read
\begin{align*}
 &\inf{A}{2}=\frac{1}{u_1^3}\left(\left(\lambda-\frac{7\kappa}{6}\right )u_2u_3
   +\left(\mu+\frac{\kappa\lambda}{3}+\frac{\kappa^2}{18}\right
   )yu_3^2\right ), \notag \\
 &\inf{C}{1}-\inf{A}{2}=\frac{\kappa u_3}{3u_1^2}, \\
 &X_4(\inf{C}{1}-\inf{A}{2})=-\left ( \frac23\kappa u_2u_3+\frac19\kappa^2 yu_3^2\right)
 \frac{1}{u_1^3}.\notag
\end{align*}
For the Chazy classes one can reduce the parameters through
\be\label{e.ch.20} \inf{A}{2}=1, \quad
\inf{C}{1}-\inf{A}{2}=\frac13,\quad
X_4(\inf{C}{1}-\inf{A}{2})=0.\ee This leads to the
four-dimensional system
\begin{equation}
 \begin{aligned}
  d\theta^1 &=\inw{a}\,\theta^1\w\theta^4-\theta^2\w\theta^4,\\
  d\theta^2 &=-2\tau\,\theta^1\w\theta^2+\inw{b}\,\theta^1\w\theta^4
  +2\inw{a}\,\theta^2\w\theta^4-\theta^3\w\theta^4,\\
  d\theta^3 &=\left (\tfrac{\lambda}{\kappa}-\tfrac{2}{3}\right )\theta^1\w\theta^2
  -3\tau\,\theta^1\w\theta^3
    +\inw{c}\in\,\theta^1\w\theta^4+\inw{b}\,\theta^2\w\theta^4+3\inw{a}\,\theta^3\w\theta^4, \\
  d\theta^4 &=\tau\,\theta^1\w\theta^4,
 \end{aligned}\label{e.ch.dth}
 \end{equation}
where \ben
 \tau=\frac{\mu}{\kappa^2}+\frac{\lambda}{6\kappa}+\frac{1}{4},
\een and $\inw{a},\,\inw{b},\,\inw{c}$ are functions on $\J^2$. We
check by direct calculations that in this case the coframe is of
order one and all the invariants are generated by $\inw{a}$ and
$\inw{a}_4=X_4(\inw{a})$ by the following formulae
\begin{align}\label{e.ch.syz}
   \inw{b}=&\frac{1}{\tau}\left (\frac{1}{3}-\frac{\lambda}{\kappa}\right)\inw{a}-\frac{1}{2\tau}
     +\frac{1}{12\tau^2}\left (\frac{1}{3}-\frac{\lambda}{\kappa} \right ),\notag \\
   \inw{c}=&\frac{1}{\tau}\left (\frac{\lambda}{\kappa}-\frac{7}{6} \right)\inw{a}_4
     +\frac{1}{2\tau}\left (\frac{\lambda}{\kappa}-\frac{7}{6} \right)\inw{a}^2
     +\left (-\frac{1}{\tau}+\frac{1}{6\tau^2}\left (\frac{\lambda}{\kappa}-\frac{7}{6} \right) \right )\inw{a}
     -\frac{1}{2\tau^2} \notag\\
     &+\frac{1}{36\tau^3}\left
     (\frac{\lambda}{\kappa}-\frac{144\nu}{\kappa^3}+\frac{3}{2}\right),\notag \\
   \inw{a}_1=&-2\tau \inw{a}-\frac{1}{6},\qquad \inw{a}_2=\tau, \qquad \inw{a}_3= 0, \\
   \inw{a}_{41}=& -3\tau \inw{a}_4+2\tau \inw{a}^2+\left (\frac{\lambda}{\kappa}-\frac{1}{6} \right )\inw{a}
          +\frac{1}{12\tau}\left (\frac{\lambda}{\kappa}-\frac{1}{3} \right )+\frac{1}{2},\notag \\
   \inw{a}_{42}=& -4\tau \inw{a}-\frac{1}{6},\qquad \inw{a}_{43}=\tau, \notag\\
   \inw{a}_{44}=& -7\inw{a}_4 \inw{a}-\frac{\inw{a}_4}{6\tau}-6\inw{a}^3+\frac{1}{\tau}\left (\frac{\lambda}{\kappa}-1 \right ) \inw{a}^2
   +\left (\frac{1}{\tau}+\frac{1}{6\tau^2}\left (\frac{\lambda}{\kappa}-\frac{1}{2} \right )\right )\inw{a}
   \notag\\
       &+\frac{1}{6\tau^2}+\frac{1}{\tau^3}\left (\frac{\nu}{\kappa^3}-\frac{1}{72} \right ),\notag
 \end{align}
where, as usual, $\inw{a}_i=X_i(\inw{a})$,
$\inw{a}_{ij}=X_j(X_i(\inw{a}))$. These algebraic formulae, when
differentiated, enable to express all other derivatives of
$\inw{a}$, $\inw{b}$ and $\inw{c}$ in terms of $\inw{a}$ and
$\inw{a}_4$. In order to do this we only must consecutively
substitute the coframe derivatives of $\inw{a}$, $\inw{b}$ and
$\inw{c}$ with the right hand side of \eqref{e.ch.syz}, for
instance
$$\inw{a}_{11}=-2\tau \inw{a}_1=4\tau^2 \inw{a}+\frac{\tau}{3},
\qquad \qquad \inw{a}_{12}=-2\inw{a}_2=-2\tau^2,$$ etc. Thereby
the non-constant components of the classifying function
$\T\colon\J^2\to\real^N$ are completely characterized by
$$(x,y,p,q)\mapsto(\inw{a},\inw{b},\inw{c},\inw{a}_1,\inw{a}_2,\inw{a}_3,
\inw{a}_4,\inw{a}_{41},\inw{a}_{42},\inw{a}_{43},\inw{a}_{44}).
$$
The graph of this function in $\real^{11}$ is parameterized by
$\inw{a}$ and $\inw{a}_4$.

Let us consider an arbitrary third-order ODE. It is locally
fibre-preserving equivalent to one of the Chazy classes if and
only if graphs of respective classifying functions overlap. This
is only possible if i) conditions \eqref{e.ch.10} are satisfied,
ii) the reduction defined by the conditions \eqref{e.ch.20} is
possible, and iii) after the reduction the equations
\eqref{e.ch.dth} and \eqref{e.ch.syz} hold. The reduction
\eqref{e.ch.20} is possible iff
 \be\label{e.ch.30}\bal
P=&\D F_{qp}-F_{qy}\neq0, \\
Q=&2W_p-\D W_q+F_qW_q\neq0. \eal \ee After the reduction to
dimension four given by \ben
 u_1 =\frac{2P^2}{Q}, \qquad u_3 =\frac{2F_q P^3}{3Q}-\frac{2P^2\D P}{Q},
 \qquad  u_3=-\frac{4P^3}{Q}
\een we get the frame \be\label{e.ch.frame}\begin{aligned}
   X_1=&\frac{Q}{P^2}\left (\frac{Q}{4W_q}-\frac{F_q}{6} \right)\partial_p
      +\frac{Q}{P^2}\left( \frac{Q^2}{16W_q^2}-\frac{F_q^2}{36}-\frac{K}{2}\right )\partial_q, \\
   X_2=&-\frac{Q^2}{4P^3}\partial_p-\frac{Q^3}{8P^3W_q}\partial_q, \\
   X_3=&\frac{Q^3}{8P^4}\partial_q, \\
   X_4=&-\frac{2P}{Q}\D.
 \end{aligned}\ee
 Equations \eqref{e.ch.dth} are satisfied if and only if
\begin{align}\label{e.ch.40}
   &2 W_{pp}-W_{qy}+F_{qp}W_q=0,\notag\\
   &W_q\D P-P\D W_q=0,\notag\\
   &P_y+\frac{1}{3}PF_{qp}=0,\notag\\
   &Q_y+\frac{1}{3}QF_{qp}-2\tau P^2=0,\\
   &K_p+\frac{1}{2}F_{qy}-\frac{5}{36}F_qF_{qp}+\frac{1}{W_q}
   \left (F_{qp}\D W_q-\frac{1}{12}F_qW_{qy}+\frac{1}{2}\D W_{qy}\right)\notag\\
   &-\frac{3W_{qy}\D W_q}{4W_q^2}+\left (\frac{2}{3}-\frac{\lambda}{\kappa}\right)P=0.\notag
 \end{align}
Finally, equations \eqref{e.ch.syz} must be satisfied by functions
\be\label{e.ch.inv}\begin{aligned}
 \inw{a}=&\frac{P}{W_q}+\frac{1}{Q}\left (4\D P-\frac{2}{3}F_qP -\frac{2PW_p}{W_q}\right )-\frac{2P\D Q}{Q^2},\\
 \inw{b}=&\left \{ \frac{5}{2W_q^2}+\frac{1}{Q}\left(-\frac{10F_q}{3W_q}-\frac{10W_p}{W_q^2}\right)
     +\frac{1}{Q^2}\left(-4F_p-4K-\frac{2}{9}F_q^2+\right.\right.  \\
    &\left.\left.+\frac{1}{W_q}\left(\frac{20}{3}W_pF_q-4\D W_p+2\D Q\right )
    +\frac{10W_p^2}{W_q^2}\right)\right\}P^2,\\
 \inw{c}=&\frac{8P^3W}{Q^3},
\end{aligned}
\ee and by their derivatives $\inw{a}_4$,
$\inw{a}_{41},\ldots,\inw{a}_{44}$ with respect to the frame
$X_1,\ldots,X_4$. Therefore, by means of theorem \ref{th.cc.clas}
we have
\begin{proposition}
An ODE is locally fibre-preserving equivalent to one of the Chazy
classes II, IV, V, VI, VII or XI for $\sigma\neq11$ in a
neighbourhood of a point $j_0\in\J^2$ if and only if i) the ODE
satisfies the conditions \eqref{e.ch.10}, \eqref{e.ch.30},
\eqref{e.ch.40}, and \eqref{e.ch.syz} with the invariants
$\inw{a}$, $\inw{b}$, $\inw{c}$, $\inw{a}_4$, and
$\inw{a}_{41},\ldots,\inw{a}_{44}$ given by \eqref{e.ch.frame} and
\eqref{e.ch.inv}, ii) the values of $\inw{a}(j_0)$ and
$\inw{a}_4(j_0)$ for the ODE and the Chazy class coincide.
\end{proposition}

 Given these criteria for the equivalence it is
interesting to find a transformation of variables transforming an
ODE
$$\frac{\der^3 y}{\der x^3}=F\left(x,y,\frac{\der y}{\der x},\frac{\der^2 y}{\der x^2}\right)$$
to its equivalent Chazy class \ben
 \frac{\der^3 \cc{y}}{\der \cc{x}^3}=\kappa \cc{y}\frac{\der^2 \cc{y}}{\der \cc{x}^2}
 +\lambda \left(\frac{\der \cc{y}}{\der \cc{x}}\right)^2
 +\mu \frac{\der \cc{y}}{\der \cc{x}}\cc{y}^2+\nu \cc{y}^4.
\een We show that the transformation \ben
\cc{y}=\phi(x,y),\qquad\qquad \cc{x}=\chi(x)
 \een
 may be easily found.

Let us apply to a Chazy class the above (arbitrary)
fibre-preserving transformation and calculate for so obtained ODE
(which is a general ODE equivalent to the Chazy classes) the
quantities $P$, $Q$ and $F_q-pF_{qp}$:
 \begin{align}\label{r1670}
 P=&-\kappa\chi_x\phi_y,\notag \\
 Q=&2\kappa^2\tau\,\chi_x^2\,\phi_y\,\phi,\\
 F_q-pF_{qp}=&\kappa\chi_x\phi+3\frac{\chi_{xx}}{\chi_x}-3\frac{\phi_{xy}}{\phi_y}.\notag
 \end{align}
 From first and second equation we get
 \begin{align*}
 (\log |\phi|)_y=&2\tau\frac{P^2}{Q},\\
 \chi_x=&-\frac{Q}{2\tau P \phi}.
 \end{align*}
Putting this into third equation of \eqref{r1670} we obtain \ben
 (\log |\phi|)_x=\frac{1}{2}\left(\log\left|\frac{Q^2}{P^3}\right|\right)_x-\frac{\kappa Q}{12\tau P}
  +\frac{1}{6}(pF_{qp}-F_q).
\een Finally, after integration of $\chi_x$ and $\phi_x$ we have
\be \label{e.ch.tr}\bal &\bar{y}=\frac{c_1Q}{|P|^{\frac{3}{2}}}
  \exp\left\{\int_{x_0}^x\left(-\frac{\kappa Q}{12\tau P}
   +\frac{1}{6}(pF_{qp}-F_q)\right)dx+2\tau\int_{y_0}^y\left.\frac{P^2}{Q}\right|_{x=x_0}dy\right\},\\
&\bar{x}=\frac{1}{2\tau}\int^x_{x_0}\frac{Q}{P\phi}dx+c_2. \eal\ee
We summarize this calculation as follows.
\begin{proposition}
If there exists a fibre-preserving transformation from an equation
$y'''=F(x,y,y',y'')$ to a reduced Chazy type II, IV -- VII or XI
$\sigma \neq 11$, then it is given by the inverse of eq.
\eqref{e.ch.tr}, where $P$ and $Q$ are calculated for
$y'''=F(x,y,y',y'')$ according to formulae \eqref{e.ch.30}.
\end{proposition}

\ifthenelse{\boolean{@twoside}}{%
\fancyhead[CO]{\iffloatpage{\headfont{tables}}{\rightmark}}
\fancyhead[CE]{\iffloatpage{}{\leftmark}}}{%
\fancyhead[CO]{\iffloatpage{}{\leftmark}}}

\addtocounter{table}{1}
\begin{sidewaystable}
 \caption{Equations admitting large contact
symmetry groups}\label{t.cc.1}
\renewcommand*{\arraystretch}{1.5}
\begin{tabular}
{|c|c|c|c|c|c|c|c|}
  \hline
  & \multicolumn{2}{|c|}{Equation} & \multicolumn{3}{|c|}{Characterization}
  & Symmetry algebra $\g$ & $\dim\g$ \\ \hline \hline

 I & $F=0$ & & \multicolumn{3}{|c|}{$W=0,\qquad F_{qqqq}=0$}  & $\o(3,2)$ & $10$ \\  \hline

 \multirow{3}{*}{II} & \multirow{3}{*}{$F=-2\mu p+y$} & \multirow{3}{*}{$\mu\in\real$} &
 \multicolumn{3}{|c|}{\multirow{3}{*}{$W\neq 0,\quad\inc{a}{}=\mu,\quad\inc{k}{}=0$}} &
 $[V_1,V_4] =-\mu V_2+V_3,\, [V_1,V_5] =V_1$, & \multirow{3}{*}{$5$} \\
 & & & \multicolumn{3}{|c|}{} & $[V_2,V_4] =V_1-\mu V_3,\, [V_2,V_5] =V_2$, & \\
 & & & \multicolumn{3}{|c|}{} & $[V_3,V_4] =V_2,\, [V_3,V_5] =V_3$, & \\ \hline

\multirow{2}{*}{III} & \multirow{2}{*}{$F=-2\mu(x)p+(1-\mu'(x))y$}
& & \multicolumn{3}{|c|}{\multirow{2}{*}{$W\neq 0,\quad
\inc{a}{}=\mu(x),
\quad\inc{b}{}=\inc{e}{}=\inc{h}{}=\inc{k}{}=0$}} &
$[V_1,V_4]=V_1,\, [V_2,V_4]=V_2$, & \multirow{8}{*}{$4$} \\
& & & \multicolumn{3}{|c|}{} & $[V_3,V_4]=V_3$, & \\ \cline{1-7}

\multirow{2}{*}{IV} & \multirow{2}{*}{$F=q^{3/2+1/(2\sqrt{\mu})}$}
& \multirow{2}{*}{$\mu>0,\,\neq\tfrac{1}{9}$} &
\multirow{2}{*}{$0<1+4\inv{\epsilon}{1}(\inc{I}{1})^{-2}\neq\tfrac{1}{9}$}
& \multirow{2}{*}{$\mu=1+4\inv{\epsilon}{1}(\inc{I}{1})^{-2}$} &
\multirow{6}{*}{\rotatebox{270}{$\inc{I}{1},\inc{I}{2},\inc{I}{3},\inc{I}{4}=const$}
\rotatebox{270}{$W\neq 0,\quad\inv{\epsilon}{1}\neq 0,$}}  &
{$[V_1,V_4]=2V_1,\, [V_2,V_4]=(1+\sqrt{\mu})V_2,$} & \\
& & & & & & {$[V_2,V_3]=V_1, \,[V_3,V_4]=(1-\sqrt{\mu})V_3$} & \\
\cline{1-5} \cline{7-7}

\multirow{2}{*}{V} & \multirow{2}{*}{$F=(q^2+1)^{\frac{3}{2}}\,
   \exp\left (\frac{{\rm arc\,tg}\,q}{\sqrt{\mu}}\right)$} &
   \multirow{2}{*}{$\mu>0$} & \multirow{2}{*}{$\inv{\epsilon}{1}=-1,\,1-4(\inc{I}{1})^{-2}<0$} &
\multirow{2}{*}{$\mu=4(\inc{I}{1})^{-2}-1$} &  & $[V_1,V_4]=2V_1,\, [V_2,V_4]=V_2-\sqrt{\mu}V_3$ & \\
& & & & & & $[V_2,V_3]=V_1,\, [V_3,V_4]=\sqrt{\mu} V_2+V_3$ & \\
\cline{1-5} \cline{7-7}

\multirow{2}{*}{VI} & \multirow{2}{*}{$F=\exp q$} &  &
\multirow{2}{*}{$\inv{\epsilon}{1}=-1,\,\inc{I}{1}=-2$} &
 &  & $[V_1,V_4]=2V_1,\, [V_2,V_4]=V_2+V_3$ & \\
& & & & & & $[V_2,V_3]=V_1,\, [V_3,V_4]=V_3$ & \\ \hline

\end{tabular}
\end{sidewaystable}

\newpage
\ifthenelse{\boolean{@twoside}}{\fancyhead[CE]{\leftmark}}{\fancyhead[CO]{\leftmark}}

\addtocounter{table}{-1}
\begin{sidewaystable}
\label{t.cc.2}\caption{Equations admitting large contact symmetry
groups}
\renewcommand*{\arraystretch}{1.5}
\begin{tabular}
{|c|c|c|c|c|c|c|c|}
  \hline
  & \multicolumn{2}{|c|}{Equation} & \multicolumn{3}{|c|}{Characterization} & Symmetry algebra $\g$
  & $\dim\g$ \\ \hline

\multirow{2}{*}{VII}  &
$F=\mu\left(\frac{q^2}{1-p^2}-p^2+1\right)^{3/2}$ &
\multirow{2}{*}{$\mu>0$} & $\inv{\epsilon}{2}=1$, &
\multirow{4}{*}{$\mu=\sqrt{\left|\frac{9(\inc{I}{7})^3}{9(\inc{I}{7})^3-2}\right|}$}
& \multirow{14}{*}{\rotatebox{270}{$W=0,\quad F_{qqqq}\neq 0,\quad
\inc{I}{5},\inc{I}{6},\inc{I}{7},\inc{I}{8}=const$}}
& \multirow{2}{*}{$\u(2)$} & \multirow{14}{*}{$4$} \\
 & $-3\frac{q^2p}{1-p^2}+p^3-p^2$ &  &  $0<\inc{I}{7}<\frac{\sqrt[3]{6}}{3}$, & & & & \\[0.5ex]
 \cline{1-4} \cline{7-7}

 \multirow{3}{*}{VIII} & \multirow{3}{*}{$F=\mu\frac{(2qy-p^2)^{3/2}}{y^2}$}  & $0<\mu<1$
 & $\inv{\epsilon}{2}=1,\,\inc{I}{7}<0$
&
& & \multirow{8}{*}{$\gl(2,\real)$} & \\ \cline{3-4}
   &  &  $\mu>1$ & $\inv{\epsilon}{2}=-1,\,\inc{I}{7}>\frac{\sqrt[3]{6}}{3}$ & & & & \\  \cline{3-5}
   &  &  $\mu=1$ & $\inv{\epsilon}{2}=\inc{I}{7}=0,\,\inc{I}{8}=1$ & & & & \\ \cline{1-5}

\multirow{2}{*}{IX} &
\multirow{2}{*}{$F=4\mu(q-p^2)^{3/2}+6qp-4p^3$} &
\multirow{2}{*}{$\mu>0$} &
\multirow{2}{*}{$\inv{\epsilon}{2}=-1,\,0<\inc{I}{7}<\frac{\sqrt[3]{6}}{3}$}
&
\multirow{4}{*}{$\mu=\sqrt{\left|\frac{9(\inc{I}{7})^3}{9(\inc{I}{7})^3-2}\right|}$} & & & \\
& & & & & & & \\ \cline{1-4}

\multirow{3}{*}{X} &
\multirow{3}{*}{$F=\mu\left(\frac{q^2}{p^2}+p^2\right)^{3/2}+3\frac{q^2}{p}+p^3$}
& $0<\mu<1$ & $\inv{\epsilon}{2}=-1,\,\inc{I}{7}<0$
&
& & & \\
\cline{3-4}
   &  &  $\mu>1$ & $\inv{\epsilon}{2}=1,\,\inc{I}{7}>\frac{\sqrt[3]{6}}{3}$ & & & & \\  \cline{3-5}
   &  &  $\mu=1$ & $\inv{\epsilon}{2}=\inc{I}{7}=0,\,\inc{I}{8}=-1$ & & & & \\ \cline{1-5} \cline{7-7}

\multirow{2}{*}{XI} & \multirow{2}{*}{$F=(q^2+1)^{3/2}$} &&
\multirow{2}{*}{$\inv{\epsilon}{2}=1,\,\inc{I}{7}=\frac{\sqrt[3]{6}}{3}$} &&& $[V_2,V_4]=-V_3,\,\,[V_2,V_3]=V_1,$ & \\
 &&&&&& $[V_3,V_4]=V_2$ & \\ \cline{1-5} \cline{7-7}

\multirow{2}{*}{XII} & \multirow{2}{*}{$F=q^{3/2}$} &&
\multirow{2}{*}{$\inv{\epsilon}{2}=-1,\,\inc{I}{7}=\frac{\sqrt[3]{6}}{3}$} &&& $[V_2,V_4]=V_2,\,\,[V_2,V_3]=V_1,$ & \\
 &&&&&& $[V_3,V_4]=-V_3$ & \\ \hline
\end{tabular}
\end{sidewaystable}

\begin{sidewaystable}
\caption{Equations admitting large point symmetry
groups}\label{t.pp.1}
\renewcommand*{\arraystretch}{1.5}
\begin{tabular}
{|c|c|c|c|c|c|c|c|}
  \hline
  & \multicolumn{2}{|c|}{Equation} & \multicolumn{3}{|c|}{Characterization}
  & Symmetry algebra $\g$ & $\dim\g$
  \\ \hline \hline
\multirow{2}{*}{I.1} & \multirow{2}{*}{$F=0$} & &
\multirow{2}{*}{$F_{qq}^2+6F_{qqp}=0$} &
\multicolumn{2}{|c|}{\multirow{4}{*}{$F_{qq}K+\frac{1}{3}F_qF_{qp}-F_{qy}+\frac12F_{pp}=0$}}
 & \multirow{2}{*}{$\co(2,1)\semi{.}\real^3$}
 & \multirow{2}{*}{$7$} \\ & & & & \multicolumn{2}{|c|}{}  & & \\
 \cline{1-4} \cline{7-8}

\multirow{2}{*}{I.2} & \multirow{2}{*}{$F=\frac32\frac{q^2}{p}$} &
& \multirow{2}{*}{$F_{qq}^2+6F_{qqp}<0$} &
\multicolumn{2}{|c|}{\multirow{4}{*}{$W=0,\qquad F_{qqq}=0$}} &
\multirow{2}{*}{$\o(2,2)$}
 & \multirow{2}{*}{$6$} \\
& & & & \multicolumn{2}{|c|}{} & & \\
 \cline{1-4} \cline{7-8}

\multirow{2}{*}{I.3} & \multirow{2}{*}{$F=\frac{3q^2p}{1+p^2}$} &
&\multirow{2}{*}{$F_{qq}^2+6F_{qqp}>0$} & \multicolumn{2}{|c|}{} &
\multirow{2}{*}{$\o(4)$} & \multirow{2}{*}{$6$} \\
 & & & & \multicolumn{2}{|c|}{} & & \\ \hline

\multirow{2}{*}{I.4} & \multirow{2}{*}{$F=q^3$} &
\multirow{2}{*}{} &
\multicolumn{3}{|c|}{\multirow{2}{*}{$W=F_{qqqq}=0,\quad
\inp{I}{9}=-\frac25,\quad \inp{I}{10}=\frac{1}{25},\quad
\inp{I}{11}=0 $}} &
{$[V_1,V_4]=2V_1,\, [V_2,V_4]=\tfrac43 V_2,$} & \multirow{2}{*}{4}\\
& & & \multicolumn{3}{|c|}{} & {$[V_2,V_3]=V_1, \,[V_3,V_4]=-\tfrac23 V_3$} & \\
\hline

\multirow{2}{*}{II.1} & \multirow{2}{*}{$F=-2\mu p+y$} &
\multirow{2}{*}{$\mu\in\real$} &
\multicolumn{3}{|c|}{\multirow{2}{*}{$W\neq 0,\quad
\inc{a}{}=\mu,\quad \inp{e}{}=\inp{k}{}=0 $}} &
\multirow{2}{*}{as in table \ref{t.cc.1}} & \multirow{2}{*}{5} \\
& & & \multicolumn{3}{|c|}{} &  & \\ \hline

\multirow{2}{*}{II.2} & \multirow{2}{*}{$F=\mu\frac{q^2}{p}$} &
\multirow{2}{*}{$\mu > \tfrac32,\neq3$} &
\multirow{2}{*}{$\inp{I}{2}\notin[\,0,\sqrt[3]{4}\,]$} &
\multirow{2}{*}{$\inp{I}{2}=\sqrt[3]{\frac{(2\mu-3)^2}{\mu(\mu-3)}}$}
& \multirow{2}{*}{\raisebox{-2ex}{$W\neq0,\,\,\inp{I}{1}=-2$}} &
\multirow{2}{*}{$[V_1,V_2]=V_1,\, [V_3,V_4]=V_3$} & \multirow{3}{*}{4}\\
& & & & & &  & \\  \cline{1-5}\cline{7-8}

\multirow{2}{*}{II.3} &
\multirow{2}{*}{$F=\frac{3p+\mu}{p^2+1}q^2$} &
\multirow{2}{*}{$\mu>0$} &
\multirow{2}{*}{$\inp{I}{2}\in(0,\sqrt[3]{4})$} &
\multirow{2}{*}{$\inp{I}{2}=\sqrt[3]{\frac{4\mu^2}{\mu^2+9}}$} &
\multirow{2}{*}{\raisebox{2ex}{$\inp{I}{2},\,\inp{I}{3},\,\inp{I}{4}=const$}}
 & $[V_1,V_2]=V_3,\, [V_3,V_1]=V_2$ & \multirow{3}{*}{4}\\
& & & & & & $[V_3,V_4]=V_3,\, [V_2,V_4]=V_2$ & \\ \hline
\end{tabular}
\end{sidewaystable}

\begin{sidewaystable}
\addtocounter{table}{-1} \label{t.pp.2} \caption{Equations
admitting large point symmetry groups}
\renewcommand*{\arraystretch}{1.5}
\begin{tabular}
{|c|c|c|c|c|c|c|c|}
  \hline
  & \multicolumn{2}{|c|}{Equation} & \multicolumn{3}{|c|}{Characterization}
  & Symmetry algebra $\g$ & $\dim\g$
  \\ \hline

\multirow{2}{*}{III} & \multirow{2}{*}{$F=-2\mu(x)p+(1-\mu'(x))y$}
& & \multicolumn{3}{|c|}{\multirow{2}{*}{$W\neq 0,\quad
\inc{a}{}=\mu(x), \quad \inp{b}{}=\inp{e}{}=\inp{h}{}=\inp{k}{}=0
$}} &
\multirow{12}{*}{as in table \ref{t.cc.1}} & \multirow{12}{*}{4} \\
& & & \multicolumn{3}{|c|}{} &  & \\ \cline{1-6}

\multirow{2}{*}{IV} & \multirow{2}{*}{$F=q^{\mu}$} &
\multirow{2}{*}{$\mu\neq0,1,\tfrac32,3$} &
\multirow{2}{*}{$\inp{I}{1}\neq-3$} &
\multirow{2}{*}{$\inp{I}{1}=\frac{4-3\mu}{\mu-1}$}
&\multirow{2}{*}{\raisebox{-1.7ex}{$W\neq0$}} &
 &  \\
& &  & & & & & \\
\cline{1-5}

\multirow{2}{*}{VI} & \multirow{2}{*}{$F=\exp q$} & &
\multicolumn{2}{|c|}{\multirow{2}{*}{$\inp{I}{1}=-3$}}
 & \multirow{2}{*}{\raisebox{1.7ex}{$\inp{I}{1},\inp{I}{2},\,\inp{I}{3},\,\inp{I}{4}=const$}}
 &  &  \\
& & & \multicolumn{2}{|c|}{} & &  & \\ \cline{1-6}

\multirow{2}{*}{VIII} &
\multirow{2}{*}{$F=\mu\frac{(2qy-p^2)^{3/2}}{y^2}$}  &
\multirow{2}{*}{$\mu>0$} & \multirow{2}{*}{$\inp{I}{8}>-\tfrac32$}
&\multirow{2}{*}{$\inp{I}{8}=-\tfrac32+\tfrac{2}{\mu^2}$} &
\multirow{3}{*}{\raisebox{-2ex}{$F_{qqqq}\neq 0,\,\,W=0$}} &
& \\
 & & & & & & & \\  \cline{1-5}

\multirow{2}{*}{IX} &
\multirow{2}{*}{$F=4\mu(q-p^2)^{3/2}+6qp-4p^3$} &
\multirow{2}{*}{$\mu>0$} & \multirow{2}{*}{$\inp{I}{8}<-\tfrac32$}
&\multirow{2}{*}{$\inp{I}{8}=-\tfrac32-\tfrac{2}{\mu^2}$} &
\multirow{3}{*}{} & &  \\
 & & & & &
 \multirow{3}{*}{\raisebox{3ex}{$\inp{I}{5},\inp{I}{6},\,\inp{I}{7},\,\inp{I}{8}=const$}} & &
 \\  \cline{1-5}

\multirow{2}{*}{XII} & \multirow{2}{*}{$F=q^{3/2}$} & &
\multicolumn{2}{|c|}{\multirow{2}{*}{$\inp{I}{8}=-\tfrac32$}} &
&  &  \\
 && & \multicolumn{2}{|c|}{} & & & \\ \hline

\end{tabular}
\end{sidewaystable}

\ifthenelse{\boolean{@twoside}}{%
\renewcommand{\chaptermark}[1]{\markboth{\headfont{#1}}{\headfont{#1}}}}{%
\renewcommand{\chaptermark}[1]{\markboth{\headfont{#1}}{}}}

\end{document}